\documentclass[12pt,leqno]{article}
\usepackage{amsmath,german,amssymb}
\usepackage[all]{xy}
\newcommand{\A}{\mathbb{A}}

\newcommand{\R}{\mathbb{R}}
\newcommand{\F}{\mathbb{F}}

\newcommand{\fC}{\mathbb{C}}

\newcommand{\Z}{\mathbb{Z}}

\newcommand{\oo}{\mathfrak{o}}
\newcommand{\pp}{\mathfrak{p}}

\newcommand{\sA}{{\cal A}}
\newcommand{\sB}{{\cal B}}

\newcommand{\sF}{{\cal F}}

\newcommand{\sR}{{\cal R}}

\newcommand{\sU}{{\cal U}}

\newcommand{\ve}{\varepsilon}
\newcommand{\vf}{\varphi}

\newcommand{\OO}{{\Omega}}

\newcommand{\wh}{\widehat}
\newcommand{\ol}{\overline}

\newcommand{\wt}{\widetilde}

\newcommand{\bea}{\begin{eqnarray}}
\newcommand{\eea}{\end{eqnarray}}

\begin{document}

\centerline{\large EXTENDIBLE FUNCTIONS AND LOCAL ROOT NUMBERS}
\centerline{\small Remarks on a paper of R.P.Langlands}
\centerline{by Helmut Koch and Ernst-Wilhelm Zink}

{\small {\bf Abstract:} This paper refers to Langlands' big set of notes [L] devoted to the question if the (normalized) local Hecke-Tate root number $\Delta=\Delta(E,\chi)$, where $E$ is a finite separable extension of a fixed non-archimedean local field $F$, and $\chi$ a quasicharacter of $E^\times$, can be appropriately extended to a local $\ve$-factor $\ve_\Delta=\ve_\Delta(E,\rho)$ for all virtual representations $\rho$ of the corresponding Weil group $W_E.$ Whereas Deligne [D] has given a relatively short proof by using the global Artin-Weil L-functions, the proof of Langlands is purely local and splits into two parts: the {\bf algebraic part} to find a minimal set of relations for the functions $\Delta$, such that the existence (and unicity) of $\ve_\Delta$ will follow from these relations; and the more extensive {\bf arithmetic part} to give a direct proof that all these relations are actually fulfilled. Our aim is to cover the algebraic part of Langlands' notes which can be done completely in the framework of representations of solvable profinite groups, where two modifications of Brauer's theorem play a prominent role.\\
Introduction / 1. The notion of extendible functions / 2. The kernel of the Brauer map and its generating relations / 3. A criterion for the extendibility of functions / 4. Recovering Theorem 3.1 for the case of local root numbers / Appendices A1. Proving Brauer 3 and Brauer 4 / A2. On type-III-groups.}\\

 \centerline{{\bf Introduction}}
 We begin from sketching the background of Langlands' paper [L] where for convenience we refer to [T3], where amongst others that background (for [L] and [D]) has been explained. Thus let $F$ be a global field and $\vf:W_F\rightarrow G_F$ the absolute Weil- and Galois group over $F$. Depending on an $F$-homomorphism $i_\nu:\ol F\rightarrow \ol F_\nu$ we have (for all places $\nu$ of $F$) the local-global relationship
 $$\xymatrix{
F_\nu^\times \ar[d] & W_{F_\nu} \ar[l] \ar[d]^{\theta_\nu} \ar[r] & G_{F_\nu}\ar[d]^{i_\nu^*} \\
C_F &  W_F \ar[l] \ar[r] & G_F  }$$
connecting the absolute Weil- and Galoisgroups resp., where $C_F=I_F/F^\times$ is the idele class group and the left square is given by class field theory. On $F_\nu^\times$ and $C_F$  we have the normalized absolute values $|\,|_{F_\nu}$ and $|\,|$ resp., which pull back to absolute values $\omega_\nu$ and $\omega$ on $W_{F_\nu}$ and $W_F$ resp. such that $\omega_\nu =\omega\circ \theta_\nu.$\\
 Let now $(\rho,V)$ be a representation of the global Weil group $W_F$, and let $V_\nu$ be the action of $W_{F_\nu}$ on $V$ which is obtained via $\theta_\nu.$ Then we have $(V\otimes\omega^s)_\nu = V_\nu\otimes \omega_\nu^s$ for all $s\in\fC$, and the Artin-Weil L-function
 $$  L(V,s):= L(V\otimes\omega^s)=\prod_\nu L((V\otimes\omega^s)_\nu),$$
 where the product is over all places $\nu$ of $F$ with local factors $L((V\otimes\omega^s)_\nu)$ as defined by E.Artin ([T3],(3.3)).
On the infinite product $L(V\otimes\omega^s)$ we know (loc.cit.(3.5.3)):\\
{\bf Theorem:} {\it The product converges for $s$ in some right half-plane and defines a function $L(V,s)$ which is meromorphic in the whole $s$-plane and satisfies the functional equation:
 $$ L(V,s) =\ve(V,s)\cdot L(V^\vee,1-s),\qquad \ve(V,s)\cdot\ve(V^\vee,1-s)=1,$$
 where $\ve(V,s):=\ve(V\otimes\omega^s),$ and $V^\vee$ is the contragredient of $V.$}\\
 For 1-dimensional $V$ corresponding to a {\bf quasicharacter $\chi$ of the idele-class-group $C_F$} this result was proved by Hecke and in modern version by Tate [T1]. More precisely
 $$ L(\chi,s) = \ve(\chi,s)\cdot L(\chi^{-1},1-s),\qquad  \ve(\chi,s) =\ve(\chi|\,|^s)= \prod_\nu\ve(\chi_\nu|\,|_{F_\nu}^s,\psi_\nu,dx_\nu)$$
 where now, due to loc.cit., we have a decomposition of $\ve(\chi,s)=\ve(\chi|\,|^s)$ as a product of local factors which is obtained as follows: Consider $\A_F$ the ring of $F$-adeles and $dx$ the Tamagawa measure on the additive group $\A_F.$ Then (see [W],VII,\S2) {\bf depending on} a fixed non-trivial additive character $\psi$ of $\A_F/F$ we obtain a decomposition
 $dx = \prod_\nu dx_\nu$ where $dx_\nu$ is the Haar measure on $F_\nu$ which is selfdual with respect to the local component $\psi_\nu:F_\nu^+\rightarrow \fC^\times$ of $\psi$. And for the pair $(\psi_\nu, dx_\nu)$ we have {\bf Tate's local functional equation [T3],(3.2.1)} which defines $\ve(\chi_\nu|\,|_{F_\nu}^s,\psi_\nu,dx_\nu)$ for the local components $\chi_\nu:F_\nu^\times\rightarrow\fC^\times$ of $\chi.$\\
The problem was, to establish a similar product decomposition
$$  \ve(V,s) =\ve(V\otimes\omega^s)= \prod_\nu \ve(V_\nu\otimes\omega_\nu^s,\psi_\nu,dx_\nu)$$
also for general (virtual) representations of $W_F.$ Since a local functional equation for higherdimensional $V_\nu$ is not available (Tate's approach rests on class field theory $W_{F_\nu}^{ab}\cong F_\nu^\times$) {\it the idea of {\bf [L]} was (purely locally) to extend} $\chi_\nu\mapsto \ve(\chi_\nu,\psi_\nu, dx_\nu)$  to a function $V_\nu\in R(W_{F_\nu}) \mapsto \ve(V_\nu,\psi_\nu,dx_\nu)$ on the {\bf Grothendieck group $R(W_{F_\nu})$} of virtual representations, which for field extensions $E_w|F_\nu$ is {\it inductive in dimension zero} (see below).\\ Then via {\bf [T3],(2.3.5)} $V\mapsto \prod_\nu\ve(V_\nu,\psi_\nu, dx_\nu)$ is inductive in dimension zero with respect to global extensions $E|F$, and is an extension of $\chi\mapsto \ve(\chi)=\prod_\nu\ve(\chi_\nu,\psi_\nu,dx_\nu).$  Because $V\mapsto \ve(V)$ has this property too (it is actually inductive in the unconditional sense), the two extensions of $\chi\mapsto \ve(\chi)$ must agree.\\
{\bf From now we restrict to the local problem with respect to the nonarchimedean completions $F_\nu$}, that means $F$ is now a nonarchimedean local field, $\psi$ a nontrival character of $F^+$ and $\chi$ a quasicharacter of $F^\times.$  And corresponding to the above decomposition of the Tamagawa measure, [L] restricts to functions $\ve(\chi,\psi):= \ve(\chi,\psi,dx)$ where $dx$ is the $\psi$-selfdual Haar measure on $F^+$. Then Tate's local functional equation for the 1-dimensional case ({\bf loc.cit.}) yields
$$ \ve(\chi,\psi)\cdot\ve(\chi^{-1}|\,|_F,\psi)=\chi(-1),\qquad \Delta(\chi,\psi)\Delta(\chi^{-1},\psi)=\chi(-1),$$
{\it if we take (as {\bf [L]} does) the normalization} $\Delta(\chi,\psi):=\ve(\chi|\,|_F^{1/2},\psi).$ By class field theory we have $W_E^{ab}\cong E^\times$ for all finite extensions $E|F$ and [L] considers the function\\
{\bf (*)}
$$  (E,\chi_E) \mapsto \Delta(\chi_E,\psi_{E|F})\in \fC^\times,\quad\textrm{where}\;\psi_{E|F}:=\psi_F\circ Tr_{E|F},$$
for all finite separable extensions $E|F$ and quasicharacters $\chi_E$ of $E^\times.$\\ {\bf Problem:} Is it possible to extend this to a function
$$  (E,\rho_E)\mapsto \ve_\Delta(\rho_E,\psi_{E|F})$$
which is defined for all $E|F$ and virtual representations $\rho_E\in R(W_E)$, such that

$(E,\rho_1 \oplus\rho_2)$ maps to
$\ve_\Delta(\rho_1\oplus\rho_2, \psi_{E|F})=
\ve_\Delta(\rho_1,\psi_{E|F})\ve_\Delta(\rho_2,\psi_{E|F})$, and moreover
$$  \ve_\Delta(Ind_{E'|E} (\rho_{E'}),\psi_{E|F}) = \ve_\Delta(\rho_{E'},\psi_{E'|F})$$
if $\rho_{E'}\in R(W_{E'})$ is virtual {\bf of dimension zero}. The answer is yes, and going backwards we take $\ve(V_\nu\otimes\omega_\nu^s,\psi_\nu,dx_\nu)= \ve_\Delta(V_\nu\otimes \omega_\nu^{s-\frac{1}{2}},\psi_\nu)$ to obtain the decomposition of the global $\ve(V,s).$
The paper {\bf [L]} splits into two parts:
\begin{itemize}
 \item {\bf Algebraic part:} Find a minimal set of relations for the functions $\Delta=\Delta(\chi_E,\psi_{E|F})$ such that the existence of a (uniquely determined) extension $\ve_\Delta=\ve_\Delta(\rho_E,\psi_{E|F})$ will follow from these relations.
 \item {\bf Arithmetic part:} Give a direct proof that all these relations are actually fulfilled.
\end{itemize}
{\bf In this paper we will cover only the algebraic part}, but from a purely group theoretical point of view which presents the arguments in a more conceptual way. Instead of local Weil groups $W_F$ we will consider solvable profinite groups $\OO$ as for instance the absolute Galois groups $G_F.$ This comes close because we have the continuous embedding $W_F\subset G_F$ with dense image.\\
Thus let $\OO$ be a solvable profinite group. By a {\it subgroup $H$ of $\OO$, denoted $H\le\OO,$ we usually understand an open subgroup}, or, equivalently, a closed subgroup of finite index. Besides that we will consider the {\it closed commutator subgroups} $[H,H]\le\OO$ which in general will not be open anymore. By definition we have $\OO=\varprojlim_N \OO/N$ as a projective limit over the open normal subgroups. Together with $N$ also $[N,N]$ is normal in $\OO$ and therefore also\\
{\bf (1)}
$$\OO = \varprojlim_N \OO/_{[N,N]},\qquad R(\OO) =\varinjlim_N R(\OO/_{[N,N]}),$$
where $R(G)$ denotes the {\bf Grothendieck group} of virtual {\bf continuous} representations of $G$  and $R(\OO/_{[N,N]})\hookrightarrow R(\OO)$ is naturally embedded.

Instead of pairs $(E,\chi_E)$ as above we consider now:\\
$R_1(\le\OO)\ni (H,\chi_H)$ {\it the set of pairs such that $H\le\OO$ and $\chi_H:H\rightarrow\fC^\times$ is a continuous 1-dimensional character}, for short we will write: $\chi_H\in H^*.$\\ From that point of view {\bf the question is}, when a function $\Delta =\Delta(H,\chi)$ on $R_1(\le\OO)$ can be extended to a function $\sF=\sF(H,\rho)$ on $R(\le\OO):=\bigsqcup_{H\le\OO} R(H)$ such that
$$ \sF(H,\rho_1+\rho_2) = \sF(H,\rho_1)\sF(H,\rho_2),\qquad \sF(H',Ind_H^{H'}(\rho)) =\sF(H,\rho),$$
if $H\le H'$ and $\rho\in R(H)$ is {\bf of dimension zero}.\\
Our main tool will be Brauer's Theorem and it's variations: The induction map $(H,\chi_H)\mapsto Ind_H^\OO(\chi_H)\in  R(\OO)$ is actually a map which depends only on $[H,\chi_H]:= \OO$-conjugacy class of $(H,\chi_H)$ and induces the {\bf Brauer map}\\
{\bf (2)}
$$  \vf_\OO: R_+(\le\OO)\rightarrow R(\OO),\qquad \sum n_{[H,\chi_H]}[H,\chi_H] \mapsto \sum n_{[H,\chi_H]}Ind_H^\OO(\chi_H),$$
where $R_+(\le\OO)$ {\it is the free abelian group, generated by the $\OO$-conjugation classes} $[H,\chi_H].$
Then we have: (see {\bf [S], section 10.5} as an easy accessible reference)

{\bf Brauer's Theorem 1:} {\it Let $\OO$ be a profinite group and $\rho$ a virtual representation of $\OO.$ Then there are pairs $(H_i,\chi_i)$ consisting of open subgroups $H_i$ and 1-dimensional characters $\chi_i\in H_i^*$, and integers $n_i$ such that
$$  \rho \cong \sum_i n_i Ind_{H_i}^\OO(\chi_i)=\vf_\OO\left(\sum_i n_i[H_i,\chi_i]\right).$$
Consequently $\vf_\OO:R_+(\le\OO)\rightarrow R(\OO)$ is a surjective homomorphism of additive groups, and actually ({\bf [D],1.9}, see below) it is a homomorphism of commutative rings.}\\
{\bf Brauer's Theorem 2:} (see {\bf [S],exercise 10.6.} or {\bf [D],Prop.1.5}) {\it If $\rho\in R(\OO)$ is virtual of dimension 0 then it admits a Brauer presentation:
$$  \rho = \sum_i n_i Ind_{H_i}^\OO(\chi_{H_i} - 1_{H_i})= \vf_\OO\left(\sum_i n_i([H_i,\chi_{H_i}] - [H_i, 1_{H_i}])\right).$$}\\
If $\OO$ is finite, the groups $H_i$ can be taken elementary, in particular nilpotent, subgroups of $\OO$. For profinite groups this fails but anything else is maintained because $R(\OO) = \varinjlim_N R(\OO/N)$ for open normal subgroups $N$, yields a reduction to the finite case.

In order to simplify the {\bf chapters 15-19 of [L]} and to present them in a more conceptual way by going down to the purely group theoretical background, the manuscript {\bf [K]} had been built on these two theorems and on the description of $Ker(\vf_\OO)$ as in {\bf [D],\S1}. But at a certain point this approach failed, because it didn't take care of the role of relative Weil-groups $W_{E|F}=W_F/[W_E,W_E]$ in the original paper {\bf [L]}. From the purely group theoretical point of view using relative Weil groups means to approach $R(\OO)$ via {\bf (1)} and to use the following two more modifications of Brauer's theorem:\\
If $N\le\OO$ is an open normal subgroup, then we will write $R_1(N\le\OO)$ {\bf for the set of pairs $(H,\chi_H)$ such that $H\ge N,$} hence $\chi_H$ is actually a character of $H/_{[N,N]}$ and $Ind_H^\OO(\chi_H) \in R(\OO/_{[N,N]})$ because the commutator subgroup $[N,N]$ is normal in $\OO.$

{\bf Brauer's Theorem 3:} {\it Let $\OO$ be a profinite group, and $N\le\OO$ an open normal subgroup. Consider $R_+(N\le\OO)\subseteq R_+(\le\OO)$ the $\Z$-submodule which is generated by all conjugacy classes $[H,\chi_H]$ such that $H\ge N.$ Then the Brauer map $\vf_\OO$ induces a surjection\\
{\bf (3)}
$$  \vf_{N\le\OO}:R_+(N\le\OO)\twoheadrightarrow R(\OO/_{[N,N]}) \;\subset\; R(\OO)$$
such that any virtual representation $\rho$ of $\OO/_{[N,N]}$ can be presented as integral combination of monomial representations
$$  \rho \cong \sum_i n_i Ind_{H_i}^\OO(\chi_i),\quad
\textrm{where all $H_i$ are open subgroups containing $N$.}$$}
{\bf Brauer's Theorem 4:} {\it And if $\rho\in R(\OO/_{[N,N]})$ is virtual of dimension 0 then it admits a Brauer presentation:
$$  \rho = \sum_i n_i Ind_{H_i}^\OO(\chi_{H_i} - 1_{H_i}) = \vf_{N\le\OO}\left(\sum_i n_i([H_i,\chi_{H_i}] - [H_i, 1_{H_i}])\right),$$
where all occurring subgroups $H_i$ contain $N.$}

{\bf Remark:} Since $\vf_{N\le \OO}$ is nothing else than the restriction of $\vf_\OO,$ it may happen that we sometimes omit the $N.$ Similar as for Brauer 1 and 2, the case of profinite groups can be reduced to the case of finite ones. Some further remarks we will add in {\bf section 4} where $\OO/_{[N,N]}$ will be replaced by a relative Weil group $W_{K|F}= W_F/_{[W_K,W_K]}.$

The easiest example is the case $N=\OO$ where $\vf_{\OO\le\OO}:R_+(\OO\le\OO)\rightarrow R(\OO/_{[\OO,\OO]})$ is the identity because both sides represent the free $\Z$-module over $\OO^*.$  Proofs of Brauer 3 and 4 will be given in {\bf Appendix 1}.\\

Now we return to the question when a function $\Delta =\Delta(H,\chi)$ on $R_1(\le\OO)$ can be extended to a function $\sF=\sF(H,\rho)$ on $R(\le\OO):=\bigsqcup_{H\le\OO} R(H)$ as above. \\
Langlands' approach in {\bf [L]}, as we interpret it here, is, for all finite subquotients $\OO'/N'$ of $\OO$ to consider the restrictions $\Delta'$ of $\Delta$ to $R_1(N'\le\OO') =\bigsqcup_{H;\;N'\le H\le\OO'} H^*$ and to prove in an inductive way that $\Delta'$ can be extended to $\sF'$ on $R(N'\le\OO'):= \bigsqcup_{H;\;N'\le H\le\OO'} R(H/_{[N',N']}).$ Thus fixing a natural number $n$ we assume that $\Delta$ extends from $R_1(N'\le\OO')$ onto $R(N'\le\OO')$ for all subquotients such that $(\OO':N')< n,$ and ask for conditions to conclude that $\Delta$ will then extend also for the cases where $(\OO':N')=n.$ (Of course it can happen that the assumption $(\OO':N')=n$ is empty, then we may go directly to $n+1.$) And since the arguments will not depend on taking $\OO$ or $\OO'$ as our absolute group, we come down to the problem of extending $\Delta$ from $R_1(N\le\OO)$ onto $R(N\le\OO),$ if we assume this for all proper subquotients $\OO'/N'$ of $\OO/N.$  As a general criterion which is motivated by {\bf [L], Theorem 2.1} we have:

{\bf 1.4 Proposition:} {\it Let $(\OO,N)$ be a profinite group and an open normal subgroup, and $\Delta$ a function on $R_1(N\le\OO)$ with values in the multiplicative abelian group $\sA.$ Then $\Delta$ is extendible to $\sF$ on $R(N\le\OO)=\bigsqcup_{\OO';\;N\le \OO'\le\OO} R(\OO'/_{[N,N]})$ {\bf if and only if} for all subgroups $H\le \OO$ containing $N$  there is a function
$$  H\in \sU_N(\OO) \mapsto \lambda_H^\OO(\Delta) \in \sA,$$
such that:  {\bf (c1)}  $.\quad \lambda_\OO^\OO(\Delta) =1,$ and\\
{\bf (c2)}   Any relation of the form  $\sum_{i=1}^r n_i\cdot Ind_{H_i}^\OO(\chi_i)\cong 0\;\in R(\OO/_{[N,N]})$ (with subgroups $H_i\ge N$ and multiplicities $n_i\in\Z$)  will imply the relation
 $$   \prod_{i=1}^r \Delta(H_i,\chi_i)^{n_i}\lambda_{H_i}^\OO(\Delta)^{n_i} =1\quad\in\sA.$$
According to Brauer 3 and 4 the extension
$$\rho'= \sum_{i=1}^r n_i\cdot Ind_{H_i}^{\OO'}(\chi_i)
\in R(\OO'/_{[N,N]})\mapsto \sF(\OO',\rho')=\prod_{i=1}^r \Delta(H_i,\chi_i)^{n_i}\lambda_{H_i}^{\OO'}(\Delta)^{n_i}\;\in\sA$$ where
$\lambda_{H_i}^{\OO'}(\Delta):= \lambda_{H_i}^\OO(\Delta)\cdot \lambda_{\OO'}^\OO(\Delta)^{-(\OO':H_i)},$
is then well defined and uniquely determined.}\\

Now in order to check {\bf 1.4 (c2)} we {\bf need generators for the kernel of the restricted Brauer map $\vf_{N\le\OO}$} as in {\bf (3)}. Specifying such a set of generators for the case that $\OO$ is (pro){\bf solvable} is the main result of {\bf section 2.} Actually this is a variation of {\bf [D],\S1} where the case $\OO$ finite and $N=\{1\}$ has been dealt with.\\
Then in {\bf section 3} we proceed in three steps:\\
{\bf Step 1:} Definition of $H\in\sU_N(\OO) \mapsto \lambda_H^\OO(\Delta):=\lambda_{H/N}^{\OO/N}(\Delta)$ where we will work in the finite quotient $G=\OO/N.$\\
{\bf Step 2:} Assuming  that $\Delta$ extends from $R_1(N'\le\OO')$ onto $R(N'\le\OO')$ for all proper subquotients $G'=\OO'/N'$ of $G=\OO/N$, we establish the tower relation
$$  \lambda_{H'}^{\OO'}(\Delta)= \lambda_{H'}^H(\Delta)\cdot\lambda_H^{\OO'}(\Delta)^{(H:H')}$$
for subgroups $N\le H'\le H\le\OO'\le\OO.$\\
{\bf Step 3:} Specifying {\bf three types} of special relations {\bf (c2)} which have still to be verified in order to ensure that $\Delta$ extends from $R_1(N\le\OO)$ onto $R(N\le\OO)$ if the extension for proper subquotients $N'\le\OO'$ is already known to exist. More precisely these three types of relations occurring in the {\bf Main Theorem 3.1} are:
\begin{itemize}
 \item[{\bf (1)}] a Davenport-Hasse relation relative to subquotients $B/K$ of $\OO$ which are cyclic of prime order $\ell,$
 \item[{\bf (2)}] a Heisenberg identity relative to subquotients $B/[Z,B]$ which are 2-step nilpotent and such that $B/Z\cong \Z/\ell\times\Z/\ell$ is bicyclic for some prime $\ell,$ and $[B,B]/[Z,B]\stackrel{\sim}{\leftarrow} B/Z \wedge B/Z$ is cyclic of order $\ell,$
 \item[{\bf (3)}] a generalized Davenport-Hasse relation relative to subquotients\\ $B/K =H/K \ltimes C/K$ where $H< B$ is a subgroup which is {\bf maximal} but {\bf not normal}, and $K=\bigcap_{b\in B} bHb^{-1}$ is the corresponding normal subgroup. Actually this means that $B/K$ is a type-III-group as it is explained in {\bf Appendix 2.}
\end{itemize}

{\bf The idea to produce special relations} is to consider the case $N=C$ of a commutative normal subgroup, hence $\vf_{C\le\OO}:R_+(C\le\OO)\twoheadrightarrow R(\OO),$ and to establish a projector
$$ \Phi_C:R_+(\le\OO)\twoheadrightarrow R_+(C\le\OO) \hookrightarrow R_+(\le\OO),\qquad \Phi_C\circ\Phi_C =\Phi_C,$$
such that $\vf_{C\le\OO}\circ\Phi_C = \vf_\OO.$ This yields $x-\Phi_C(x)\in Ker(\vf_\OO)$ for all $x\in R_+(\le\OO).$  The projector $\Phi_C$ is obtained as a corollary to {\bf Lemma 2.2} which essentially goes back to {\bf [L],15.1} and {\bf [D],1.11} resp.\\
Concerning steps 1 and 2 we will basically follow {\bf [L],chapt.16}, and then our step 3 will cover {\bf chap.19}.\\
Altogether this covers the purely algebraic part of {\bf [L]}. We will finish with a short introduction to the arithmetic part, where $\OO=G_F$ is the absolute Galois group of a $p$-adic field $F$ (or the corresponding Weil group $W_F\subset G_F$) and $\Delta(H,\chi) =\Delta(E,\chi_E)$ is (as in {\bf (*)} above) the local Hecke-Tate root number for $E|F$ a finite extension and $\chi_E:E^\times\rightarrow\fC^\times$ a quasicharacter. Then the relations {\bf 3.(1) - (3)} turn into identities {\bf 4.(1) - (3)} for these local root numbers which are sufficient to ensure the existence of local $\ve$-factors for higherdimensional representations. Verifying these arithmetic relations is the other main topic of {\bf [L]} which is not touched here. In {\bf section 4} we do nothing else than translating the algebraic relations {\bf 3.(1) - (3)} into the arithmetic relations of {\bf [L]} and adding some first simplifications which occur in the arithmetic context.\\

{\bf 1. The notion of extendible functions}\\
The following approach is a generalization which includes the corresponding considerations in {\bf [T2]} and in {\bf [L],chap.2} as well.\\

Let $(\OO,N)$ be a pair consisting of a profinite group $\OO$ and an open normal subgroup $N$ (equivalently: $N$ is closed and of finite index).
We denote by $$R(N\le\OO):= \bigsqcup_{H;\;N\le H\le\OO} R(H/_{[N,N]})$$ the set of all pairs $(H,\rho)$ where:
\begin{itemize}
 \item $H\le\OO$ is a subgroup {\bf containing $N$},
 \item $\rho\in R(H/_{[N,N]})$ is the equivalence class of a virtual representation of $H/_{[N,N]}.$\\
\end{itemize}
Since we consider equivalence classes, the group $\OO/N$ acts on $R(N\le\OO)$ by means of
\begin{itemize}
 \item  $(H,\rho)^g = (H^g,\rho^g),\qquad g\in\OO/N,\; H^g:= g^{-1}Hg,$
\item $\rho^g(x) = \rho(gxg^{-1}) =:(g^{-1}\rho g)(x),\qquad\textrm{for}\; x=g^{-1}yg\in H^g.$
\end{itemize}
We will also consider {\bf the case $N= \{1\}$ where we will write: $R(\le\OO):=\bigsqcup_{H;H\le\OO} R(H).$}\\
Note also that $R(\OO\le\OO) = R(\OO/_{[\OO,\OO]})$ are the virtual representations of the abelian quotient $\OO/_{[\OO,\OO]}.$\\
Above we have considered already
$$R_1(N\le\OO)=\bigsqcup_{H;\;N\le H\le\OO} R_1(H/_{[N,N]}) \quad\subset R(N\le\OO)$$
the subset of pairs $(H,\chi)$ with $H\ge N$ and $\chi\in H^*,$ hence $\chi$ is trivial on $[N,N].$  
  {\bf Therefore by the very definition we have:}
$$  R_1(N\le\OO) = R_1(\ol N\le \ol \OO),\qquad R(N\le\OO) = R(\ol N\le\ol\OO),$$
where $\ol N:= N/_{[N,N]}$ is an abelian normal subgroup of $\ol\OO:=\OO/_{[N,N]}.$ More general we will write $\ol H=H/_{[N,N]}$ for all $H\ge N.$\\
Also note that for $N\le N'\le \OO'\le \OO$ where $N'$ is a normal subgroup of $\OO',$ we have natural embeddings
\begin{itemize}
 \item[{\bf (e1)}] $\quad R_1(N'\le \OO')\subset R_1(N\le \OO') \subset R_1(N\le \OO),$
 \item[{\bf (e)}] $\quad R(N'\le \OO')\subset R(N\le \OO') \subset R(N\le \OO),$
\end{itemize}
where in each case the first embedding is by inflation contravariant to $H/_{[N,N]}\twoheadrightarrow H/_{[N',N']},$ and the second embedding considers $R(N\le \OO')$ as the subset of pairs $(H,\rho)$ such that $H\le \OO'.$ \\

And the direct limit for $N\rightarrow \{1\}$ yields
$$  \varinjlim_{N\rightarrow \{1\}} R(N\le\OO) = R(\le\OO),\qquad\varinjlim_{N\rightarrow \{1\}} R_1(N\le\OO) = R_1(\le\OO).  $$
In {\bf [T2]} and also in {\bf [K]} only the absolute versions $R(\le\OO),$ $R_1(\le\OO)$ and the limits with respect to $\OO=\varprojlim_N\OO/N$ have been dealt with. But for the aims of this paper it is necessary two follow {\bf [L]}, who works with relative Weil-groups, more closely. From the group theoretical point of view, working with relative Weil groups means to consider the limit $\OO=\varprojlim_N \OO/_{[N,N]}.$\\

In the following we will focus to the relative case where a lower bound $N\le \OO$ will be fixed. But obviously all considerations take over to the case $N=\{1\}.$\\
Suppose we have a function $\Delta$ defined on $R_1(N\le\OO)$ taking values in some multiplicative abelian group $\sA$ with properties\\
{\bf (1)}
$$   \Delta(H,1_H) =1,$$
{\bf (2)}
$$  \Delta(H^g, \chi^g) =\Delta(H,\chi)$$
for all $(H,\chi)\in R_1(N\le\OO),$ where $1_H$ denotes the trivial representation of $H.$\\

{\bf\large (!)} From now, when speaking of {\it functions $\Delta$ on $R_1(N\le\Omega)$} we will always mean functions with the properties {\bf (1)} and {\bf (2)}.\\

{\bf 1.1 Definition:} We will say that the function  $(H,\chi)\mapsto \Delta(H,\chi)$ on $R_1(N\le\OO)$ is {\bf extendible}, if it can be extended to an $\sA$-valued function $(H,\rho)\mapsto\sF(H,\rho)$ on $R(N\le\OO)$ satisfying:
$$ \sF(H,\chi)= \Delta(H,\chi),\quad\textrm{if}\; (H,\chi)\in R_1(N\le\OO),$$
\begin{itemize}
\item[{\bf (3)}]
$ \sF(H,\rho_1 +\rho_2) = \sF(H,\rho_1)\cdot\sF(H,\rho_2),\qquad\forall \; (H,\rho_i)\in R(N\le\OO),$ that means $\sF$ restricts to a homomorphism of groups $\sF:R(\ol H)\rightarrow\sA$ for each $H\ge N,$
\item[{\bf (4)}]  $\sF(G,Ind_H^G(\rho))=\sF(H,\rho),$ if $N\le H\le G\le \OO$ and $\rho\in R(\ol H)$ such that $\bf dim(\rho)=0.$ More general this implies (5) below, if $dim(\rho)\ne 0.$
\end{itemize}

{\bf 1.2 Uniqueness:} {\it If $\Delta$ is extendible, then the extension $\sF$ satisfying {\bf (3)} and {\bf (4)} will be uniquely determined.}\\

{\bf Proof:} Let be $(G,\rho)\in R(N\le\OO).$ Actually $\rho$ is a virtual representation of $\ol G=G/_{[N,N]}$ and therefore we may apply {\bf Brauer 4} which gives us:
$$  \rho -dim(\rho)1_G = \sum _i n_i Ind_{\ol H_i}^{\ol G}(\chi_{H_i} - 1_{H_i}),$$
for groups $H_i\ge N.$ And
using {\bf (1)} and the property {\bf (4)} this implies:
$$  \sF(G,\rho) = \sF(G,\rho - dim(\rho)\cdot 1_G) = \prod_i \sF(Ind_{H_i}^G(\chi_i -1_{H_i}))^{n_i} =$$
$$ =\prod_i \sF(H_i,\chi_i-1_{H_i})^{n_i} = \prod_i\Delta(H_i,\chi_i)^{n_i},$$
which computes $\sF(G,\rho)$ in terms of the original $\Delta$ on $R_1(N\le\OO)$. Uniqueness follows because for any other extension $\sF'$ we could do the same computation. {\bf qed.}\\

{\bf 1.3 $\lambda$-Function:} {\it Let $\Delta$ be a function on $R_1(N\le\OO)$ which extends to $\sF$ on $R(N\le\OO).$ Then, in the situation {\bf (4)} but for $(H,\rho)\in R(N \le\OO)$ of arbitrary dimension, we will have:\\
{\bf (5)}
$$ \sF(G, Ind_H^G(\rho)) =\lambda_H^G(\Delta)^{dim(\rho)}\cdot \sF(H,\rho),\qquad \sF(G, Ind_H^G(\chi)) =\lambda_H^G(\Delta)\cdot \Delta(H,\chi),$$
where by {\bf definition:} $\lambda_H^G(\Delta):=  \sF(G, Ind_H^G(1_H))$ is independent of $\rho.$}\\

Indeed, we consider the virtual representation $\rho -dim(\rho)\cdot 1_H \in R(H/_{[N,N]})$ and then
due to {\bf (4)} we have
$$  \sF\left(G, Ind_H^G(\rho - dim(\rho)\cdot 1_H)\right) = \sF(H,\rho- dim(\rho)\cdot 1_H).$$
Using $Ind_H^G(\rho - dim(\rho)\cdot 1_H)= Ind_H^G(\rho) - dim(\rho)\cdot Ind_H^G(1_H)$ and
the property {\bf (3)}, this rewrites as
$$  \frac{\sF(G, Ind_H^G(\rho))}{\sF(G,Ind_H^G(1_H))^{dim(\rho)}} = \frac{\sF(H,\rho)}{\Delta(H,1_H)^{dim(\rho}},$$
such that {\bf (5)} follows from {\bf (1)}.\\
The function $H\le G \mapsto \lambda_H^G(\Delta)$, for {\bf subgroups $H\le G$ between $N$ and $\OO$}, is Langlands' $\lambda$-function. (We need not to refer here to the extension $\sF$ because it is unique.)\\

We mention {\bf four properties of the $\lambda$-functions} which are attached to extendible functions $\Delta$ and which easily follow from the definitions.\\
{\bf (6)}
$$  \lambda_{H^g}^G(\Delta) =\lambda_H^G(\Delta),\quad \textrm{if}\;g\in G,$$
{\bf (7)}
$$ \lambda_{H'}^G(\Delta) =\lambda_{H'}^H(\Delta)\cdot\lambda_H^G(\Delta)^{[H:H']},\quad\textrm{if}\; G\ge H\ge H',$$
which follows by transitivity of induction.\\
{\bf (8)}
$$  \lambda_H^G(\Delta) = \lambda_{H/H'}^{G/H'}(\Delta_{G/H'}),$$
if $H'\le H\le G$ such that $H'$ is normal in $G$, hence $Ind_H^G(1_H)=Ind_{H/H'}^{G/H'}(1_{H/H'}).$ (For subquotients  $G/H'$ of $\OO/N$ the functions $\Delta_{G/H'}$ and $\sF_{G/H'}$ are defined via the inflation map; see the remark behind {\bf 1.4} and {\bf (F3)} of section 2.)
In particular this applies if $H'=H$ is normal in $G.$  
Therefore {\bf (8)} follows because for $N\le N'\le\OO'\le\OO$ the extendible $\Delta$ will in particular extend from $R_1(N'\le \OO')$ onto $R(N'\le \OO'),$ and from
$R_1(\le \OO'/N')$ onto $R(\le \OO'/N')$, and the extensions will be uniquely determined. Moreover we have:\\
{\bf (9)}
$$  \lambda_H^G(\Delta) = \prod_{\chi\in(G/H)^*} \Delta(G,\chi),$$
if $G/H$ is an abelian subquotient of $\OO/N,$ because $Ind_H^G(1) = \sum_{\chi\in(G/H)^*} \chi$ and therefore {\bf (3)} applies.\\

Actually the existence of a $\lambda$-function $\lambda_H^G(\Delta)$ with appropriate behaviour can be turned into a criterion for the extendibility of $\Delta$. Compared to {\bf [K]}, Lemma 3.2, we deal here with the refined problem of extending $\Delta$ from $R_1(N\le\OO)$ to $R(N\le\OO)$ which is motivated by {\bf [L],Theorem 2.1}.\\

{\bf\large 1.4 Proposition (Criterion):  (i)} {\it Let $(\OO,N)$ be a profinite group and an open normal subgroup, and $\Delta$ a function on $R_1(N\le\OO)$ with values in the multiplicative abelian group $\sA.$ Then $\Delta$ is extendible to $\sF$ on $R(N\le\OO)=\bigsqcup_{H;\;N\le H\le\OO} R(H/_{[N,N]})$ {\bf if and only if} for all subgroups $H\le \OO$ containing $N$  there is a function
$$  U\in \sU_N(H) \mapsto \lambda_U^H(\Delta)=\lambda_{U/N}^{H/N}(\Delta)\in \sA,$$
which is defined on the set $\sU_N(H)$ of subgroups of $H$ which contain $N,$ such that:
\begin{itemize}
 \item[{\bf (c1)}]   $.\quad \lambda_H^H(\Delta) =1$
 \item[{\bf (c2)}]   Any relation of the form  $\sum_{i=1}^r n_i\cdot Ind_{U_i}^H(\chi_i)\cong 0$  will imply the relation
 $$   \prod_{i=1}^r \Delta(U_i,\chi_i)^{n_i}\lambda_{U_i}^H(\Delta)^{n_i} =1.$$
\end{itemize}
Moreover the {\bf functions} $U\in \sU_N(H)\mapsto \lambda_U^H(\Delta)$ are uniquely determined by these requirements and satisfy the tower relation:
$$ \lambda_U^G(\Delta) = \lambda_U^H(\Delta)\cdot\lambda_H^G(\Delta)^{(H:U)},\quad\textrm{for subgroups}\;N\le U\le H\le G\le \OO.$$
{\bf\large (ii)}: For $\Delta$ on $R_1(N\le\OO)$ being extendible to $\sF$ on $R(N\le\OO)$, it is already sufficient to have a function
$$  U\in \sU_N(\OO) \mapsto \lambda_U^\OO(\Delta)=\lambda_{U/N}^{\OO/N}(\Delta)\in \sA,$$
such that {\bf (c1), (c2)} are fulfilled in the particular case $H=\OO,$ (and then using $\lambda_U^H(\Delta):= \lambda_U^\OO(\Delta)\lambda_H^\OO(\Delta)^{-(H:U)}$).}\\

{\bf Remark:} Consider $R_1(\le\OO/N)=\{(H,\chi)\;|\; \chi\in (H/N)^*\}\subset R_1(N\le\OO)$ and
\begin{itemize}
 \item let $\Delta_0$ be the restriction of $\Delta$ down to $R_1(\le\OO/N).$
\end{itemize}
If $\Delta$ extends to $\sF$ on $R(N\le\OO) =\bigsqcup_{H\ge N} R(\ol H),$ then the restriction $\sF_0$ of $\sF$ down to $R(\le \OO/N)=\bigsqcup_{H\ge N} R(H/N)$ will be an extension of $\Delta_0$ and therefore $Ind_H^G(1_H)=Ind_{H/N}^{G/N}(1_{H/N})$ will imply:
$ \lambda_H^G(\Delta) = \lambda_{H/N}^{G/N}(\Delta_0),\quad\forall\; N\le H\le G\le\OO.$\\   A priori Langlands distinguishes $\lambda$-functions $U\in \sU_N(\OO)\mapsto \lambda_U^\OO(\Delta),$ from {\bf weak} $\lambda$-functions:
$U/N\in \sU(\OO/N) \mapsto \lambda^{\OO/N}_{U/N}(\Delta),$ where the relations {\bf (c2)} are requested only for characters $\chi$ of $U/N.$\\
{\bf Proof of the Criterion 1.4 (i):}  If $\Delta$ is extendible to $\sF,$ then, as we see from {\bf (5)}, the definition $\lambda_U^H(\Delta) := \sF(H, Ind_U^H(1))$ will have the requested properties.\\
{\bf Conversely} suppose that for all $H\le \OO$ containing $N$ the functions $\lambda_U^H(\Delta)$ together with the displayed properties {\bf (c1), (c2), do exist}. Let now $\rho\in R(\ol H)$ be a virtual representation. Then by {\bf Brauer 3} there are subgroups $U_i\ge N$ of $H$, characters $\chi_i\in U_i^*$ and multiplicities $n_i\in\Z$ such that:\\
{\bf (10)}
$$   \rho\cong \sum_i n_i\cdot Ind_{U_i}^H(\chi_i).$$
Then we define:\\
{\bf (11)}
$$ (H,\rho)\in R(N\le\OO) \mapsto  \sF(H,\rho):= \prod_i \Delta(U_i,\chi_i)^{n_i}\cdot\lambda_{U_i}^H(\Delta)^{n_i},$$
{\bf in particular:} $\quad \sF(H,Ind_U^H(\chi)):=\Delta(U,\chi)\cdot\lambda_U^H(\Delta),$  $\quad \sF(H,Ind_U^H(1)):=\lambda_U^H(\Delta).$\\
Our assumptions on $U\in \sU_N(H) \mapsto \lambda_U^H(\Delta)\in \sA$ guarantee that $\sF(H,\rho)$ is independent  of the choice of the presentation {\bf (10)} of $\rho.$ And obviously the relation {\bf (3)} will hold.  It remains to show {\bf (5)}.\\
Let $G$ be a subgroup of $\OO$ containing $H$ and consider $(H,\rho)\in R(N\le\OO)$ as before. Then {\bf (10)} implies:
$$ dim(\rho)=\sum_i n_i\cdot(H:U_i),\qquad Ind_H^G(\rho) \cong \sum_i n_i\cdot Ind_{U_i}^G(\chi_i).$$
Hence by definition
$$   \sF(G,Ind_H^G(\rho)) =\prod_i \Delta(U_i,\chi_i)^{n_i}\cdot \lambda_{U_i}^G(\Delta)^{n_i}.$$
Now dividing this by {\bf (11)} we obtain
$$  \frac{\sF(G,Ind_H^G(\rho))}{\sF(H,\rho)} =\prod_i \left(\frac{\lambda_{U_i}^G(\Delta)}{\lambda_{U_i}^H(\Delta)}\right)^{n_i},$$
thus we see that in order to prove {\bf (5)} it is sufficient to verify the property {\bf (7)}:
$$   \lambda_{U_i}^G(\Delta) = \lambda_{U_i}^H(\Delta)\cdot \lambda_H^G(\Delta)^{[H:U_i]}.$$
To establish that property we remark:\\
{\bf Remark 1.} For a fixed $H\le G$ the function $U\in\sU_N(H) \mapsto \lambda_U^H(\Delta)$ is by the requirements {\bf (c1), (c2)} {\bf uniquely determined.} Indeed, by {\bf Brauer 2} with respect to the group $H/N$, we have a relation:
$$  Ind_U^H(1_U) - [H:U]\cdot 1_H \cong \sum_i n_i\cdot Ind_{U_i}^H(\chi_i - 1_{U_i})=\sum_i n_iInd_{U_i}^H(\chi_i) -\sum_i n_iInd_{U_i}^H(1_{U_i}),$$
for certain $N\le U_i< H,$ $\chi_i\in (U_i/N)^*, n_i\in \Z.$  Now using {\bf (c2)} this implies the relation:
$$ \frac{\Delta(U,1_U)\cdot\lambda_U^H(\Delta)}{\Delta(H,1_H)^{[H:U]}} =\frac{\prod_i\left(\Delta(U_i,\chi_i)\lambda_{U_i}^H(\Delta)\right)^{n_i}}{\prod_i\left(\Delta(U_i,1_{U_i})\lambda_{U_i}^H(\Delta)\right)^{n_i}},$$
and using the property {\bf (1)} of $\Delta$ this comes down to:
$$   \lambda_U^H(\Delta) = \prod_i \Delta(U_i,\chi_i)^{n_i}.$$
Thus the $\lambda$-function can be computed in terms of the original $\Delta$, hence it will be uniquely determined if it exists.\\

{\bf Remark 2.} For $N\le H\le G\le \OO$ consider the function
$$  U\in \sU_N(H) \mapsto \lambda_U^H:=\lambda_U^G(\Delta)\lambda_H^G(\Delta)^{-[H:U]}.$$
This will be a function satisfying (c1), (c2): because for any relation $\sum_i n_i Ind_{U_i}^H(\chi_i)\cong 0$ we have $\sum_i n_i[H:U_i]=0$ and $\sum_i n_i Ind_{U_i}^G(\chi_i)\cong 0$ and therefore $\prod_i \Delta(U_i,\chi_i)^{n_i}\lambda_{U_i}^G(\Delta)^{n_i} =1,$ if we apply {\bf (c2)} for the function $U\mapsto \lambda_U^G(\Delta)$. But then we obtain:
$$\prod_i \Delta(U_i,\chi_i)^{n_i}(\lambda_{U_i}^H)^{n_i} = \prod_i \Delta(U_i,\chi_i)^{n_i}\left(\lambda_{U_i}^G(\Delta)\lambda_H^G(\Delta)^{-[H:U_i]}\right)^{n_i}=\prod_i \Delta(U_i,\chi_i)^{n_i}\lambda_{U_i}^G(\Delta)^{n_i} =1,$$
hence $U\mapsto \lambda_U^H$ has the properties {\bf (c1), (c2)} and uniqueness implies
$$ \lambda_U^H(\Delta) = \lambda_U^H= \lambda_U^G(\Delta)\lambda_H^G(\Delta)^{-[H:U]}.$$
{\bf This finishes} the proof of the first part of the Proposition.\\
Finally we prove that {\bf the criterion 1.4 (ii), is already sufficient}:\\
If $H\le \OO$ is any other subgroup containing $N$ then we may define:
$$   \lambda_U^H(\Delta):= \lambda_U^\OO(\Delta)\cdot \lambda_H^\OO(\Delta)^{-[H:U]}.$$
The arguments from remark 2 above show that $U\in\sU_N(H)\mapsto \lambda_U^H(\Delta)\in \sA$ has the properties {\bf (c1), (c2)} too, and remark 1 shows that it is uniquely determined by these properties. Thus we have the full set of functions $\lambda_U^H(\Delta)$ at our disposal, and this is sufficient as we have seen already. {\bf qed.}\\

Still we note the following compatibility which justifies the further approach:\\

{\bf 1.5 Proposition:}  {\it For $N\le N'\le \OO'\le\OO$ such that $N'$ is normal in $\OO'$, we consider the embeddings {\bf (e1), (e)} from above, and we assume that $\Delta$ on $R_1(N\le\OO)$ extends to $\sF$ on $R(N\le \OO).$ Then the restriction $\Delta'$ of $\Delta$ onto $R_1(N'\le\OO')$ will extend to the restriction $\sF'$ of $\sF$ onto $R(N'\le\OO').$}\\

Our intention is to go in the opposite direction:

{\bf 1.6 Further approach to the extension problem:} We are given a function $\Delta$ on $R_1(\le\OO)$ (that means we assume {\bf (1), (2)}) and we want criteria
for $\Delta$ being extendible to a function $\sF$ on $R(\le\OO).$ The idea is to consider the restrictions of $\Delta$ to the subsets $R_1(N'\le\OO')$ for finite subquotients $\OO'/N'$ of $\OO,$ and then via {\bf criterion 1.4} to search for conditions to be put on $\Delta$ such that an extension onto $R(N'\le\OO')$ can be obtained. This will be done in an {\bf inductive way} that means for any natural number $n$ we assume that $\Delta$ can always be extended if $(\OO':N') < n,$ and then we are going to prove that the criteria of {\bf Theorem 3.1} are sufficient to obtain the extensions of $\Delta$  also for subquotients of index $n$. In particular we get then extensions $\sF_N$ from $\sR_1(N\le\OO)$ onto $\sR(N\le\OO)$ for all open normal subgroups $N\le\OO$, and these extensions will fit together to the extended function $\sF$ on $R(\le\OO).$\\
The induction begins from cases $N'=\OO'$ where obviously $R_1(\OO'\le\OO')= {\OO'}^*$ and $R_+(\OO'\le\OO')=R(\OO'/_{[\OO',\OO']})$ is the free $\Z$-module over ${\OO'}^*.$ Therefore in this case the extension $\sF$ of $\Delta$ is already determined by property {\bf (3)} and obviously does exist. Condition {\bf (4)} plays no role because for $N'=\OO'$ we have no induction such that we are left with {\bf (c1)}.\\

.\\

{\bf\large 2. The kernel of the Brauer map and its generating relations}\\
To apply our criterion {\bf 1.4} we need a generating system of relations {\bf (c2)}.
Here we follow essentially Deligne [D],\S1, but with some modifications because we have to deal with pairs $(\OO,N).$\\
As already explained {\bf we consider the Brauer maps (2), (3)} of the introduction:
$$\vf_\OO: R_+(\le\OO)\twoheadrightarrow R(\OO),\qquad  \vf_{N\le\OO}:R_+(N\le\OO)\twoheadrightarrow R(\OO/_{[N,N]}),$$
where $\vf_{N\le\OO}$ is the restriction of $\vf_\OO.$ (Therefore the extra notation $\vf_{N\le\OO}$ will be sometimes omitted.)\\
In view of {\bf 1.4, (c2)}  what we need are generators of $Ker(\vf_{N\le\OO}) = Ker(\vf_\OO)\cap R_+(N\le\OO)$ for all open normal subgroups $N\le\OO.$ Then, {\it in order to verify the {\bf Criterion}, it will be enough to check these generating relations.}\\

{\large\bf A)} First of all we recall {\bf the multiplication in} $R_+(\le\OO)\supseteq R_+(N\le\OO)$  which turns the Brauer maps $\vf_\OO$ and $\vf_{N\le\OO}$ into homomorphisms of commutative rings:\\
{\bf (1)}
$$ [H_1,\chi_1]\cdot [H_2,\chi_2] =\sum_{(g_1,g_2)\in H_1\backslash\OO\times H_2\backslash\OO} \left( H_1^{g_1}\cap H_2^{g_2},\; \chi_1^{g_1}\chi_2^{g_2}\right) = \sum_{g\in H_1\backslash\OO/H_2} [H_1^g\cap H_2, \chi_1^g\chi_2],$$
where the brackets in the middle indicate the pairs which will occur and {\it which have still put together into equivalence classes},
such that\\ $[H_1^{g_1}\cap H_2^{g_2},\; \chi_1^{g_1}\chi_2^{g_2}] = [H_1^{g_1g_2^{-1}}\cap H_2, \chi_1^{g_1g_2^{-1}} \chi_2] =[H_1\cap H_2^{g_2g_1^{-1}}, \chi_1\chi_2^{g_2g_1^{-1}}]$, because $[.,.]$ denotes $\OO$-equivalence classes. Therefore writing also the right side of {\bf (1)} as a sum of equivalence classes it is then
only over a subset of pairs $(g_1,g_2)$ representing the different equivalence classes, which means that $g_1g_2^{-1}\in H_1\backslash\OO/H_2$ or $g_2g_1^{-1}\in H_2\backslash\OO/H_1$ cover the different double cosets. Then indeed:
$$ \vf_\OO([H_1,\chi_1]\cdot [H_2,\chi_2]) = \vf_\OO([H_1,\chi_1])\otimes \vf_\OO([H_2,\chi_2]) = Ind_{H_1}^\OO(\chi_1)\otimes Ind_{H_2}^\OO(\chi_2),$$
because (see {\bf [S],7.3}):
$$  Ind_{H_1}^\OO(\chi_1)\otimes Ind_{H_2}^\OO(\chi_2) = Ind_{H_1}^\OO(\chi_1\otimes Res_{H_1}Ind_{H_2}^\OO(\chi_2))=$$
$$= Ind_{H_1}^\OO\left(\chi_1\;\otimes \sum_{g\in H_2\backslash\OO/H_1} Ind_{H_2^g\cap H_1}^{H_1}(\chi_2^g|_{H_2^g\cap H_1})\right)= \sum_{g\in H_2\backslash\OO/H_1} Ind_{H_2^g\cap H_1}^\OO\left( \chi_1|_{H_2^g\cap H_1}\cdot \chi_2^g|_{H_2^g\cap H_1}\right),$$
and a symmetric calculation for
$Ind_{H_2}^\OO(Res_{H_2} Ind_{H_1}^\OO(\chi_1)\otimes\; \chi_2).$\\
As a particular case we have\\
{\bf Character twist:}  If $H_2=\OO,$ hence $[H_2,\chi_2]=(\OO,\chi_2),$ then {\bf (1)} comes down to\\
{\bf  (2) }
$$  [H_1,\chi_1]\cdot (\OO,\chi_2) = [H_1,\chi_1\cdot res(\chi_2)]$$
corresponding to $Ind_{H_1}^\OO(\chi_1) \otimes\chi_2 = Ind_{H_1}^\OO(\chi_1\otimes res(\chi_2)).$ We call this the {\bf character twist} of $[H_1,\chi_1]$ by $\chi_2\in\OO^*.$\\
{\bf Remark:} Obviously the additive subgroup $R_+(N\le \OO)\subset R_+(\le\OO)$ which is generated by classes $[H,\chi]$ such that $H\ge N$, is stable under that multiplication such that $R_+(N\le\OO)$ turns into a subring, and $\vf_\OO$
restricts to a ring homomorphism $\vf_{N\le\OO}.$\\

Thus $Ker(\vf_{N\le\OO})$ is actually an ideal in $R_+(N\le\OO)$, and {\it our aim in this section is to determine a set of generators for $Ker(\vf_{N\le\OO})$ considered as additive group}. The main result will be {\bf Theorem 2.7}.\\

{\large\bf B)}  As in {\bf [D]},1.9. we recall {\bf some functoriality properties of the Brauer maps} $\vf_\OO.$ It is easy to adapt this to the restricted maps $\vf_{N\le\OO}$ which we leave to the reader.\\
{\bf (F1):}  If $\OO'\le\OO$ is an open subgroup, then the map
$$Ind_{\OO'}^\OO:R_+(\le \OO')\rightarrow R_+(\le \OO),\qquad [H',\chi]_{\OO'}\mapsto [H',\chi]_\OO,$$
which is induced by the inclusion $R_1(\le\OO')\subseteq R_1(\le\OO)$ and takes $\OO'$-classes into $\OO$-classes, is called induction. This is justified by the commutative diagram:
$$\xymatrix{
R_+(\le \OO') \ar[d]^{\vf_{\OO'}} \ar[r]^{Ind} & R_+(\le \OO)\ar[d]^{\vf_\OO} \\
R(\OO') \ar[r]^{Ind} & R(\OO)  }$$\\
{\bf (F2):}   Let $u:\OO'\rightarrow\OO$ be a continuous homomorphism with finite cokernel. Then we define the corresponding restriction $u^*:R_+(\le \OO)\rightarrow R_+(\le \OO')$ by:
$$  u^*\circ [H,\chi]_\OO = \sum_{g\in H\backslash\OO/u(\OO')} [u^{-1}(H^g),\;\chi^g\circ u]_{\OO'},$$
 where an $\OO$-orbit is mapped to a sum of $\OO'$-orbits, and we will write $u^*= Res_{\OO'}^\OO$ if $u$ is an inclusion.
Again under $\vf$ this is transformed into the restriction map $u^*:R(\OO)\rightarrow R(\OO'),$ that means: $u^*\circ\vf_\OO= \vf_{\OO'}\circ u^*:R_+(\le\OO)\rightarrow R(\OO'),$ as we see from {\bf [S], 7.3}:
 $$ Ind_H^\OO(\chi)\circ u = \left(Res_{u(\OO')}Ind_H^\OO(\chi)\right)\circ u =\left(\sum_{g\in H\backslash\OO/u(\OO')} Ind_{H^g\cap u(\OO')}^{u(\OO')}(\chi^g|_{H^g\cap u(\OO')})\right)\circ u =$$
 $$=\sum_{g\in H\backslash\OO/u(\OO')} Ind_{u^{-1}(H^g)}^{\OO'}(\chi^g\circ u|_{u^{-1}(H^g)}) = \vf_{\OO'}\circ u^*([H,\chi]_\OO).$$\\
 {\bf (F3):}   If $u:\OO'\twoheadrightarrow\OO$ is {\bf surjective}, then $H\backslash\OO/u(\OO')$ is only one element, and we obtain:
 $$ u^*:R_+(\le \OO)\rightarrow R_+(\le \OO'),\qquad u^*\circ [H,\chi]_\OO = [u^{-1}(H), \chi\circ u]_{\OO'},$$
 which is called the {\bf inflation map} $u^*=Inf_\OO^{\OO'}.$\\

{\bf 2.1 Proposition:} ([D],(1.9.6)) {\it
If $\rho\in R_+(\le \OO)$ is any element, then the product with a generator $[H,\chi]\in R_+(\le \OO)$ can be rewritten as:
$$ \rho\cdot [H,\chi] = Ind_H^\OO(Res_H(\rho)\cdot (H,\chi))\quad\in R_+(\le \OO) ,$$
where $Res_H(\rho)\cdot (H,\chi)$ is the product in $R_+(\le H),$  and $Ind_H^\OO : R_+(\le H)\rightarrow R_+(\le \OO)$.}\\

For the proof we may assume that $\rho= [H_1,\chi_1]$ is a generator of $R_+(\le \OO)$ and then use {\bf (F2)} where $[H,\chi]$ is replaced by $[H_1,\chi_1]$ and $u^*=Res_H^\OO.$ This brings us down to the twist operation {\bf (2)} on the level of $H.$\\

We recall that $\rho = \sum_i n_i [H_i,\chi_i] \in Ker(\vf_\OO)$ means a relation
$  \sum_i n_i Ind_{H_i}^\OO (\chi_i) \cong 0\in R(\OO).$

If $\OO'$ is a subquotient of $\OO$ which means if $\OO'\stackrel{u}{\twoheadleftarrow}\OO'' \le \OO,$ then our functoriality properties fit together to a commutative diagram
$$\xymatrix{
R_+(\le \OO') \ar[d]^{\vf_{\OO'}} \ar[r]^{Inf} & R_+(\le \OO'')\ar[d]^{\vf_{\OO''}} \ar[r]^{Ind} & R_+(\le \OO)\ar[d]^{\vf_\OO} \\
R(\OO') \ar[r]^{Inf} & R(\OO'') \ar[r]^{Ind} & R(\OO)  }$$

Moreover, since  $Ker(\vf)$ is always an ideal and taking into account  the character twist operation {\bf (2)}, we see from the diagram that\\
{\bf (3)}
$$   \sigma\in R_+(\le \OO') \mapsto        Ind_{\OO''}^\OO(\chi\cdot Inf_{\OO'}^{\OO''}(\sigma)),\quad\textrm{where}\;\chi\in {\OO''}^*,$$
will take $Ker(\vf_{\OO'})$ into $Ker(\vf_{\OO}).$\\

{\large\bf C)} Our next aim is to {\bf produce elements $\rho\in Ker(\vf_\OO)$ in a systematic way}.\\ Here we restrict to the case where $N=C$ is a {\it commutative normal subgroup}, hence $[C,C]=\{1\}$ and $\vf_{C\le\OO}:R_+(C\le\OO)\rightarrow R(\OO).$  Then we are going to define a {\bf projector}
$$  \Phi_C:R_+(\le\OO) \twoheadrightarrow R_+(C\le\OO)\subset R_+(\le\OO),\qquad \Phi_C\circ\Phi_C =\Phi_C,$$
such that $ \vf_\OO(\rho) = \vf_{C\le\OO}\circ\Phi_C(\rho),$ hence
$\rho-\Phi_C(\rho)\in Ker(\vf_\OO)\quad\forall\;\rho\in R_+(\le\OO),$
{\it and together with $\vf_\OO$ also the restriction $\vf_{C\le\OO}$ is surjective}.
\begin{itemize}
 \item We begin from $[H,\chi]\in R_+(\le\OO)$ and produce a certain decomposition of $Ind_H^\OO(\chi)$ {\bf which depends on $C$:}
\end{itemize}
$H$ acts by conjugation on $C$ which induces an action $\mu\in C^*\mapsto \mu^h=h^{-1}\mu h,$ that means: $\mu^h(c)=\mu(hch^{-1}).$ The subset
$$   S=S(\chi):=\{\mu\in C^*\;|\; \mu|_{H\cap C} = \chi|_{H\cap C}\}, \qquad \#S = (C:H\cap C),\qquad S(\chi)=\mu_0\cdot S(1),$$
is stable under that action, and for $\mu\in S$ {\bf the function $\chi\mu$:}
$$ hc\in HC \mapsto  (\chi\mu)(hc):=\chi(h)\mu(c)$$
is well defined because $hc=h'c'$ implies $h^{-1}h'= c{c'}^{-1}\in H\cap C,$ hence $\chi(h^{-1}h')=\mu(c{c'}^{-1}).$\\
But we want $\chi\mu$ to be a character, in particular it should be trivial on commutators. Since $hch^{-1}c^{-1}\in C$ a necessary condition is:
$$  1=(\chi\mu)(hch^{-1}c^{-1}) = \mu(hch^{-1}c^{-1}) = \mu^{h-1}(c),\quad\forall c\in C,$$
hence $h$ must be in the isotropy group $H_\mu:=\{h\in H\;|\;\mu^h =\mu\}$. But this is also sufficient that means {\it the restriction $\chi\mu\in (H_\mu C)^*$ is actually a character.} And conjugation by $g\in\OO$ turns it into the character $(\chi\mu)^g = \chi^g \mu^g$ of $g^{-1}(H_\mu C)g = (H^g)_{\mu^g}\cdot C.$ In particular:
\begin{itemize}
 \item $(\chi\mu)^h =\chi\cdot\mu^h$ is a character of $(H_\mu C)^h = H_{\mu^h}\cdot C.$
\end{itemize}
Finally let $T(\chi)\subset S=S(\chi)$ be a system of representatives for the orbits $[\mu]_H=\{\mu^h\}_{h\in H}\subseteq S.$\\

{\bf 2.2 Lemma:} (for (i) see: {\bf [D] 1.11., [L] 15.1} has the version for $\chi\equiv 1.$)\\
{\it Depending on a {\bf commutative normal subgroup} $C\le \OO$ and with notations as we have fixed them, we will have the relation\\
{\bf (i)}
$ \quad  Ind_H^\OO(\chi)= Ind_{HC}^\OO\left(Ind_H^{HC}(\chi)\right) = \sum_{\mu\in T(\chi)} Ind_{H_\mu C}^\OO (\chi\mu),$
$$[H,\chi] -\sum_{\mu\in T(\chi)} [H_\mu C,\chi\mu]\quad \in  Ker(\vf_\OO),$$
where on the right side the induction is via groups $H_\mu C$ such that $HC\ge H_\mu C\ge C.$ In the particular case where $H\ge C$ the relation is simply the identity because $T(\chi)=\{\chi|_C\}$ is then a single element and $HC=H.$ And for $\chi\equiv 1$ we obtain:
$T(1)=(C/H\cap C)^*/H$  and $\chi\mu=\mu'\in (H_\mu C)^*$ is the extension of $\mu$ by $1$.\\
{\bf (ii)}  If $H\le \OO$ is normalized by $C$, then:
$$ \quad Ind_H^\OO(\chi) = Ind_{HC}^\OO\left( Ind_{H'C}^{HC}(\chi \mu_0)\otimes (\sum_{\nu\in (C/H\cap C)^*/\sim} \nu')\right),$$
$$[H,\chi] - \sum_{\nu\in (C/H\cap C)^*/\sim} [H'C, \chi(\mu_0\nu)]\quad \in Ker (\vf_\OO),$$
where $H'=Stab_H(\mu)$ is for all $\mu\in S=S(\chi)$  the projective kernel of $Ind_H^{HC}(\chi),$\\
$\mu_0\in S$ is any fixed constituent, and\\
$\nu'\in (HC/H)^*$ is the extension of $\nu\in (C/H\cap C)^*$ for a system of representatives $\nu$ with respect to the equivalence relation: $\nu_1\sim\nu_2$ if and only if $\nu_1\nu_2^{-1}=\mu_0^{h-1}$ for some $h\in H.$}\\

{\bf Proof: (i)} Since $H_\mu C\subseteq HC,$ we may rewrite our assertion as
$$  Ind_{HC}^\OO\circ Ind_H^{HC}(\chi)= Ind_{HC}^\OO\left(\sum_{\mu\in T(\chi)} Ind_{H_\mu C}^{HC} (\chi\mu)\right),$$
so it is sufficient to prove the Lemma  in the case $\OO=HC.$
The $\OO$-module $Ind_H^{HC}(\chi)$ is spanned by complex functions $f:HC\rightarrow\fC,$ such that
$$  f(hx)=\chi(h)f(x),\;\forall\;h\in H,\;x\in HC\qquad (gf)(x):= f(xg)\;\textrm{for}\;g\in\OO.$$
Because of $\OO=HC$ and $S=S(\chi)$ we have $\#S=(C:H\cap C)=(HC:H)=(\OO:H),$ and the elements $\mu\in S$ give rise to
a distinguished basis $\{f_\mu:=\chi\mu\}_{\mu\in S}$ for that space of functions:  $f_\mu(hc):= \chi(h)\mu(c)$ is well defined as we have seen above already.\\
The function $f_\mu:HC\rightarrow\fC$ extends the character $\mu$ and each of the functions $f_\mu$ has the property  $f_\mu(hx) =\chi(h)f_\mu(x)$. To see this put $x=h'c'.$\\
Moreover $f_\mu(hc\cdot c') = f_\mu(hc)\mu(c'),$ hence $f_\mu(gc)=f_\mu(g)\mu(c),$ for any $g\in HC$ and $c\in C.$  Thus we obtain:
$$ (gf_\mu)(hc) = f_\mu(hcg)= f_\mu(h\cdot g\cdot g^{-1}cg)=\chi(h) f_\mu(g) \mu(g^{-1}cg) = f_\mu(g)\cdot f_{g\mu g^{-1}}(hc).$$
We see that the action of $HC$ on the basis $\{f_\mu\}_{\mu\in S}$ is monomial that means up to scalar factors the basis is permuted:\\
{\bf (4)}
$$   gf_\mu = f_\mu(g)\cdot f_{g\mu g^{-1}},\qquad (hc)f_\mu = \chi(h)\mu(c)\cdot f_{h\mu h^{-1}},$$
where $g\mu g^{-1}=h\mu h^{-1},$ because the action of $HC$ on $S\subseteq C^*$ factors to an action of $HC/C=H/C\cap H.$ Therefore for each $H$-orbit $[\mu]_H\subseteq S$ we have the subspace $V_{[\mu]_H}$ with basis $\sB_{[\mu]_H}=\{f_\mu=\chi\mu\}_{\mu\in [\mu]_H},$ such that
$$ Ind_H^{HC}(\chi) = \bigoplus_{[\mu]_H\in S/H} V_{[\mu]_H}$$
is a decomposition into $HC$-submodules. Now the assertion {\bf (i)} will follow from $V_{[\mu]_H} \cong Ind_{H_\mu C}^{HC}(\chi\mu)$ by verifying the character identity\\
{\bf (5)}
$$  Tr(hc; V_{[\mu]_H}) = \chi(h)\sum_{\mu;\;h\in H_\mu} \mu(c)=Tr(hc,\; Ind_{H_\mu C}^{HC}(\chi\mu)).$$

{\bf (ii):} If $H$ is normalized by $C$, then $[H,C]\subseteq H\cap C$ and therefore the inner action of $H$ on $C/H\cap C$ will be trivial and
$$  \mu^h(c)= \mu(hch^{-1}) =\mu([h,c]\cdot c) = \chi([h,c])\cdot \mu(c), \qquad \mu^{h-1}(c) = \chi([h,c]),$$
such that $\mu^{h-1}\in C^*$ is for all $\mu\in S$ the same character.
Therefore independently of $\mu$ we will have $\mu^h =\mu$ if and only if $\chi([h,c])=\chi^{1-c}(h)=1$ for all $c\in C.$ This means the isotropy group $H_\mu=:H'$ and the group $H_\mu C = H'C$ {\bf will not depend on $\mu\in S$.} Moreover
$$  Ind_H^{HC}(\chi)|_H = \sum_{c\in C/H\cap C} \chi^c$$
implies that $H' =\{h\in H\;|\; \chi^c(h)=\chi(h)\;\forall\;c\}$ is the projective kernel of $Ind_H^{HC}(\chi).$\\

{\bf Sublemma:} {\it Assume $H\le HC$ normal and let $H'= Stab_H(\mu)$ for all $\mu\in S=S(\chi).$ Then fixing any $\mu_0\in S$
 we obtain an injection $h\in H/H'\mapsto  \mu_0^{h-1}\in (C/H\cap C)^*.$}\\
 {\bf Proof:} Since $\mu_0^h$ and $\mu_0$ are both in $S$ we obtain $\mu_0^{h-1}\in (C/H\cap C)^*$, hence $\mu_0^{(h_1-1)(h_2-1)} \equiv 1,$ because the inner action of $H$ on $(C/H\cap C)$ is trivial, and this implies:
 $$  h_1h_2 \mapsto \mu_0^{h_1h_2 -1} = \mu_0^{h_1-1}\cdot \mu_0^{h_2-1}.$$

Now returning to the proof of {\bf (ii)} we use that $S = \mu_0\cdot(C/H\cap C)^*,$ for any fixed constituent $\mu_0\in S,$ and equivalence $\mu_0\nu_1\sim\mu_0\nu_2$ with respect to the inner action of $H,$ means:\\
$\mu_0\nu_1 =(\mu_0\nu_2)^h= \mu_0^h\nu_2$, such that $\nu_1\nu_2^{-1}=\mu_0^{h-1}$ is in the image of the Sublemma-map. This turns relation {\bf (i)} into relation {\bf (ii)}. {\bf qed.}\\

{\bf 2.3 Corollary:} (The projector $\Phi_C$):\\ {\it {\bf (i)} Via Lemma 2.2 we obtain (for a fixed abelian normal subgroup $C\le\OO$) a well defined additive homomorphism
$$ \Phi_C: R_+(\le \OO)\twoheadrightarrow R_+(C\le \OO)\subset R_+(\le\OO),\qquad [H,\chi]\mapsto \sum_{\mu\in T(\chi)} [H_\mu C,\chi\mu]$$
which is actually a projector: $\Phi_C^2 =\Phi_C,$ and such that with respect to the Brauer map:
$$  \vf_\OO = \vf_{C\le \OO}\circ \Phi_C,\qquad \vf_\OO\circ(Id -\Phi_C)=\vf_{C\le \OO}\circ\Phi_C\circ(Id -\Phi_C)=0.$$
Moreover $\Phi_C$ is a homomorphism of modules over the ring $R_+(\OO\le \OO)=R(\OO/_{[\OO,\OO]})$ of the abelian quotient of $\OO.$ \\
{\bf (ii)} If $C$ is in the center of $\OO$, then $[H,H]=[HC,HC]$ that means $\chi\in H^*$ can be extended to a character of $HC,$ hence $Ind_H^{HC}(\chi) =\sum_{\mu\in S(\chi)} \chi\mu$ is the sum of all extensions of $\chi$, and:
$$ Ind_H^\OO(\chi)=\sum_{\mu\in S(\chi)} Ind_{HC}^\OO(\chi\mu),\qquad  \Phi_C([H,\chi]) = \sum_{\mu\in S(\chi)} [HC,\chi\mu].$$
{\bf (iii)} If $C'\ge C$ is a larger abelian normal subgroup of $\OO$, then $R_+(C\le\OO) \supseteq R_+(C'\le\OO)$ and
$$  \Phi_{C'}\circ\Phi_C =\Phi_{C'}.$$
{\bf (iv)} If $C\le\OO'\le\OO$ then the relative projector $\Phi_C^{\OO'}$ fits into the commutative diagram
$$\xymatrix{
R_+(\le \OO') \ar[d]^{\Phi_C^{\OO'}} \ar[r]^{Ind_{\OO'}^\OO}  &  R_+(\le\OO) \ar[d]^{\Phi_C^\OO}\\
R_+(C\le \OO') \ar[r]^{Ind_{\OO'}^\OO} &  R_+(C\le\OO) },$$
where $Ind_{\OO'}^\OO$ means to replace $\OO'$-conjugacy classes by the corresponding $\OO$-conjugacy classes.}\\

{\bf Proof: (i):} From the proof of the previous Lemma it is easy to see that $[H,\chi]\mapsto [HC,\chi\mu]$ is a well defined map of $\OO$-conjugacy classes which does not depend on the choice of a representative $\mu\in T(\chi).$ Thus $\Phi_C([H,\chi])$
is well defined and extends to arguments in $R_+(\le \OO)$ which by definition is the free $\Z$-module over the conjugacy classes $[H,\chi].$ And in Lemma 2.2 we have mentioned already that $\Phi_C([H,\chi])=[H,\chi]$ if $H\ge C,$  hence
$ \Phi_C\circ\Phi_C([H,\chi]) = \Phi_C([H,\chi]).$ Also the property $\vf_\OO = \vf_{C\le\OO}\circ\Phi_C$ is a direct consequence of {\bf 2.2}.\\
Finally in $R_+(\le \OO)$ we have (see {\bf (2)}) the twist operation:  $(\OO,\eta)\cdot [H,\chi] = [H,\eta_H\chi],$ and therefore:
$$ \Phi_C((\OO,\eta)\cdot [H,\chi]) = \Phi_C([H,\eta_H\chi]) = \sum_{\mu\in T(\eta_H\chi)} [H_\mu C,\eta_H\chi\mu].$$
Since $\eta\in\OO^*$ we have $S(\eta_H\chi) =\eta_C\cdot S(\chi),$ $\mu = \eta_C\mu_1$ and $H_\mu = H_{\mu_1}.$ Therefore the last displayed term rewrites as:
$$  =\sum_{\mu_1\in T(\chi)} [H_{\mu_1}C, (\eta_H\chi)(\eta_C\mu_1)],$$
where $(\eta_H\chi)(h)\cdot(\eta_C\mu_1)(c)=\eta(hc)\cdot \chi(h)\mu_1(c)=\eta(hc)(\chi\mu_1)(hc),$ and we may proceed:
$$ ... =\sum_{\mu_1\in T(\chi)} [H_{\mu_1}C, Res_{H_{\mu_1}C}(\eta)\chi\mu_1]= (\OO,\eta)\cdot\left(\sum_{\mu_1\in S(\chi)/H} [H_{\mu_1}C,\chi\mu_1]\right) =(\OO,\eta)\cdot \Phi_C([H,\chi]).$$
{\bf (ii)} is obvious.\\
{\bf (iii):} If $C'\ge C$ are both central in $\OO$ then this follows from {\bf (ii)} because extending $\chi\in H^*$ to $HC'$ is the same as first extending it from $H$ to $HC$ and then from $HC$ to $(HC)C'=HC'.$ The general case is less obvious. By definition we have:
$$ \Phi_C([H,\chi])= \sum_{\mu\in S_C(\chi)/H} [H_\mu C,\chi\mu],\qquad
\Phi_{C'}([H,\chi])= \sum_{\mu'\in S_{C'}(\chi)/H} [H_{\mu'} C',\chi\mu'],$$ where $S_{C'}(\chi)=\{\mu'\in(C')^*\;|\; \mu'|_{H\cap C'}=\chi|_{H\cap C'}\}.$ Since $C'$ is commutative too, we have the surjective restriction map:
$$ res_C^{C'}: S_{C'}(\chi)\twoheadrightarrow S_C(\chi)$$
which is compatible with $H$-conjugation. Thus we are left to show that:
$$ \Phi_{C'}([H_\mu C,\chi\mu]) = \sum_{\mu'\in S_{C'}(\chi)/H;\,\mu'|_C=\mu} [H_{\mu'}C',\chi\mu'].$$
As to the left side, by definition we have:
$$ \Phi_{C'}([H_\mu C,\chi\mu]) = \sum_{\mu'\in S_{C'}(\chi\mu)/(H_\mu C)} [(H_\mu C)_{\mu'} C', (\chi\mu)\mu'],$$
where $S_{C'}(\chi\mu)=\{\mu'\in (C')^*\;|\; \mu'|_{H_\mu C\cap C'}= (\chi\mu)|_{H_\mu C\cap C'}\}.$ Since $H_\mu C\cap C'=(H_\mu\cap C')C$
and (see before {\bf 2.2} above):  $\chi\mu =\chi$ on $H_\mu\cap C'\le H,$ and  $\chi\mu=\mu$ on $C,$ we conclude:\\ $\mu'\in S_{C'}(\chi\mu)$ implies $\mu'|_{H_\mu\cap C'}=\chi|_{H_\mu\cap C'}$ and $\mu'|_C =\mu.$\\
Moreover $\mu'|_C = \mu$ implies $H_{\mu'}\le H_\mu,$ and together with
$H\cap C'\le H_{\mu'}$ we obtain:
$$  H_\mu\cap C'\le H\cap C'\le H_{\mu'}\cap C' \le H_\mu\cap C',$$
such that all inclusions are equalities. Therefore for $\mu\in S_C(\chi)$ we will have:\\
$\mu' \in S_{C'}(\chi\mu)$ if and only if $\mu'\in S_{C'}(\chi)$ and $\mu'|_C=\mu.$\\
Finally from $H_{\mu'}\le H_\mu$ we see that $(H_\mu C)_{\mu'}=H_{\mu'}C$ and therefore $(H_\mu C)_{\mu'} C' = H_{\mu'}C'$ and $(\chi\mu)\mu'=\chi\mu'.$ This proves {\bf (iii)},\\
and the assertion {\bf (iv)} is easy to see.  {\bf qed.}\\

As an immediate consequence for arbitrary normal subgroups we note:\\

{\bf 2.4 Corollary:} (The projector $\Phi_N$.)\\
{\bf (i)} {\it If $N\le \OO$ is any open normal subgroup, then we denote
 $$ \Phi_N:R_+(\le\OO/_{[N,N]}) \rightarrow R_+(N/_{[N,N]} \le \OO/_{[N,N]})= R_+(N\le\OO)$$
 the projector which is induced modulo $[N,N]$, such that $\vf_\OO(\rho) = \vf_{N\le \OO}\circ\Phi_N(\rho)\in R(\OO/_{[N,N]}).$\\
{\bf (ii)} If $N=\OO$ then $R_+(\le\OO/_{[\OO,\OO]})$ is generated by pairs $[H,\chi]$ such that $H\ge [\OO,\OO]$, $\chi\in(H/_{[\OO,\OO]})^*,$ and this implies
$$ \Phi_\OO([H,\chi]) =\sum_{\wt\chi} (\OO,\wt\chi)\in R_+(\OO\le\OO),$$
where $\wt\chi\in\OO^*$ runs over the extensions of $\chi.$}\\

{\large\bf D)} A list of basic relations.\\
{\bf 2.5 Basic relations:} Now we give a list of basic relations $\rho\in Ker(\vf_\OO)\subset R_+(\le\OO)$ or more special $\rho\in Ker(\vf_{N\le\OO})\subset R_+(N\le\OO),$ the members of which will actually be generators. Some additional comments will follow below. We begin by presenting those generators which refer directly to the group $\OO$:
\begin{itemize}
 \item[{\bf I}] For $K<\OO$ a normal subgroup of prime index $\ell$ and $\chi\in \OO^*$ let
 $$ \rho(K,\chi):= [K,\chi_K] -\Phi_\OO([K,\chi_K]) =
 [K,\chi_K] - \sum_{\mu\in(\OO/K)^*} (\OO,\chi\mu),$$
 $$ \vf_\OO(\rho(K,\chi)) = Ind_K^\OO(\chi_K) - \sum_{\mu\in(\OO/K)^*} \chi\mu\quad \cong 0.$$
 \item[{\bf II}] For $Z<\OO$ a normal subgroup and $\eta\in Z^*$ with the properties
 \begin{itemize}
  \item  $\OO/Z\cong \Z/\ell \times \Z/\ell$ is abelian bicyclic for some prime $\ell,$ hence $[\OO,\OO]\le Z$,
  \item  the commutator induces an isomorphism
 $[.,.]:\OO/Z\wedge \OO/Z \cong [\OO,\OO]/[Z,\OO]$ of cyclic groups of order $\ell,$ (in particular $[Z,\OO]\ne [\OO,\OO]$),
 \item  $\eta:Z/[Z,\OO]\rightarrow \fC^\times$ is a character which is non-trivial on $[\OO,\OO]/[Z,\OO].$
 \end{itemize}
 consider pairs $(H,\eta^H)$ where $Z< H<\OO$ is of index $\ell$ and $\eta^H:H/[Z,\OO]\rightarrow\fC^\times$ extends $\eta.$ {\bf Then put}
 $$  \rho(Z,\eta): = [H_1,\rho^{H_1}] - [H_2,\rho^{H_2}] = [H_1,\rho^{H_1}] - \Phi_{H_2}([H_1,\rho^{H_1}]),$$
 $$  \vf_\OO(\rho(Z,\eta)) = Ind_{H_1}^\OO(\chi_{H_1}) - Ind_{H_2}^\OO(\chi_{H_2})\quad\cong 0.$$
 \item[{\bf III}] For $H<\OO$ a (proper) {\bf maximal} subgroup which is {\bf not normal}, that means $K:=\bigcap_{g\in\OO} gHg^{-1}\ne H$ and therefore $\OO/K$ is a non-degenerate type-III-group (as defined in Appendix 2), let $C<\OO$ be the uniquely determined normal subgroup such that:
 $$   H\cdot C=\OO,\qquad H\cap C =K,\quad(\textrm{hence}\;C/K\;\textrm{ is abelian}),$$
 and for characters $\chi\in \OO^*$ put:
 $$\rho(H,\chi):= [H,\chi_H] - \Phi_C([H,\chi_H])=[H,\chi_H]-\sum_{\mu\in(C/K)^*/H} [H_\mu C,\chi_{H_\mu C}\cdot\mu']$$
 referring to the inner action of $H$ on $(C/K)^*$ and $\mu'$ the trivial extension of $\mu\in (C/K)^*$ onto the stabilizer $H_\mu C$, which implies
 $$ \vf_\OO(\rho(H,\chi))= Ind_H^\OO(\chi_H) - \sum_{\mu\in(C/K)^*/H} Ind_{H_\mu C}^\OO(\chi_{H_\mu C}\cdot\mu')\quad\cong 0.$$
\end{itemize}
{\bf Further comments:} In {\bf II} the group $H/Z$ is cyclic, hence $H/Z\wedge H/Z =\{1\},$ and therefore $[H,H]\le [Z,\OO],$ whereas $[H,\OO]=[\OO,\OO]$ because $\OO/Z=H/Z \times H'/Z$. Thus $\eta^H$ is nontrivial on $[H,\OO]/[Z,\OO]$ which means that the $\OO$-conjugates of $\eta^H$ are all different, hence $Ind_H^\OO(\eta^H)$ will be irreducible. This also means that the conjugacy class $[H,\eta^H]_\OO$ consists of all extensions of $\eta$ onto $H >Z$.\\
Nevertheless $\rho(Z,\eta)$ is not a single element because we can form it for each pair $(H_1, H_2)$ of intermediate groups $Z< H_i <\OO.$\\
In {\bf III} the normal subgroup $C<\OO$ will usually not be abelian. Nevertheless $\Phi_C([H,\chi_H])\in R_+(C/_{[C,C]}\le \OO/_{[C,C]})$ is well defined because $H> K\ge [C,C]$ and $\chi_H$ is trivial on $H\cap[\OO,\OO]>[C,C],$ hence $[H,\chi_H]\in R_+(\le\OO/_{[C,C]}).$\\

Similarly, {\bf for all open subgroups} $B\le\OO$ we may consider elements
$$ \rho_{B}(K,\chi),\;\rho_{B}(Z,\eta),\;\rho_{B}(H,\chi)\quad\in Ker(\vf_{B})\subset R_+(\le B),$$
referring now to subgroups $K,\;Z,\;H$ of $B$, {\bf and the complete set of generators consists of }:
\begin{itemize}
 \item[{\bf I}] elements $\rho(K,B,\chi):=Ind_{B}^\OO(\rho_{B}(K,\chi)),$
 \item[{\bf II}] elements $\rho(Z,B,\eta):=Ind_{B}^\OO(\rho_{B}(Z,\eta)),$
 \item[{\bf III}] elements $\rho(H,B,\chi):=Ind_{B}^\OO(\rho_{B}(H,\chi)),$
\end{itemize}
where $Ind_{B}^\OO:R_+(\le B)\rightarrow R_+(\le\OO)$ is the enlargement of conjugacy classes which according to {\bf (F1)} takes $Ker(\vf_{B})$ into $Ker(\vf_\OO).$\\

{\bf 2.5/N  Basic relations $\rho\in Ker(\vf_{N\le\OO})$}:\\  If moreover $N\le\OO$ is a fixed open normal subgroup then we consider the\\
subsets {\bf I/N, II/N, III/N} consisting of those generators which are in $R_+(N\le\OO),$ hence they are in $Ker(\vf_\OO)\cap R_+(N\le\OO) = Ker(\vf_{N\le\OO}).$ This means that all occurring subgroups $K,\;Z,\;H,\;B$, and for {\bf III/N} also $K=\bigcap_{g\in B} gHg^{-1}$ must contain $N.$\\

{\bf 2.6 Definition:} We let $\sR(\le\OO)\subseteq Ker(\vf_\OO)$ and $\sR(N\le\OO)\subseteq Ker(\vf_{N\le\OO})$ be the {\it additive subgroups} which are generated by the above elements of types {\bf I - III} and {\bf I/N - III/N} resp. This means that:
$$  \sR(N\le\OO) \subseteq \sR(\le\OO) \cap R_+(N\le\OO).$$

{\bf Remarks:} {\bf 2.6.1  If $\OO$ is abelian then the types II and III will not occur} because they refer to non-abelian groups.
{\bf If $\OO$ is nilpotent then types III will not occur} because in a nilpotent group $B$ maximal subgroups $H <B$ will always be normal subgroups. Moreover in accordance with {\bf (3)} above, we have:\\
{\bf 2.6.2} {\bf Twisting} basic relations with a character $X\in \OO^*$ will stabilize this set and therefore the additive groups $\sR(\le\OO)$ and $\sR(N\le\OO)$ will be modules over the character ring $R_+(\OO\le\OO) =R(\OO/_{[\OO,\OO]}).$  {\it And type preserving we have:}\\
{\bf 2.6.3} ({\bf Induction})
$$  Ind_H^\OO\circ \sR(\le H)\subseteq \sR(\le\OO),\qquad Ind_H^\OO\circ \sR(N\le H)\subseteq \sR(N\le\OO)$$
if $H\le\OO$ is an open subgroup, and if $H\in\sU_N(\OO)$ resp.\\
{\bf 2.6.4} ({\bf Inflation}) If $u:\OO'\twoheadrightarrow\OO$ is a surjection, and $u(N')\le N,$ then the corresponding inflation map $[H,\chi]_\OO\mapsto [u^{-1}(H),\chi\circ u]_{\OO'}$ (see (F3) above)
will induce:
$$  u^*\circ \sR(\le\OO)\subseteq \sR(\le\OO'),\qquad u^*\circ \sR(N\le\OO)\subseteq \sR(u^{-1}(N)\le\OO')\subseteq \sR(N'\le\OO').$$

{\large\bf E)}  The main result of this section, and reducing it to the case of a central extension:\\
{\bf 2.7 Theorem:} {\it If $\OO$ is a solvable profinite group then the inclusions $\sR(\le\OO)\subseteq Ker(\vf_\OO)$ and $\sR(N\le\OO)\subseteq Ker(\vf_{N\le\OO})$ are actually equalities.}\\

{\bf Proof:} Obviously it is enough to prove that $\sR(N\le\OO)= Ker(\vf_{N\le\OO})$ because, due to\\ $R_+(\le\OO)=\bigcup_N R_+(N\le\OO),$ hence $\sR(\le\OO) = \bigcup_N\sR(N\le\OO)$ and $Ker(\vf_\OO)=\bigcup_N Ker(\vf_{N\le\OO}),$ the equality $\sR(\le\OO)= Ker(\vf_\OO)$ will follow if the $N$-equalities hold for all open normal subgroups $N\le\OO.$ {\bf We begin by reducing this to cases where $N$ is in the center of $\OO.$}\\

For the orbits $\xi=[N,\mu]\in N^*/\OO$ we have a natural decomposition\\
{\bf (6)}
$$   R_+(N\le\OO) = \bigoplus_{\xi\in N^*/\OO} R_+(N\le\OO)_\xi,\qquad r= \sum_\xi r_\xi,$$
where {\bf $R_+(N\le\OO)_\xi$ is the free $\Z$-module over all $[H,\chi]\ge \xi =[N,\mu],$} that means $\chi|_N$ is a conjugate of $\mu.$ Because of $(N^*)^\OO=(N/_{[N,\OO]})^*$ we obtain the  partial sums:
$$  R_+(N/_{[N,\OO]}\le \OO/_{[N,\OO]}) = \bigoplus_{\xi=\mu \in(N^*)^\OO} R_+(N\le\OO)_\xi,\qquad R_+(\le \OO/N) = R_+(N\le\OO)_{\xi\equiv 1}.$$
If a {\bf basic relation} $\rho\in \sR(\le\OO)$ is contained in $R_+(N\le\OO)$ for some normal subgroup $N\le\OO,$ that means if $\rho$ is a member of the above subsets {\bf I/N, II/N, III/N} then all the constituents $[H,\chi]$ of $\rho$ will have the same $\xi=\xi([H,\chi])\in N^*/\OO,$ such that actually  $\rho\in R_+(N\le\OO)_\xi$ for a uniquely determined $\xi=\xi(\rho)\in N^*/\OO.$ More precisely for $\rho=\rho(K,B,\chi),\;\rho(H,B,\chi)$ of type I/N, III/N resp. the conjugacy class $\xi(\rho)$ will be generated by the restriction $\chi|_N,$ whereas for $\rho=\rho(Z,B,\eta)$ of type {\bf II/N} the restriction of $\eta\in (Z^*)^B$ to $N\le Z$ will do it.  Therefore the decomposition {\bf (6)} induces\\
{\bf (7)}
$$  \sR(N\le\OO) =\bigoplus_{\xi\in N^*/\OO} \sR(N\le\OO)_\xi,$$
where {\bf $\sR(N\le\OO)_\xi$ is the additive group over all basic relations $\rho$ with $\xi(\rho)=\xi.$}\\ And as partial sums we have:
$$ \sR(N/_{[N,\OO]}\le \OO/_{[N,\OO]}) = \bigoplus_{\xi=\mu \in(N^*)^\OO} \sR(N\le\OO)_\xi,\qquad \sR(\le \OO/N) = \sR(N\le\OO)_{\xi\equiv 1}.$$
Let $N\le\OO$ be a normal subgroup and $\ol\OO=\OO/_{[N,N]}.$
The Grothendieck-group $R(\ol\OO)$ has the set of irreducible representations $\rho$ as ON-basis. And by restriction to the abelian normal subgroup $\ol N$ each $\rho$ determines an $\OO/N$-orbit $\xi(\rho) =[N,\mu]\in N^*/\OO$ (with respect to conjugation). This gives us an orthogonal decomposition\\
{\bf (8)}
$$ R(\ol\OO) = \bigoplus_{\xi\in N^*/\OO} R(\ol\OO)_\xi $$
into submodules generated by all irreducible representations $\rho$ with a fixed $\xi(\rho) =[N,\mu].$\\

{\bf 2.8 Lemma:} {\it The Brauer map $\vf_{N\le\OO}: R_+(N\le\OO) \twoheadrightarrow R(\ol\OO)$ respects the decompositions {\bf (6)} and {\bf (8)} resp. that means we obtain $\vf_\xi:R_+(N\le\OO)_\xi\rightarrow R(\ol\OO)_\xi$ such that
$$ \vf_{N\le\OO} =  \sum_{\xi\in N^*/\OO} \vf_\xi,\qquad
  Ker(\vf_{N\le\OO}) = \bigoplus_{\xi\in N^*/\OO} Ker(\vf_\xi),$$
and the Theorem equivalently means that the inclusions $\sR(N\le\OO)_\xi\subseteq Ker(\vf_\xi)$ are equalities for all $\xi\in N^*/\OO$.\\
In particular for $\xi\equiv 1$ we have $\vf_{\xi\equiv 1}=\vf_{\OO/N}:R_+(\le\OO/N)\rightarrow R(\OO/N)$ and $\sR(N\le\OO)_{\xi\equiv 1}=\sR(\le\OO/N),$ such that $\sR(\le\OO/N) = Ker(\vf_{\OO/N})$ is a necessary condition.}\\

{\bf Proof:} $[H,\chi]\in R_+(N\le\OO)_\xi$ means $[H,\chi]\ge \xi=[N,\mu]$, and this implies
$$   \vf_{N\le\OO}([H,\chi])|_N = Ind_H^\OO(\chi)|_N = \textrm{a multiple of}\;[N,\mu],\quad\textrm{hence:}\quad \vf_{N\le\OO}([H,\chi])\in R(\ol\OO)_\xi.$$
{\bf qed.}

Now we turn to the maps $\vf_\xi$:\\
{\bf 2.9 Lemma:} (see {\bf [L],chapt.17, [D],1.13.2})\\
{\it Fix $\mu\in N^*$ a representative of the conjugacy class $\xi\in N^*/\OO$ and let $\OO_\mu=Stab_\OO(\mu)\ge N$ be the stabilizer of $\mu$ under conjugation. Then we obtain the commutative diagram
$$\xymatrix{
Ker(\vf_\mu)\ar[d]^{Ind} \ar[r]^{\subset} & R_+(N\le \OO_\mu)_\mu \ar[d]^{Ind} \ar[r]^{\vf_\mu}  &  R(\ol\OO_\mu)_\mu \ar[d]^{Ind}\\
Ker(\vf_\xi)\ar[r]^{\subset} & R_+(N\le\OO)_\xi \ar[r]^{\vf_\xi} &  R(\ol\OO)_\xi },$$
where $\ol\OO_\mu = \OO_\mu/_{[N,N]}$ and
the vertical induction-arrows are bijections, in particular
$$  Ker(\vf_\xi) = Ind_{\OO_\mu}^\OO(Ker(\vf_\mu)).$$
Replacing $\OO$ by $\OO_\mu$ and $\vf_\xi$ by $\vf_\mu$ this brings us down to the case where $\xi=\{\mu\}$ is a single element $\mu\in (N^*)^\OO = (N/[N,\OO])^*.$}\\

{\bf Proof:} We consider the group extension $\ol N\hookrightarrow \ol\OO\twoheadrightarrow \OO/N$ with abelian kernel and let $\rho\in R(\ol\OO)$ be an irreducible representation that means a generator. Then it is well known that $\rho$ determines an orbit $\xi(\rho)\in N^*/\OO$ and, for a fixed $\mu\in \xi(\rho)$, an irreducible representation $\rho_\mu\in R(\ol\OO_\mu)$ such that $Ind_{\OO_\mu}^\OO(\rho_\mu) =\rho.$ Therefore the right hand vertical is a bijection.\\
As to the left hand vertical consider $\sum_{i\in I} n_i[H_i,\chi_i]_\OO \;\in Ker(\vf_\xi)\subseteq Ker(\vf_{N\le\OO}),$ that means\\ {\bf (*)}: $\quad \sum_{i\in I} n_i Ind_{H_i}^\OO(\chi_i) \cong 0,\quad$ where we may assume for all $i$ that $\chi_i|_N =\mu$ is a fixed representative from $\xi.$
Thus $\mu\in N^*$ extends to $\chi_i\in H_i^*$, hence $H_i\le\OO_\mu$ and:
$$  Ind_{H_i}^\OO(\chi_i)|_{\OO_\mu} = \sum_{s\in\;H_i\backslash \OO/\OO_\mu} Ind_{H_i^s\cap \OO_\mu}^{\OO_\mu}(\chi_i^s).$$
The direct sum on the right contains the term $Ind_{H_i}^{\OO_\mu}(\chi_i)$ for $s\in H_i\OO_\mu$ and additional components for $H_i s\OO_\mu \ne H_i\OO_\mu.$ And restricting further to the normal subgroup $N \le H_i^s\cap\OO_\mu$ we see that those other components are disjoint from $Ind_{H_i}^{\OO_\mu}(\chi_i)$ because $s\notin \OO_\mu.$ Therefore from {\bf (*)} we may conclude:
$$ \sum_{i\in I} n_i Ind_{H_i}^{\OO_\mu}(\chi_i) \cong 0, \qquad \sum_{i\in I} n_i[H_i,\chi_i]_{\OO_\mu} \in Ker(\vf_\mu),$$
which gives us a uniquely determined preimage in $Ker(\vf_\mu).$\\
The converse map $R_+(N\le\OO)_\xi\rightarrow R_+(N\le\OO_\mu)_\mu$ is obtained by selecting from the conjugacy class $[H,\chi]_\OO$ all members $\chi'$ such that $\chi'|_N =\mu.$ Then $\OO$-conjugation shrinks to $\OO_\mu$-conjugation.
{\bf qed.}\\
As a reformulation for a commutative $N=C$ we note:\\
{\bf 2.10 Reformulation:} {\it Let $C\le\OO$ be an abelian normal subgroup and assume $\rho=\sum_i n_i[H_i,\chi_i]\in Ker(\vf_{C\le\OO})\subset R_+(C\le\OO).$ Then for any fixed $\mu\in C^*$ we will have:
 $$  \rho_\mu:=\sum_{\chi_i;\;\chi_i|_C=\mu} n_i [H_i,\chi_i]_{\OO_\mu} \in  Ker(\vf_{C/_{[C,\OO_\mu]} \le \OO_\mu/_{[C,\OO_\mu]}})\subseteq Ker(\vf_{C\le\OO_\mu}),$$
 and fixing one representative $\mu$ in each conjugacy class $\xi\in C^*/\OO$ we recover $\rho=\sum_\mu Ind_{\OO_\mu}^\OO(\rho_\mu).$}\\

{\bf 2.11 Corollary:} {\it {\bf (i)} It is enough to prove $\sR(N\le\OO)_\xi = Ker(\vf_\xi)$ in the case where $\xi=\mu\in (N^*)^\OO= (N/[N,\OO])^*$.\\
{\bf (ii)} It is enough to prove $\sR(N\le\OO) = Ker(\vf_{N\le\OO})$ in the case where the open normal subgroup $N$ is in the center of $\OO.$}\\

{\bf Proof:} {\bf (i)}
In terms of Lemma {\bf 2.9} we will have $ Ker(\vf_\mu)=\sR(N\le\OO_\mu)_\mu$ {\bf if we assume} equality for cases $\xi=\mu,$ and therefore {\bf 2.9} and {\bf 2.6.3} yield
$$  Ker(\vf_\xi) =_{2.9} Ind_{\OO_\mu}^\OO(Ker(\vf_\mu)) = Ind_{\OO_\mu}^\OO(\sR(N\le\OO_\mu)_\mu)\subseteq \sR(N\le\OO)_\xi .$$
{\bf (ii)} According to {\bf (i)} we may in {\bf (6)} and {\bf (7)} restrict to the partial sums for $N/_{[N,\OO]} \le \OO/_{[N,\OO]}.$  {\bf qed.}\\

{\large\bf F)} Proving {\bf 2.7} in the case of a central extension.\\
{\bf 2.12 Proposition:} {\it If $(\OO,Z)$ is a pair consisting of a solvable profinite group $\OO$ and an open subgroup $Z\le\OO$ which is  central in $\OO,$ then
$\sR(Z\le \OO)\subseteq Ker(\vf_{Z\le\OO})$ is always an equality.}\\

We begin from the case $Z=\{1\}$ where {\bf $\OO$ is finite} and prove in three steps:\\
{\bf 2.12.(i)} (Revised version of {\bf [D],1.13,1.14}): {\it We have $\sR(\le\OO) =Ker(\vf_\OO)$ if $\OO$ is a finite solvable group.}\\

{\bf Step 1:} The case where $\OO$ is commutative.\\
If $\OO$ is commutative then the projector $\Phi_\OO:R_+(\le\OO)\twoheadrightarrow R_+(\OO\le\OO)$ is well defined. Moreover $\rho\in Ker(\vf_\OO)$ implies $\Phi_\OO(\rho)\in Ker(\vf_{\OO\le\OO})=\{0\},$ because $\Phi_{\OO\le\OO}:R_+(\OO\le\OO)\rightarrow R(\OO)$ is the identity. Therefore
$Ker(\vf_\OO)= Ker(\Phi_\OO) = Im(Id -\Phi_\OO)$ and we are left to show that
$$  [H,\chi] -\Phi_\OO([H,\chi])\in \sR(\le\OO)$$
for all generators $[H,\chi]\in R_+(\le\OO).$ We do this by induction over $(\OO:H).$ If $(\OO:H)=\ell$ is a prime then this is true by definition {\bf 2.6}. Otherwise consider $H<H'<\OO$ where $(H':H)=\ell$ is again a prime. Then we have
$$  [H,\chi]-\Phi_\OO([H,\chi])=Ind_{H'}^\OO\left([H,\chi]_{H'}-\Phi_{H'}^{H'}([H,\chi]_{H'})\right) +\Phi_{H'}([H,\chi])-\Phi_\OO([H,\chi]),$$
where the first term is in $\sR({\bf I},\le\OO)$, and
$$ \Phi_{H'}([H,\chi])-\Phi_\OO([H,\chi])=\Phi_{H'}([H,\chi])-\Phi_\OO\circ\Phi_{H'}([H,\chi])$$
rewrites in terms of $[H'',\chi] - \Phi_\OO([H'',\chi])$ for $H''\ge H',$ hence it is in $\sR(\le\OO)$ by induction assumption. {\bf qed.}\\

{\bf Step 2: ([L], chapt.18)} {\it The case where $\OO$ is a finite nilpotent group.}\\

In this case $\OO$ has a {\bf nontrivial center} $Z\ne\{1\},$ and we argue by induction over $(\OO:Z),$ where $(\OO:Z)=1$ has been done in step 1. Thus let $Z$ be a proper subgroup.\\  Then we consider $\rho\in Ker(\vf_\OO),$ and write
$  \rho = \Phi_Z(\rho) + (\rho -\Phi_Z(\rho)).$ As to the second summand it is enough to consider
$$ [H,\chi]-\Phi_Z([H,\chi])= Ind_{HZ}^\OO \left([H,\chi]_{HZ} - \Phi_Z^{HZ}([H,\chi])\right).$$
We have $[HZ,HZ]=[H,H]$ that means $A:=HZ/[H,H]$ is abelian and $\chi\in H^*$ extends to a character of $A.$ Therefore modulo $[H,H]$ we obtain $[H,\chi]_{HZ} - \Phi_Z^{HZ}([H,\chi])\in Ker(\vf_A)= \sR(\le A)=\sR({\bf I},\le A)$  and via inflation and induction this yields $[H,\chi]-\Phi_Z([H,\chi])\in \sR({\bf I},\le\OO).$\\
Thus we are left to show that $\rho\in Ker(\vf_\OO)$ implies $\rho': =\Phi_Z(\rho)\in \sR(\le\OO).$ To this end we choose $\OO> C> Z$
such that $C$ is a commutative normal subgroup of $\OO,$ and $(C:Z)=\ell$ is a prime, and $C/Z$ is in the center of $\OO/Z.$ Such a $C$ exists because $\OO$ is nilpotent. (Consider $Z_2/Z =$ center of $\OO/Z,$ and the commutator map $Z_2/Z\wedge Z_2/Z \rightarrow Z.$ Then $C/Z$ of prime order inside $Z_2/Z$ will do it).  Then we consider
$$ \rho' = \Phi_C(\rho') + (\rho' - \Phi_C(\rho'))\quad\textrm{for}\;\rho'=\Phi_Z(\rho)\in Ker(\vf_{Z\le\OO}).$$
Here we use {\bf 2.10} for $\Phi_C(\rho')=\Phi_C(\rho)\in Ker(\vf_{C\le\OO})$ (instead of $\rho$ there). This gives us
$\Phi_C(\rho)=\sum_\mu Ind_{\OO_\mu}^\OO(\Phi_C(\rho)_\mu)$ for
$\Phi_C(\rho)_\mu\in Ker(\vf_{\ol C\le \ol\OO_\mu})$ where $\ol C= C/_{[C,\OO_\mu]}$ is in the center of $\ol\OO_\mu=\OO_\mu/_{[C,\OO_\mu]}$ with quotient $(\ol\OO_\mu:\ol C)= (\OO_\mu :C) < (\OO:Z).$\\
Therefore for all $\mu$ we have $\Phi_C(\rho)_\mu\in Ker(\vf_{\ol\OO_\mu})=\sR(\le \ol\OO_\mu)$ by induction assumption.
By inflation we are then in $\sR(\le\OO_\mu)$ and induction brings us to $\Phi_C(\rho')=\Phi_C(\rho)\in\sR(\le\OO).$\\
Finally we need to handle the difference $\rho' - \Phi_C(\rho')\in R_+(Z\le\OO)$, and here it is enough to see:\\
{\bf 2.13 Lemma:} {\it For generators $[H,\chi]\in R_+(Z\le\OO)$ we always have:
$$  [H,\chi] - \Phi_C([H,\chi]) \in \sR(Z\le\OO).$$}\\
{\bf Proof:} Obviously:
$[H,\chi] - \Phi_C([H,\chi]) = Ind_{HC}^\OO\left([H,\chi]_{HC} - \Phi_C^{HC}([H,\chi]_{HC}\right),$
such that it is enough to see that $$[H,\chi]_{HC} - \Phi_C^{HC}([H,\chi]_{HC})\;\in \sR(Z\le HC).$$
If $HC\ne\OO$ then $(HC:Z)< (\Omega:Z),$ hence $ Ker(\vf_{HC})= \sR(\le HC)$ by induction assumption.\\
{\bf Now we assume $HC=\OO.$} {\it Then $H$ is a normal subgroup of $\OO$, of index $(\OO:H)\mid (C:Z)=\ell$.} Since the case $\OO=H\ge C$ is trivial, we are left with the case  where $H<\OO$ is normal of index $\ell$ and $H\cap C =Z$ which means $C/Z\stackrel{\sim}{\rightarrow}\OO/H.$  In particular $H$ is normalized by $C$ such that due to {\bf 2.2(ii)} and {\bf 2.3} we have
$$ \Phi_C([H,\chi]) = \sum_{\nu\in (C/Z)^*/\sim} [H'C,\chi\mu_0]\nu',$$
where $H'\;=H_\mu$ for all $\mu\in S(\chi),$ is the projective kernel of $Ind_H^{HC}(\chi)=Ind_H^\OO(\chi).$ And $\nu'\in (\OO/H)^*$ is the natural extension of $\nu\in (C/Z)^*,$ where the $\nu$ are running over a system of representatives for the cokernel of\\
{\bf (11)}
$$   H/H'\hookrightarrow (C/Z)^*,\qquad  h\mapsto \mu_0^{h-1}.$$
Since $(C/Z)^*$ is of order $\ell,$ a prime, we have only two cases to deal with:\\
{\bf Case 1:} The map {\bf (11)} is surjective, that means $(H:H')=\ell,$ and the cokernel of {\bf (11)} is trivial, hence can be represented by the trivial character $\nu\equiv 1$. Then we obtain:
$$  \Phi_C([H,\chi]) = [H'C,\chi\mu_0],\qquad Ind_H^\OO(\chi) = Ind_{H'C}^\OO(\chi\mu_0),$$
where $(\OO:H'C) = \ell = (\OO:H)$ and $H\cap H'C =H'(H\cap C) = H'Z =H'$ because $Z=H\cap C$ has trivial inner action on $C^*,$ hence $Z\le H'=Stab_H(\mu)$ and therefore $Z\le H'\cap C\le H\cap C =Z.$ Thus we see that $H$ and $H'C$ are both normal of index $\ell$ in $\OO$ and $H\cap H'C=H'(H\cap C) =H'.$ This gives us:
\begin{itemize}
 \item $\OO/H'\cong \Z/\ell\times \Z/\ell,$ in particular $[\OO,\OO]\le H',$ hence $[[\OO,\OO],\OO]\le [H',\OO]$ is in the kernel of our representation (because $H'$ is the projective kernel),
 \item $\chi|_{H'} = (\chi\mu_0)|_{H'} =:\eta\quad$ is a character of $H'/[H',\OO]$ which is non-trivial on $[\OO,\OO]/[H',\OO]$ because otherwise $\chi$ could be extended to a character $\wt\chi\in\OO^*$ and $Ind_H^\OO(\chi)=\wt\chi\otimes Ind_H^\OO(1),$ hence the projective kernel would be $H'=H.$
\end{itemize}
Altogether we see that in {\bf case 1} we obtain:
$$  [H,\chi] - \Phi_C([H,\chi])= [H,\chi]- [H'C,\chi\mu_0]= \rho(H',\OO,\eta) \in\sR({\bf II},Z\le\OO).$$
{\bf Case 2:} $H=H'$ that means the image of {\bf (11)} is trivial, hence $(C/Z)^*$ is the cokernel and $C/Z\stackrel{\sim}{\rightarrow}\OO/H=\OO/H'.$
Therefore the character $\chi\mu_0\in\OO^*$ is an extension of $\chi\in H^*$,
$$  \Phi_C([H,\chi]) =  \sum_{\nu\in (C/Z)^*} (\OO,\chi\mu_0\nu'),$$
$$ [H,\chi] -\Phi_C([H,\chi]) = \rho(H,\OO,\chi\mu_0)  \in \sR({\bf I},H\le\OO)\subseteq \sR({\bf I},Z\le\OO).$$
This ends the proof of {\bf 2.13} and of {\bf Step 2}.\\

{\bf Step 3:} {\it For any finite solvable group $\OO$ we will have $\sR(\le\OO) = Ker(\vf_\OO).$}\\

We prove this by induction over $\#\OO.$ Thus for proper subgroups $H<\OO$ we will have $\sR(\le H)=Ker(\vf_H)$ which is an ideal in $R_+(\le H)$, and this will imply that:\\
{\it Also $\sR(\le\OO)$ is an ideal in $R_+(\le\OO).\quad$} Indeed, let $\rho\in\sR(\le\OO)$ and $[H,\chi]\in R_+(\le\OO)$ a generator. If $H=\OO,$ then $\rho\cdot (\OO,\chi)$ belongs to $\sR(\le\OO)$ by definition {\bf 2.6}. If $H\ne \OO$, then we may use Proposition {\bf 2.1}, hence
$$ \rho\cdot [H,\chi] = Ind_H^\OO(Res_H^\OO(\rho)\cdot(H,\chi)).$$
Moreover, applying {\bf (F2)} we obtain $\vf_H(Res_H^\OO\,\rho) = Res_H^\OO(\vf_\OO(\rho)) = 0.$
Therefore   $Res_H^\OO(\rho)\in Ker(\vf_H)=\sR(\le H)$ and then obviously
$Ind_H^\OO(Res_H^\OO\;\rho\cdot(H,\chi)) \in \sR(\le\OO).$  {\bf qed.}\\

Now let $Z$ be the center of $\OO.$ By Step 2 we can assume that $\OO$ is {\bf not nilpotent} hence it might happen that $Z=\{1\}.$\\
If $Z\ne\{1\}$, we use {\bf Brauer 1} for $\OO/Z$:  There are {\bf nilpotent} subgroups $H_i\supseteq Z$ and characters $\chi_i$ of $H_i/Z$ such that
$$  1_{\OO/Z}  = \sum_i n_i\cdot Ind_{H_i/Z}^{\OO/Z}(\chi_i)\quad\textrm{with}\quad n_i\in \Z.$$
\begin{itemize}
 \item Together with $\OO$ also $\OO/Z$ is not nilpotent, but the $H_i/Z$ are nilpotent, hence they are {\bf proper subgroups} of $\OO/Z$ which ensures that we obtain here a non-trivial relation.
\end{itemize}

Since we do induction over the order of $\OO$, we may assume $Ker(\vf_{\OO/Z}) = \sR(\le\OO/Z),$ and by inflation we obtain
$$  \sigma: = 1_\OO - \sum_i n_i\cdot [H_i,\chi_i]_\OO \in \sR(\le\OO).$$ Now we use that $\sR(\le\OO)\subset R_+(\le\OO)$ is an ideal, in particular $\rho\sigma\in \sR(\le\OO)$ for all $\rho\in Ker(\vf_\OO).$ And by Proposition {\bf 2.1} we have:
$$  \rho\sigma = \rho - \sum_i n_i Ind_{H_i}^\OO(Res_{H_i}^\OO(\rho)\cdot (H_i,\chi_i)).$$
Moreover $Res_{H_i}^\OO(\rho)\in Ker(\vf_{H_i})=\sR(\le H_i)$ as we have seen already, and therefore
$$ \rho = \rho\sigma +\sum_i n_i Ind_{H_i}^\OO(Res_{H_i}^\OO(\rho)\cdot (H_i,\chi_i))\in\sR(\le\OO),$$
{\bf which proves our assertion in the case $Z\ne\{1\}.$}\\
If $Z=\{1\}$ then let $C$ be a minimal commutative normal subgroup of $\OO$ which exists because $\OO$ is solvable. And for $\rho\in Ker(\vf_\OO)$ write again:
$$  \rho= \Phi_C(\rho) + (\rho - \Phi_C(\rho)),$$
where $\Phi_C(\rho)\in Ker(\vf_{C\le\OO})\subset R_+(C\le\OO)$. Thus by {\bf 2.10} we obtain
$$\Phi_C(\rho) =\sum_\mu Ind_{\OO_\mu}^\OO(\Phi_C(\rho)_\mu),$$
where $\mu\in C^*$ runs over a system of representatives for the conjugacy classes $C^*/\OO.$\\
If $\mu\equiv 1$ then $\Phi_C(\rho)_1\in Ker(\vf_{\OO/C}) = \sR(\le\OO/C).$\\
If $\mu\not\equiv 1$ then $\OO_\mu\ne\OO$ because otherwise $\mu\in(C/[C,\OO])^*$ which contradicts with $[C,\OO]=C$ (because $C$ is minimal normal and the center of $\OO$ is trivial). Therefore $\Phi_C(\rho)_\mu\in Ker(\vf_{\OO_\mu})=\sR(\le\OO_\mu).$
Altogether this yields $\Phi_C(\rho)\in \sR(\le\OO).$\\
Thus we are left with the following variant of {\bf 2.13}:\\
{\bf 2.14 Lemma:} {\it For generators $[H,\chi]\in R_+(\le\OO)$ we always have:
$$  [H,\chi] - \Phi_C([H,\chi]) \in \sR(\le\OO).$$}\\
Since we argue by induction on $\#\OO,$ we may as in {\bf 2.13} restrict to the particular case where  $HC =\OO.$\\
Thus assume $HC=\OO$, hence $H\cap C$ is normal in $\OO$ and therefore $H\cap C =\{1\}$ or $=C,$ because $C$ was minimal. The case $C\le H$ is trivial such that we have to consider only the case where $H\cap C=\{1\}$, {\bf hence $\OO=H\ltimes C$ semidirect. Then $C$ minimal normal in $\OO$, implies that $H<\OO$ is a maximal subgroup} (see Lemma 3 in Appendix 2). Moreover $\chi\in H^*$ lifts to $\wt\chi\in(\OO/C)^*$ and therefore\\
{\bf (12)}
$$[H,\chi] - \Phi_C([H,\chi]) =(\OO,\wt\chi)\cdot\left([H,1] - \Phi_C([H,1])\right).$$
{\bf Suppose there is a nontrivial subgroup $H_1\le H$ which is normal in $\OO.$} Then $[H_1,C]\subseteq H_1\cap C\subseteq H\cap C =\{1\}$. Therefore
$H_1\le H_\mu$ for all $\mu$ such that under\\
$\OO\twoheadrightarrow \OO/H_1=(H/H_1)\ltimes C$ we obtain
$$[H,1] - \Phi_C([H,1]) = Inf_{\OO/H_1}^\OO\left([H/H_1,1] -
\Phi_C([H/H_1,1])\right)\quad\in \sR(\le\OO),$$
because $Ker(\vf_{\OO/H_1})=\sR(\le \OO/H_1)$ by induction assumption.\\
Thus we are left with the case where $H<\OO$ is maximal and $H$ does not contain normal subgroups of $\OO=H\ltimes C$. But then $\OO> H$ is a type-III-group (as we see from Lemma 3 of Appendix 2) and by definition {\bf 2.6} we have then {\bf (12)} in $\sR({\bf III},\le\OO)$\\
{\bf This ends the proof of 2.12.(i)}.\\

To complete the proof of {\bf 2.12}, we are left with the assertion:\\
{\bf 2.12.(ii):} {\it $\sR(Z\le\OO) = Ker(\vf_{Z\le\OO})$ in all cases where $Z$ is a {\bf non-trivial} open subgroup which is contained in the center of the solvable profinite group $\OO.$}\\

{\bf Proof:} We prove this by induction over $(\OO:Z).$ For $Z=\OO$ the group is abelian and we obtain an isomorphism
$$  \vf_{\OO\le\OO}:R_+(\OO\le\OO)\cong R(\OO),\qquad Ker(\vf_{\OO\le\OO}) =\{0\},$$
such that there is nothing to show.\\
Now let $Z\ne\OO$ and {\bf consider} $\rho \in R_+(Z\le\OO)\cap Ker(\vf_\OO) = Ker(\vf_{Z\le\OO}).$
We are going to prove $\rho\in \sR(Z\le\OO),$ where we distinguish three cases:\\
{\bf a)} If $Z$ is not the whole center of $\OO$ then we find $C$ in the center of $\OO$ such that $(C:Z)=\ell$ is a prime.\\
{\bf b)} If $Z$ is the center of $\OO$, and if $\OO/Z$ has nontrivial center $Z_2/Z\ne\{1\},$ then we will consider $Z<C\le Z_2$ such that $(C:Z)=\ell$ is a prime and $C<\OO$ is an abelian normal subgroup (as in step 2 above).\\
{\bf c)} If $Z$ is the center of $\OO$ and moreover $\OO/Z$ has trivial center, then we will reduce to $\sR(\le G)= Ker(\vf_G)$ for the finite quotient $G=\OO/Z.$\\

In cases {\bf a)} and {\bf b)} resp. we write $\rho = \Phi_C(\rho) + (\rho - \Phi_C(\rho)),$  and we are going to see that our assumption on  $\rho$ implies $\Phi_C(\rho)\in \sR(Z\le\OO)$ and moreover $[H,\chi]-\Phi_C([H,\chi])\in \sR(Z\le\OO)$ for all generators $[H,\chi]\in R_+(Z\le\OO).$\\
{\bf In case a)} we obtain:
$\Phi_C(\rho)\in Ker(\vf_{C\le\OO}) = \sR(C\le\OO)$ by induction assumption (because $C$ is central and $(\OO:C)<(\OO:Z)$),  and obviously $\sR(C\le\OO)\subseteq\sR(Z\le\OO).$\\
On the other hand:  $[H,\chi]-\Phi_C([H,\chi])= 0$  if $H\ge C$. Otherwise we have $C/Z\cong HC/H$ because $C/Z$ is of prime order, and, since $C$ is central in $\OO,$ Corollary {\bf 2.3 (ii)} yields
$$  [H,\chi]-\Phi_C([H,\chi])= [H,\chi] -\sum_{\chi'} [HC,\chi'] = \rho(H,HC,\wt\chi)\in \sR({\bf I},Z\le\OO),$$
where $\chi'\in (HC)^*$ runs over the extensions of $\chi$ and where $\wt\chi$ is any of these extensions.\\

{\bf In case b)}, as in step 2 above, we take $Z<C\le Z_2$ such that $(C:Z)=\ell,$ and then we can argue in the same way:\\
Indeed due to {\bf 2.10} we have $\Phi_C(\rho)=\sum_{\mu\in C^*/\OO} Ind_{\OO_\mu}^\OO(\Phi_C(\rho)_\mu)$, where
$$ \Phi_C(\rho)_\mu \in Ker(\vf_{C/_{[C,\OO_\mu]} \le \OO_\mu/_{[C,\OO_\mu]}}).$$
Since $C/_{[C,\OO_\mu]} \hookrightarrow \OO_\mu/_{[C,\OO_\mu]}\twoheadrightarrow \OO_\mu/C$ is a central extension with $(\OO_\mu:C)< (\OO:Z),$ we can apply the induction assumption and conclude that $\Phi_C(\rho)\in\sR(C\le\OO)\subseteq \sR(Z\le\OO).$ And concerning   $[H,\chi] -\Phi_C([H,\chi])$ for $[H,\chi]\in R_+(Z\le\OO)$ we can use Lemma {\bf 2.13}.\\

{\bf In case c)} the quotient $\OO/Z$ has trivial center, hence it is not nilpotent. Therefore we can argue precisely the same way as in the case $Z\ne\{1\}$ of step 3 above:\\
{\bf 2.15 Lemma:} {\it Let $Z\ne\{1\}<\OO$ be a central open subgroup such that $\OO/Z$ is not nilpotent.  If $Ker(\vf_{Z\le\OO'}) =\sR(Z\le\OO')$ for all $Z\le\OO'<\OO$ and if $Ker(\vf_{\OO/Z})=\sR(\le\OO/Z),$ then $\sR(Z\le\OO)$ is an ideal in $R_+(Z\le\OO),$ and
$$  Ker(\vf_{Z\le\OO}) =\sR(Z\le\OO)$$
will follow.}\\
The assumptions hold, because we do induction on $(\OO:Z)$ and because $\OO/Z$ is finite resp. {\bf This finishes the proof of 2.12.(ii) and of Theorem 2.7}

.\\

{\bf\large 3. A criterion for the extendibility of functions}\\

Now we come to the goal of this paper. We want to prove the following criterion {\bf 3.1} for the extendibility of functions, in the sense of {\bf Definition 1.1}. The criterion makes it more precise how to use {\bf Theorem 2.7} in order to see that a function $\Delta$ on $R_1(\le\OO)$ with values in a multiplicative abelian group, extends to a function $\sF$ on $R(\le\OO).$ As already explained in {\bf 1.6} we are going to see by induction, that for all subquotients $N'\le\OO'$ the restriction of $\Delta$ on $R_1(N'\le\OO')$ can be extended onto $R(N'\le\OO').$ And without loss of generality we may restrict here to $\OO'=\OO,$ because all arguments are left unchanged if we replace $\OO$ by an open subgroup $\OO'.$\\

{\bf Theorem 3.1:}\\
{\it Let $\OO$ be a solvable profinite group and let $\Delta=\Delta(H,\chi)$ be a function on $R_1(\le\OO)$ with values in a multiplicative abelian group $\sA$.
Then $\Delta$ is {\bf extendible} to $\sF$ on $R(\le\OO)=\bigsqcup_{H\le\OO} R(H)$ {\bf if and only if} the following three conditions are fulfilled for all open subgroups $B\le \OO:$ \\
{\bf I. Davenport-Hasse:}  Let $K<B$ be normal of prime index $\ell$ and $\chi\in B^*.$ Then $Ind_K^B(\chi_K)= \bigoplus_{\mu\in(B/K)^*}\chi\mu,$ and we should have:\\
{\bf (1)}
$$ \Delta(K,\chi_K)\cdot\prod_{\mu\in(B/K)^*}\Delta(B,\mu) =\prod_{\mu\in(B/K)^*} \Delta(B,\chi\mu).$$

{\bf II. Heisenberg identity:} Let $Z< B$ be normal with abelian quotient $B/Z\cong \Z/\ell\times\Z/\ell,$ where $\ell$ is a prime and such that the commutator induces an isomorphism
$$ [.,.]:B/Z\wedge B/Z \cong [B,B]/[Z,B]\quad (\cong \Z/\ell),$$
and let $\eta:Z/[Z,B]\rightarrow \fC^\times$ be a character which is non-trivial on $[B,B].$ Depending on $(Z,\eta)$ consider pairs $(H,\eta^H)$ where $H> Z$ is of index $\ell$ in $B,$ hence $[H,H]\le [Z,B]$ and where $\eta^H:H/[Z,B]\rightarrow \fC^\times$ extends $\eta.$ Then $Ind_H^B(\eta^H)$ is for all such pairs the same Heisenberg representation, $\ell\cdot Ind_H^B(\eta^H)= Ind_Z^B(\eta),$   and we should have:\\
{\bf (2)}
$$  \Delta(H,\eta^H)\cdot\prod_{\mu\in(B/H)^*}\Delta(B,\mu) \;\textrm{ does not depend on the pair}\; (H,\eta^H).$$

{\bf III. Generalized Davenport-Hasse:} Let $H<B$ be a {\bf maximal} subgroup which is {\bf not normal}, such that $K:=\bigcap_{b\in B} bHb^{-1} \ne H,$ hence $H/K < B/K$ is non-degenerate of type III (see Appendix 2, Lemma 1), and let $C< B$ be the uniquely determined normal subgroup such that $HC =B$ and $H\cap C =K.$ Then for any $\chi\in B^*$ and with respect to the inner action of $H$ on $(C/K)^*$ we have $Ind_H^B(\chi_H) = \bigoplus_{\mu\in (C/K)^*/H} Ind_{H_\mu C}^B(\chi\mu'),$ and related to that the condition\\
{\bf (3)}
$$ \Delta(H,\chi_H)\cdot \prod_{\mu\in(C/K)^*/H}\Delta(H_\mu C,\mu') =\prod_{\mu\in(C/K)^*/H}\Delta(H_\mu C,\chi\mu')$$
is required, (where $\mu'$ denotes the trivial extension of $\mu$ onto the stabilizer $H_\mu C$).}\\

{\bf Remarks:} For the degenerate case $H=K,$ $C=B$ the type-III-relation turns into type-I.\\

{\bf Proof of Theorem 3.1} (Later we will take for $\Delta$ the local Hecke-Tate root number and then we will translate the requirements {\bf I.-III.} from groups to fields, relations which have to be confirmed in order to verify the existence of local Artin root numbers.)\\
{\bf Assume that $\Delta$ on $R_1(\le\OO)$ is a function which extends to $\sF$ on $R(\le\OO)=\bigsqcup_{H\le\OO} R(H).$} In the situation {\bf I} we have the relation
$$  Ind_K^B(\chi_K) = \sum_{\mu\in(B/K)^*} \chi\mu,\quad\textrm{hence}\quad \sF(Ind_K^B(\chi_K))=\prod_{\mu\in(B/K)^*} \Delta(B,\chi\mu),$$
and together with $\sF(Ind_K^B(\chi_K)) = \Delta(K,\chi_K)\cdot\lambda_K^B(\Delta)=\Delta(K,\chi_K)\cdot\prod_{\mu\in(B/K)^*}\Delta(B,\mu)$ this yields the condition {\bf (1)}.\\

In the situation {\bf II} let $(H,\eta^H)$ and $(H',\eta^{H'})$ be two isotropic pairs related to $(Z,\eta).$ Then using Lemma 2.1 for the commutative normal subgroup $C=H/_{[Z,B]} < B/_{[Z,B]}$ we obtain $Ind_{H'}^B(\eta^{H'}) \cong Ind_H^B(\eta^H),$ hence
$\sF(Ind_{H'}^B(\eta^{H'})) = \sF(Ind_{H'}^B(\eta^{H'})),$ which rewrites as {\bf (2)}:
$$  \Delta(H',\eta^{H'})\cdot\prod_{\mu'\in(B/H')^*}\Delta(B,\mu') = \Delta(H,\eta^H)\cdot\prod_{\mu\in(B/H)^*}\Delta(B,\mu).$$

Finally in the situation {\bf III}  we have $Ind_H^B(\chi_H) = \bigoplus_{\mu\in (C/K)^*/H} Ind_{H_\mu C}^B(\chi\mu'),$ which in particular will hold for the trivial character $\chi\equiv 1.$ Therefore:
$$  \Delta(H,\chi_H) =\frac{\sF(Ind_H^B(\chi_H))}{\sF(Ind_H^B(1))} = \frac{\prod_{\mu\in(C/K)^*/H}\Delta(H_\mu C,\chi\mu')\lambda_{H_\mu C}^B(\Delta)}{\prod_{\mu\in(C/K)^*/H}\Delta(H_\mu C,\mu')\lambda_{H_\mu C}^B(\Delta)},$$
which recovers the condition {\bf (3)}.\\

Now, {\bf conversely},  we assume that $\Delta$ is a function on $R_1(\le \OO)$ which satisfies the relations {\bf (1) - (3)} and we want to show that $\Delta$ is extendible to a function $\sF$ on $R(\le\OO).$ The basic idea (implicit in {\bf [L]}) is to proceed by induction as follows: Let $N'\le\OO'$ be any finite subquotient of $\OO$ and let $\Delta'$ on $R_1(N'\le\OO')$ be the restriction of $\Delta.$ Then, using Criterion {\bf 1.4}, we will check by induction on $(\OO':N')$  that $\Delta'$ extends to $\sF'$ on $R(N'\le\OO') = \bigcup_{H;\;N'\le H\le\OO'} R(H/_{[N',N']}).$ The way of induction has been sketched already in {\bf 1.6}. And we remark once more that without loss of generality we may restrict to quotients $N\le\OO$ because the same arguments will work if we replace $\OO$ by a subgroup $\OO'$ and $N$ by a normal subgroup $N'\le\OO'$. Then we can go to $(\OO':N')\rightarrow\infty,$ in particular $(\OO:N)\rightarrow\infty,$ and this will prove the Theorem.
Thus {\it we restrict to prove that the relations {\bf (1) - (3)} imply extendibility of $\Delta_N$ from $R_1(N\le\OO)$ onto $\sF_N$ of $R(N\le \OO)=\bigcup_{H;\;N\le H\le\OO} R(H/_{[N,N]}),$ {\bf if we assume this for all proper subquotients $\OO'/N'$ of $\OO/N.$}} As already explained in the introduction we will argue in  three steps (where concerning steps 1 and 2 we follow {\bf [L], chapt.16}, and then step 3 will cover {\bf [L], chap.19}):\\

{\bf Step 1:} Definition of $U\in\sU_N(\OO)\mapsto \lambda_U^\OO(\Delta).$\\
{\bf Fixing a minimal abelian normal subgroup} $C\le G=\OO/N$ the {\bf Lemma 2.2} gives us the isomorphism:\\
{\bf (4)}
$$ Ind_U^\OO(1)= Ind_H^G(1) \cong \sum_{\mu\in T} Ind_{H_\mu C}^G(\mu'),$$
where $H=U/N,$ and $T$ is a set of representatives for $(C/H\cap C)^*/H,$ and we use this for the\\

{\bf Definition 3.2:} Let $U\in \sU_N(\OO)$ be subgroups containing $N$ and let $\Delta_0$ on $R_1(\le G)$ be the restriction of $\Delta.$  Then, using the data of {\bf (4)}, define:
$$  \lambda_U^\OO(\Delta):=\lambda_H^G(\Delta_0):=\prod_{\mu\in T} \left(\Delta_0(H_\mu C,\mu')\cdot \lambda_{H_\mu C/C}^{G/C}(\Delta_0)\right),$$
where $\lambda_{H_\mu C/C}^{G/C}(\Delta_0)$ has been defined by induction.
\begin{itemize}
\item A priori $U\in\sU_N(\OO)\mapsto \lambda_U^\OO(\Delta)$ depends on $C\le G$, but we are going to see that it can be used to extend the function $\Delta_N =\Delta|_{R_1(N\le\OO)},$ and therefore it is actually independent.
 \item In the particular case where $H\ge C,$ Lemma 2.2 comes down to the identity $Ind_H^G(1)=Ind_H^G(1),$ where $T$ consists only of the trivial character of $C,$ and therefore our definition comes down to
$\lambda_H^G(\Delta_0):=\Delta_0(H,1)\cdot\lambda_{H/C}^{G/C}(\Delta_0)=\lambda_{H/C}^{G/C}(\Delta_0).$
\item If $G=\OO/N$ is cyclic of prime order $\ell$ then the only possibility is $C=G$ and therefore: $\lambda_N^\OO(\Delta)=\lambda_1^G(\Delta_0) = \prod_{\mu\in G^*} \Delta_0(G,\mu).$
\end{itemize}

{\bf Lemma 3.3:} {\it Keep $\OO\ge U\ge N$ and $G\ge H$ as in {\bf 3.2}.\\
{\bf (i)} The definition of $\lambda_U^\OO(\Delta)=\lambda_H^G(\Delta_0)$ is independent of the choice of the set $T$ of representatives for $(C/H\cap C)^*/H.$\\
{\bf (ii)} $\lambda_{H^g}^G(\Delta_0) = \lambda_H^G(\Delta_0)$ for all $g\in G,$ which means $\lambda_{U^g}^\OO(\Delta)=\Delta_U^\OO(\Delta)$ for all $g\in\OO$ and $U\ge N.$\\
{\bf (iii)} If $N'>N$ is a larger normal subgroup of $\OO$ such that $\OO/N'= G/N_G'$ is a quotient of $\OO/N =G,$ and if $U\in\sU_{N'}(\OO)$ such that $U/N'=H/N_G'$ is a quotient of $U/N =H,$ then
%If $N_G'\le H$ is a nontrivial normal subgroup of $G$ which sits in $H=U/N,$ and if $N<N'\twoheadrightarrow N_G'$ is the preimage in $\OO,$ that means $\OO/N'=G/N_G'$ and $U\in\sU_{N'}(\OO),$ then
$$  \lambda_H^G(\Delta_0) = \lambda_{H/N_G'}^{G/N_G'}(\Delta_0),\qquad \lambda_U^\OO(\Delta) = \lambda_U^\OO(\Delta_{N'}),$$
where $\Delta_{N'}$ is the restriction of $\Delta$ to $R_1(N'\le\OO)$ and the right sides are given by induction assumption.}\\
{\bf Remark:} The Lemma says that $H\in \sU(G) \mapsto \lambda_H^G(\Delta_0)$ is well defined and has the properties {\bf 1.(6)} and {\bf 1.(8)}.\\

{\bf Proof:} {\bf (i)} If $h\in H$ then $H_{\mu^h}C = h^{-1}(H_\mu C)h,$ and $(\mu^h)' = (\mu')^h$ is the extension of $\mu^h.$ Hence by {\bf 1.(2)}:
$$  \Delta_0(H_{\mu^h}C, (\mu^h)') = \Delta_0(h^{-1}(H_\mu C)h, (\mu')^h) =\Delta_0(H_\mu C,\mu').$$
And
$$ \lambda_{H_{\mu^h}C/C}^{G/C}(\Delta_0) = \lambda_{h^{-1}(H_\mu C/C)h}^{G/C}(\Delta_0) = \lambda_{H_\mu C/C}^{G/C}(\Delta_0),$$
because  we may apply {\bf 1.(6)} for the proper quotient $G/C.$ Therefore the definition 3.2 is left unchanged if we replace the representatives $\mu\in T$ by other representatives $\mu^h.$\\
{\bf (ii)} is the same argument as in {\bf (i)}. We have to compare {\bf (4)} with the corresponding formula for $U^g=g^{-1}Ug\twoheadrightarrow H^g=g^{-1}Hg.$ Then we obtain $T^g=(C/H^g\cap C)^*/H^g$ with representatives $\mu^g$, and $(H^g)_{\mu^g} = g^{-1}Hg,$ hence $(H^g)_{\mu^g}\cdot C= g^{-1}(H_\mu C)g$ and $(\mu^g)'=(\mu')^g,$ which implies:
$$  \lambda_{H^g}^G(\Delta_0) =\prod_{\mu\in T} \Delta_0(g^{-1}(H_\mu C)g, (\mu')^g)\cdot\lambda_{g^{-1}(H_\mu C/C)g}^{G/C}(\Delta_0) = \lambda_H^G(\Delta_0),$$
because of {\bf 1.(2)} and {\bf 1.(6)} for $G/C.$\\
{\bf (iii)} Since $C$ is minimal normal in $G$ the only possibilities are:\\
{\bf a)} $N_G'\cap C=C$ that means $H\ge N_G'\ge C$, which implies $\lambda_H^G(\Delta_0) = \lambda_{H/C}^{G/C}(\Delta_0) =\lambda_{H/N_G'}^{G/N_G'}(\Delta_0),$ where the second equality follows because $N_G'/C\le G/C$ is a normal subgroup, and for the quotient $G/C$ we assume property {\bf 1.(8)} by induction.\\
{\bf b)} $N_G'\cap C=\{1\}$ that means $C\hookrightarrow G/N_G'$ is a commutative normal subgroup. Then: $[N_G',C]\le N_G'\cap C=\{1\}$ implies that $N_G'$ and $C$ commute and therefore the inner action of $N_G'$ on $C^*$ is trivial. Thus we have $N_G'\le H_\mu$ for all $\mu\in (C/H\cap C)^*,$ and the relation {\bf (4)} rewrites as:\\
{\bf (5)}
$$ Ind_{H/N_G'}^{G/N_G'}(1) \cong \sum_{\mu\in T} Ind_{H_\mu C/N_G'}^{G/N_G'}(\mu').$$
And by assumption  $\Delta_0$ restricted to $R_1(\le G/N_G')$ extends to $\sF_0$ on $R(\le G/N_G'),$ hence:
$$ \lambda_{H/N_G'}^{G/N_G'}(\Delta_0) =\sF_0(Ind_{H/N_G'}^{G/N_G'}(1)) =_{(5)} \prod_{\mu\in T} \Delta_0(H_\mu C/N_G',\;\mu')\cdot\lambda_{H_\mu C/N_G'}^{G/N_G'}(\Delta_0).$$
Now obviously $\Delta_0(H_\mu C/N_G',\;\mu') =\Delta_0(H_\mu C,\;\mu')$, and
$$ \lambda_{H_\mu C/N_G'}^{G/N_G'}(\Delta_0) =\lambda_{H_\mu C/N_G' C}^{G/N_G' C}(\Delta_0) = \lambda_{H_\mu C/C}^{G/C}(\Delta_0),$$
because $N_G' C/N_G'$ is a normal subgroup in $G/N_G'$ which is contained in $H_\mu C$, and $N_G' C/C$ is a normal subgroup in $G/C$ which is contained in $H_\mu C/C.$ And replacing $G$ by the proper quotients $G/N_G'$ and $G/C$ resp. the relation {\bf 1.(8)} holds by induction. Finally comparing this to our definition {\bf 3.2} we see that $ \lambda_{H/N_G'}^{G/N_G'}(\Delta_0) =\lambda_H^G(\Delta_0).$ {\bf qed.}\\

{\bf Step 2:} We are going to prove that $H\in \sU(G)\mapsto \lambda_H^G(\Delta_0)$ satisfies {\bf 1.(7)} too.\\
{\bf Lemma 3.4:} ($N$-Tower lemma) {\it Let $N\le \OO$ be a normal subgroup, $\Delta$ on $R_1(N\le\OO),$ and assume that
\begin{itemize}
 \item $\Delta$ is extendible from $R_1(N'\le\OO')$ onto $R(N'\le \OO')$ for all $N\le N'\le\OO'\le\OO,$ such that $(\OO':N')< (\OO:N).$
\end{itemize}
{\bf Then:}  $\quad\lambda_{H'}^{G_0}(\Delta_0) = \lambda_{H'}^H(\Delta_0)\cdot\lambda_H^{G_0}(\Delta_0)^{(H:H')}$ for all $H'\le H\le G_0\le G=\OO/N.$}\\

{\bf Proof:} If $G_0=\OO_0/N$ is $\ne G,$ we are done by assumption, because $\#G_0=(\OO_0:N) < (\OO:N).$ Thus we are left with the   assertion\\
{\bf (6)}
$$\lambda_{H'}^{G}(\Delta_0) = \lambda_{H'}^H(\Delta_0)\cdot\lambda_H^{G}(\Delta_0)^{(H:H')}.$$
At first we reduce this to the case where $H$ is a complement of the fixed abelian normal subgroup $C< G.$ Consider\\
{\bf (7)}
$$  Ind_{H'}^{HC}(1) \cong \sum_{\mu \in S/H'} Ind_{H_\mu'C}^{HC}(\mu'),\qquad (HC:H') =\sum_\mu (HC: H_\mu'C). $$
By induction we have\\
{\bf (8)}
$$  \lambda_{H_\mu'C}^G(\Delta_0) = \lambda_{H_\mu'C}^{HC}(\Delta_0)\cdot\lambda_{HC}^G(\Delta_0)^{(HC:H_\mu'C)},$$
because here {\bf we can argue modulo $C.$} And by definition we have\\
{\bf (9)}
$$\lambda_{H'}^G(\Delta_0) =\prod_{\mu}\left(\Delta_0(H_\mu'C,\mu')\lambda_{H_\mu'C}^G(\Delta_0)\right).$$
Now substituting {\bf (8)} on the right side of {\bf (9)} and then using {\bf (7)} we conclude:\\
{\bf (10)}
$$  \lambda_{H'}^G(\Delta_0) = \lambda_{H'}^{HC}(\Delta_0)\cdot \lambda_{HC}^G(\Delta_0)^{(HC:H')}.$$
{\bf Assume $HC\ne G.$} Then:  $\lambda_{H'}^{HC}(\Delta_0) = \lambda_{H'}^H(\Delta_0)\cdot \lambda_H^{HC}(\Delta_0)^{(H:H')}$  by induction. This we may substitute on the right side of {\bf (10)}, and in the resulting equality we may use {\bf (10)} in the particular case where $H'=H.$ Then we obtain our assertion {\bf (6)} in this case.\\
Thus we are left with {\bf the case $G=HC.$} Then $H\cap C$ is a normal subgroup of $G$, because it is normalized by $H$ and is a subgroup of the abelian group $C.$ Since $C$ has been chosen minimal, abelian, normal in $G$ the only possibilities are:\\
{\bf a)}  $H\cap C =C,$ that means  $H\supseteq C,$ hence $G=HC =H$ and our assertion follows from $\lambda_G^G(\Delta_0)=1,$  or\\
{\bf b)}   $H\cap C =\{1\}$ hence {\bf $G=H\ltimes C$ is semidirect. This is the case we are going to proceed with}. But before going further we need the following:\\

{\bf Lemma 3.5.} {\it Assume that $H$ is a complement of $C\le G$ that means $G=H\ltimes C$ is a semidirect product. For $H'\le H$ and a character $\chi:H'C\rightarrow \fC^\times$ consider the relation\\
{\bf (11)}
$$ Ind_{H'}^G(\chi_{H'})= Ind_{H'C}^G\left(Ind_{H'}^{H'C}(\chi_{H'})\right) \cong \sum_{\mu\in C^*/H'} Ind_{H_\mu'C}^G(\chi\mu'),$$
which is given by Lemma {\bf 2.2} (in the particular case $H'\cap C=\{1\}$) and where $H'C\ge H_\mu'C\ge C$ and $\mu':H_\mu'C\rightarrow\fC^\times$ is the extension of $\mu$ by $1.$ Then this will imply:\\
{\bf (12)}
$$ \Delta_0(H',\chi_{H'})\cdot \lambda_{H'}^G(\Delta_0) = \prod_{\mu\in C^*/H'} \left(\Delta_0(H_\mu'C,\chi\mu')\cdot\lambda_{H_\mu'C}^G(\Delta_0)\right).$$}\\
Note that the case $H'=\{1\},$ hence $\chi=\mu_0\in C^*,$ brings us back to the formula
$$ \lambda_N^\OO(\Delta)= \lambda_e^G(\Delta_0) =\prod_{\mu\in C^*}\Delta_0(C,\mu)\cdot \lambda_C^G(\Delta_0)^{\#C}=\lambda_e^C(\Delta_0)\cdot \lambda_C^G(\Delta_0)^{\#C},$$
which is an immediate consequence of definition {\bf 3.2} for $H=\{1\}.$\\

{\bf Proof:} Put $G':=H'C.$\\  {\bf If $H'\ne H,$ equivalently $G'\ne G,$} then by induction assumption $\Delta_0$ will extend from $R_1(\le G')$ onto $R(\le G')$. This we may apply to
$$  \chi\otimes Ind_{H'}^{G'}(1)=Ind_{H'}^{G'}(\chi_{H'}) \cong \sum_\mu Ind_{H_\mu'C}^{G'}(\chi\mu'),$$
hence\\
{\bf (13)}
$$  \Delta_0(H',\chi_{H'})\cdot\lambda_{H'}^{G'}(\Delta_0) =\prod_\mu \left( \Delta_0(H_\mu'C,\chi\mu')\cdot\lambda_{H_\mu'C}^{G'}(\Delta_0)\right).$$
Now we may apply {\bf (10)} with $H=G'=H'C,$ hence $HC=G'$, which gives us:\\
{\bf (14)}
$$  \lambda_{H'}^G(\Delta_0)= \lambda_{H'}^{G'}(\Delta_0)\cdot\lambda_{G'}^G(\Delta_0)^{(G':H')},\qquad
\lambda_{H_\mu'C}^G(\Delta_0)= \lambda_{H_\mu'C}^{G'}(\Delta_0)\cdot\lambda_{G'}^G(\Delta_0)^{(G':H_\mu'C)},$$
where the second equality follows from the first replacing $H'$ by $H_\mu'C.$
And the identity $(G':H')=\# C = \sum_{\mu\in C^*/H'} (G':H_\mu'C)$ yields
$$  \lambda_{G'}^G(\Delta_0)^{(G':H')} =\prod_{\mu \in C^*/H'} \lambda_{G'}^G(\Delta_0)^{(G':H_\mu'C)}.$$
Thus multiplying {\bf (13)} by this last equality and then using the equalities {\bf (14)} will give us {\bf (12)} in the case where $H'\ne H.$\\
{\bf Now assume $H'=H$ or equivalently}  $G'=G,$ hence $\chi$ is a character of $G=H\ltimes C,$ and in {\bf (11)} we have  $Ind_{H}^G(\chi_{H})=\chi\otimes Ind_{H}^G(1).$ Therefore {\bf (11)} is now the same as {\bf (4)} tensored by the character $\chi$ of $G.$ Since $H\cap C=\{1\}$ we have Lemma 2.2(i) with $S(\chi)= C^*,$ $T(\chi)=C^*/H.$ Moreover, using {\bf Lemma 3 of Appendix 2} we have:\\
{\it If $G=H\ltimes C$ is semidirect and $C$ is minimal normal (hence commutative), then $H< G$ is a maximal subgroup, and if $\wt C\ge C$ denotes the centralizer of $C,$ then:
$$\wt C\cap H =\bigcap_{g\in G} gHg^{-1} =:K$$
is normal in $G$ and $\ol G=G/K>\ol H= H/K$ is a type-III-group.}

{\bf Now we can finish the proof of Lemma 3.5}. (in the case $H'=H$) by using our assumption {\bf (3)}. Indeed $\Delta$ on $R_1(\le\OO)$ induces $\Delta_0$ on $R_1(\le \OO/N)$, and our type-III-group $\ol G =\ol H\ltimes C$ is actually a quotient of $G=\OO/N.$ Therefore the relation {\bf (3)} of $\Delta$ for type-III-groups which are subquotients of $\OO$ will induce a similar relation of $\Delta_0$ for type-III-groups which are quotients of $G.$ Thus referring to {\bf (3)} we may replace $\Delta,\;H<B,\;K,$ and $C,$ such that $HC=B,$ $H\cap C=K,$ and $\chi\in B^*$ by:\\
$\Delta_0,\; H<G,\; K,$ and $\wt C,$ such that $H\wt C =G,$ $H\cap\wt C=K,$ and $\chi\in G^*,$\\
which yields $G/K=(H/K)\ltimes C$ a type-III-group, where we may identify $\wt C/K =C,$ and turns our assumption {\bf (3)} into:\\
{\bf (15)}
$$  \Delta_0(H,\chi_H)\cdot \prod_{\mu\in C^*/H}\Delta_0(H_\mu C, \mu')= \prod_{\mu\in C^*/H}\Delta_0(H_\mu C, \chi\mu').$$
Here we have used that the inner action of $K=H\cap\wt C$ on $C^*$ is trivial, hence $K\le H_\mu$ and $H_\mu\wt C=H_\mu C$ for all $\mu.$ Finally, by our definition of $\lambda_H^G(\Delta_0)$ the equality {\bf (15)} (which holds by assumption!) turns into {\bf (12)} for $H'=H,$ if we multiply it by $\prod_{\mu\in C^*/H} \lambda_{H_\mu C}^G(\Delta_0).$
{\bf qed.}\\

{\bf We come back to the proof of Lemma 3.4.} where we were left to show that: $\lambda_{H'}^G(\Delta_0) = \lambda_{H'}^H(\Delta_0)\cdot\lambda_H^G(\Delta_0)^{(H:H')},$ if $G=H\ltimes C$ is semidirect, $C$ is minimal normal in $G$, and $H'\le H.$
Because of the {\bf Remark 1} from the proof of {\bf Proposition 1.4} it is sufficient to show that
$$   H'\in \sU(H) \mapsto  \lambda_{H'}^H:= \lambda_{H'}^G(\Delta_0)\lambda_H^G(\Delta_0)^{-(H:H')}$$
satisfies the conditions {\bf (c1)} and {\bf (c2)}. In other words we show that the function $\Delta_0$ on $R_1(\le H)$ is extendible to $\sF_0$ on $R(\le H)$ if we make use of $\lambda_{H'}^H$ as we have defined it. Then of course
$$  \lambda_{H'}^H = \sF_0(Ind_{H'}^H(1)) =\lambda_{H'}^H(\Delta_0)$$
because the extension $\sF_0$ is unique.\\
Thus we turn to the condition {\bf (c2)} of Proposition 1.4, where we consider $\Delta_0$ on $R_1(\le H)$.  Assume a relation\\
{\bf (16)}
$$  \sum_{i=1}^r n_i Ind_{H_i'}^H(\chi_i) \cong 0,$$
for certain subgroups $H_1',...,H_r'\le H,$ multiplicities $n_i\in\Z$ and characters $\chi_i\in {H_i'}^*.$
.\\
We have to show that {\bf (16)} implies\\
{\bf (17)}
$$  \prod_i\left(\Delta_0(H_i',\chi_i)\lambda_{H_i'}^H\right)^{n_i} =1,$$
for $\lambda_{H_i'}^H$ as defined above.\\
First of all from {\bf (16)} we see: $\sum_i n_i(H:H_i') =0,$ hence $\prod_i \lambda_H^G(\Delta_0)^{-(H:H_i')n_i} =1.$ Thus in order to verify {\bf (17)} it is enough to see that {\bf (16)} will imply\\
{\bf (18)}
$$ \prod_i\left(\Delta_0(H_i',\chi_i)\lambda_{H_i'}^G(\Delta_0)\right)^{n_i} =1.$$

For each pair $(H_i',\chi_i)$ we may consider the subgroup $G_i':= H_i'\ltimes C \le G= H\ltimes C$ and interpret $\chi_i$ as character of $G_i'/C$ and to this situation apply {\bf Lemma 3.5}. Then we obtain relations of type {\bf (11)} and {\bf (12)} resp. which we will write as:\\
{\bf (19)}
$$  Ind_{H_i'}^G(\chi_i)\cong \sum_{j=1}^{r_i} Ind_{H_{ij}'}^G(\chi_i\mu_{ij}'),$$
where $\{\mu_{ij}\}_j \subset C^*$ are representatives for the cosets $C^*/H_i'$ and where $H_{ij}':=Stab_{H_i'}(\mu_{ij})C,$\\
and $\chi_i$ is understood to be trivial on $C,$ and $\mu_{ij}'$ extends $\mu_{ij}$ hence is non-trivial on $C.$
And {\bf (12)} applied to $H'=H_i'$ reads now as follows:\\
{\bf (20)}
$$  \Delta_0(H_i',\chi_i)\lambda_{H_i'}^G(\Delta_0)= \prod_{j=1}^{r_i} \left( \Delta_0(H_{ij}',\chi_i \mu_{ij}') \lambda_{H_{ij}'}^G(\Delta_0)\right).$$
Implementing this into {\bf (18)} the assertion turns into\\
{\bf (21)}
$$ \prod_{i=1}^r\left(  \prod_{j=1}^{r_i} \left( \Delta_0(H_{ij}', \chi_i\mu_{ij}') \lambda_{H_{ij}'}^G(\Delta_0)\right)\right)^{n_i} =1.$$
On the other hand, applying $Ind_H^G$ to {\bf (16)} and using {\bf (19)} we obtain:\\
{\bf (22)}
$$ \sum_{i=1}^r n_i\left(\sum_{j=1}^{r_i} Ind_{H_{ij}'}^G(\chi_i\mu_{ij}')\right) \cong 0,$$
where all $H_{ij}'\ge C$ and therefore $\lambda_{H_{ij}'}^G(\Delta_0) = \lambda_{H_{ij}'/C}^{G/C}(\Delta_0).$\\
{\it The decisive point is, that $\Delta_0$ comes from $\Delta$ on $R_1(N\le\OO).$} Via $\OO\twoheadrightarrow \OO/N=G$ let $\wh C\le \wh{H_{ij}'}$ be the full preimages of $C\le H_{ij}'.$ Then $\wh C >N> [\wh C,\wh C]$ is normal in $\OO$ and via inflation we may interpret $\chi_i\mu_{ij}'$ as $(\wh{H_{ij}'}, \chi_i\mu_{ij}') \in R_1(\wh C \le\OO).$\\
{\bf Now we may use our assumption that $\Delta$ is extendible from $R_1(\wh C\le\OO)$ onto $R(\wh C\le\OO)$ because $(\OO:\wh C)< (\OO:N).$} Therefore, interpreting {\bf (22)} as a representation of $\OO/[\wh C,\wh C]$ we obtain:
$$ \prod_{i=1}^r\left(  \prod_{j=1}^{r_i} \left( \Delta(\wh{H_{ij}'},\chi_i \mu_{ij}') \lambda_{\wh{H_{ij}'}}^\OO(\Delta)\right)\right)^{n_i} =1.$$
Here we consider $\chi_i\mu_{ij}'$ as a character of $\wh{H_{ij}'} \ge \wh C$ which is trivial on $N>[\wh C,\wh C].$   Therefore we may replace $\Delta$ by $\Delta_0$, where $\lambda_{\wh{H_{ij}'}}^\OO(\Delta)= \lambda_{\wh{H_{ij}'}/N}^{\OO/N}(\Delta_0)=\lambda_{H_{ij}'}^G(\Delta_0)$ and $\lambda_{H_{ij}'}^G(\Delta_0) =\lambda_{H_{ij}'/C}^{G/C}(\Delta_0) $ by definition,
which turns our last equality into the asserted equality {\bf (21)} and ends the proof of the $N$-Tower lemma 3.4. {\bf qed.}\\
{\bf Remark:} The relation {\bf (22)} is not concerned with representations of $G/C.$ Therefore the induction $(G:C)<\#G$ with respect to $\Delta_0$ does not work; instead we have to use the condition $\wh C >N> [\wh C,\wh C]$. On one hand $G=\OO/N$ is a quotient of $\OO/[\wh C,\wh C],$ and on the other hand  $\Delta$ is extendible from $R_1(\wh C\le\OO)$ onto $R(\wh C\le\OO)$ because $(\OO:\wh C)< (\OO:N).$\\

{\bf Step 3}: (Using the generators of $Ker(\vf_{N\le\OO})$ and the $N$-tower Lemma 3.4.)  So far we have seen that the function $\Delta$ on $R_1(N\le\OO)$ gives rise to $\Delta_0$ on $R_1(\le\OO/N)$ and a function
$$ H\in \sU(G)\mapsto \lambda_H^G(\Delta_0),\qquad\textrm{where}\; G:=\OO/N,$$
which satisfies the tower lemma. Finally in step 3 we want to verify that the definition $\lambda_{U'}^U(\Delta):= \lambda_{H'}^H(\Delta_0)$ for $N\le U'\le U\le \OO$ and $H'=U'/N\le H=U/N,$  implies the extension of $\Delta$ from $R_1(N\le\OO)$ onto $\sF$ on $R(N\le \OO).$ (So far we have only used the induction assumption that $\Delta$ extends from $R_1(\wh C\le\OO)$ onto $R(\wh C\le\OO).$ Now the assumptions {\bf (1) - (3)} will come in.)\\
For the extension onto $R(N\le\OO)$ we will use the {\bf criterion 1.4} that means we have to verify {\bf (c2)} for $H=\OO$, which comes down {\bf to check the generators I/N, II/N, III/N of the kernel of the Brauer map }
$$  \vf_{N\le\OO}: R_+(N\le \OO)\twoheadrightarrow  R(\OO/_{[N,N]}).$$

{\bf Type I/N:} The basic relation $\rho(K,B,\chi)$ refers to
$$ Ind_K^B(\chi_K) =\sum_{\mu\in (B/K)^*} \chi\mu,\qquad Ind_K^\OO(\chi_K) =\sum_{\mu\in (B/K)^*} Ind_B^\OO(\chi\mu),$$
for $\chi\in B^*,$ where $B/K$ is cyclic of prime order $\ell.$ Here {\bf our assumption  (1)} together with: $\prod_{\mu\in (B/K)^*} \Delta(B,\mu) =\lambda_K^B(\Delta),$ (because $B/K$ is cyclic of prime order $\ell$ and a subquotient of $G=\OO/N$, thus by induction we have {\bf 1.(9)} at our disposal), rewrites as
$$ \Delta(K,\chi_K) \lambda_K^B(\Delta) = \prod_{\mu\in (B/K)^*} \Delta(B,\chi\mu),$$
and multiplying this by $\lambda_B^\OO(\Delta)^\ell$ and using the $N$-tower relation {\bf 3.4} we obtain
$$ \Delta(K,\chi_K) \lambda_K^\OO(\Delta) = \prod_{\mu\in (B/K)^*} \left(\Delta(B,\chi\mu)\lambda_B^\OO(\Delta)\right),$$
which verifies the condition {\bf 1.4 (c2)} for the generating relations of type I/N.\\

{\bf Type II/N:} Here we have to verify that $Ind_{H_1}^\OO(\eta^{H_1}) \cong Ind_{H_2}^\OO(\eta^{H_2})$ implies\\
{\bf (23)}
$$  \Delta(H_1,\eta^{H_1})\lambda_{H_1}^\OO(\Delta) =  \Delta(H_2,\eta^{H_2})\lambda_{H_2}^\OO(\Delta),$$
where $\OO \ge B > H_1,\;H_2 > Z=H_1\cap H_2 \ge N,$ and $\eta:Z/[Z,B]\rightarrow \fC^\times$ is a fixed character which is nontrivial on $[B,B]/[Z,B],$ and $\eta^{H_i}\in H_i^*$ is any extension of $\eta.$\\
From {\bf 1.(9)} we know that $\lambda_{H_i}^B(\Delta) = \prod_{\mu\in(B/H_i)^*} \Delta(B,\mu)$ because $B/H_i$ is cyclic of prime order $\ell$ and is a subquotient of $G=\OO/N.$ Therefore our assumption {\bf (2)} rewrites as:
$$  \Delta(H_1,\eta^{H_1})\lambda_{H_1}^B(\Delta) =  \Delta(H_2,\eta^{H_2})\lambda_{H_2}^B(\Delta).$$
Now multiplying this by $\lambda_B^\OO(\Delta)^\ell$ and using the $N$-tower lemma {\bf 3.4} we obtain {\bf (23)}.\\

{\bf Type III/N:} We consider the basic relation
$$\rho(H,B,\chi)= Ind_B^\OO(\rho_B(H,\chi)) =[H,\chi_H]-\sum_{\mu\in(C/K)^*/H}[H_\mu C,\chi_{H_\mu C}\mu']\quad \in \sR(N\le\OO),$$
for open subgroups $H< B\le \OO$ where $H< B$ is maximal, and $K< C< B$ such that $HC=B$ and $H\cap C=K:=\bigcap_{b\in B} bHb^{-1},$ and moreover we assume here that $B/K$ is a subquotient of $G=\OO/N.$ Corresponding to $\vf_{N\le\OO}(\rho(H,B,\chi))=0$ {\bf we have $Ind_H^\OO(\chi_H)\cong\sum_\mu Ind_{H_\mu C}^\OO(\chi_{H_\mu C}\mu'),$ and in this situation we have to verify 1.4(c2)}.\\
{\bf By assumption} we have the relation\\
{\bf (3)}
$$ \Delta(H,\chi_H)\cdot \prod_{\mu\in(C/K)^*/H}\Delta(H_\mu C,\mu') =\prod_{\mu\in(C/K)^*/H}\Delta(H_\mu C,\chi\mu'),$$
and due to $Ind_H^B(1)=\sum_{\mu\in(C/K)^*/H} Ind_{H_\mu C}^B(\mu'),$ and restricting $\Delta$ to $R_1(K\le B),$ we have
$$ \lambda_H^B(\Delta)= \lambda_{H/K}^{B/K}(\Delta_0)= \prod_{\mu\in(C/K)^*/H}\Delta(H_\mu C,\mu')\cdot\lambda_{H_\mu C/C}^{B/C}(\Delta_0)$$
by induction.
Restricting further to $R_1(C\le B)$ we may write: $\lambda_{H_\mu C/C}^{B/C}(\Delta_0) = \lambda_{H_\mu C}^B(\Delta).$ The only case where we cannot apply the induction assumption is the case $B=\OO,$ $K=N,$ but then the definition {\bf 3.2} applies because $C/K$ is the uniquely determined abelian normal subgroup in $B/K.$\\
Therefore multiplying {\bf (3)} by $\prod_\mu \lambda_{H_\mu C}^B(\Delta)$ we conclude:
$$  \Delta(H,\chi_H)\lambda_H^B(\Delta) = \prod_\mu \Delta(H_\mu C,\chi\mu')\lambda_{H_\mu C}^B(\Delta).$$
Finally multiplying this last equality by $\lambda_B^\OO(\Delta)^{(B:H)} = \lambda_B^\OO(\Delta)^{\sum_\mu (B:H_\mu C)}$ and using the tower relation {\bf 3.4} for $N\le\OO,$ we arrive at
$$  \Delta(H,\chi_H)\lambda_H^\OO(\Delta) = \prod_\mu \Delta(H_\mu C,\chi\mu')\lambda_{H_\mu C}^\OO(\Delta),$$
which confirms the criterion {\bf 1.4(c2)} for the generating relation $\rho(H,B,\chi)\in Ker(\vf_{N\le\OO})$ of type III.\\
Altogether this finishes our induction argument that $\Delta$ extends from $R_1(N\le\OO)$ onto $R(N\le\OO)$ if we assume this for proper subquotients $N'\le\OO'$. {\bf qed. Theorem 3.1}.\\

{\bf 4. Recovering Theorem 3.1 for the case of local root numbers.}\\

So far, beginning from a profinite solvable group $\OO$ we have dealt with the question as to when a function $\Delta =\Delta(H,\chi)$ on $R_1(\le\OO)$ can be extended to a function
$$  \sF=\sF(H,\rho)\quad\textrm{on}\quad R(\le\OO)=\bigsqcup_{H\le\OO} R(H)$$
where $H$ runs over the open subgroups of $\OO$ and $\chi\in H^*$ is replaced by any virtual representation $\rho\in R(H).$ {\it Now we come to the arithmetic application}:\\
As our base field we fix a non-archimedean local field $F$ with finite residue field of characteristic $p$ and consider $\OO = G_F$ the absolute Galois group over $F$ or $\OO=W_F\subset G_F$ the absolute Weil group. Thus fixing a separable closure $\ol F|F$ we obtain
\begin{itemize}
 \item $R_1(\le G_F)$ the set of all pairs $(E,\chi)$ where $E|F$ is a finite subextension that means $G_E\le G_F$ is an open subgroup and $\chi\in G_E^*$ is via class field theory a continuous character of $E^\times$ which is of finite order.
 \item $R_1(\le W_F)$ is the set of all pairs $(E,\chi)$ which now means that $W_E\le W_F$ is an open subgroup of finite index and $\chi\in W_E^*$ is via class field theory a quasicharacter of $E^\times\cong W_E/[W_E,W_E].$
\end{itemize}
Let $\psi_F$ be a {\bf fixed non-trivial additive character} of $F$, and take
$$  \psi_{E|F}:= \psi_F\circ Tr_{E|F}$$
as additive characters for the extension fields $E|F.$ Then in both cases $\OO=G_F$ or $\OO=W_F$ our function on  $R_1(\le\OO)$ to begin with is:\\
{\bf (0)}
$$ \Delta(E,\chi):= \Delta(\chi,\psi_{E|F}):= \chi(c)\cdot q_E^{-a_E(\chi)/2}\sum_{x\in U_E/U_E^{a_E(\chi)}} \chi^{-1}(x) (c^{-1}\psi_{E|F})(x),$$
which is Tate's local root number assigned to the Gauss sum for the pair $(\chi,\psi_{E|F}).$ More precisely $q_E=\#\kappa_E$ is the order of the residue field,
$a_E(\chi) = min\{a\;|\;U_E^a = 1+\pp_E^a\subseteq Ker(\chi)\}$ is the exponential conductor of $\chi$, and $(c^{-1}\psi_{E|F})(x):= \psi_{E|F}(xc^{-1})$ is a translation of $\psi_{E|F}$ such that $\chi:E^\times/1+\pp_E^a \rightarrow \fC^\times$ and $c^{-1}\psi_{E|F}:E/\pp_E^a \rightarrow\fC^\times$ are given on parallel levels $a=a_E(\chi).$ Then the Gaussian sum adds the character $\chi^{-1}$ of $U_E/U_E^a=(\oo_E/\pp_E^a)^\times$ against the  character $c^{-1}\psi_{E|F}$ of $(\oo_E/\pp_E^a)^+.$\\ 
{\bf Defining the level} $\ell(\psi_{E|F}):=min\{a\;|\;\pp_E^a\subseteq Ker(\psi_{E|F})\}$, ($=-n(\psi_{E|F})$ in terms of {\bf [L]}) we have
$$ \ell(c^{-1}\psi_{E|F})= \ell(\psi_{E|F}) + \nu_E(c),$$
hence for computing $\Delta(E,\chi)$ we need $\nu_E(c)= a_E(\chi)-\ell(\psi_{E|F})$,\\ and the value $\Delta(E,\chi)$ will not change if we modify $c$ by a unit factor.\\ In terms of {\bf [L], p.5} the definition {\bf (0)} rewrites as 
$$  \Delta(E,\chi) = \chi(c)\cdot\frac{\int_{U_E} \chi^{-1}(x) (c^{-1}\psi_{E|F})(x)\; dx}{|\int_{U_E}\chi^{-1}(x) (c^{-1}\psi_{E|F})(x)\;dx|}= \chi(c)\cdot\frac{\sum_{x\in U_E/U_E^{a_E(\chi)}} \chi^{-1}(x) (c^{-1}\psi_{E|F})(x)}{|\sum_{x\in U_E/U_E^{a_E(\chi)}} \chi^{-1}(x) (c^{-1}\psi_{E|F})(x)|}$$
where the second equality follows if we consider $U_E=\oo_E^\times =\bigsqcup_{x\in(\oo_E/\pp_E^a)^\times} x+\pp_E^a$ as a disjoint union of cosets and then use that
$\chi^{-1}(x)(c^{-1}\psi_{E|F})(x)$ is constant on the cosets for $a=a_E(\chi)$.    Moreover keep $a=a_E(\chi)$ and put $R=\oo_E/\pp_E^a,$
$P=\pp_E/\pp_E^a$, $\psi:= c^{-1}\psi_{E|F}$ hence $\ell(\psi)=a.$ Then
$$ |\sum_{x\in R^\times}\chi^{-1}(x) \psi(x)|^2=\sum_{x,x_1\in R^\times}\chi(x_1x^{-1})\psi(x-x_1)=\sum_{y\in R^\times}\chi(y)\left(\sum_{x\in R}\psi(x(1-y))-\sum_{x\in P}\psi(x(1-y))\right)$$
$$ = \left(\sum_{y\in 1+P^{a-1},\,y\ne 1} \chi(y)\right)\cdot(-\# P) + 1\cdot(\#R -\#P)= \#R = q_E^{a_E(\chi)},$$
where we have used $P^a=\{0\},$ $1+P^a=\{1\}$ and $\ell((1-y)\psi)>1$ if $y\notin 1+P^{a-1}.$ Therefore {\bf (0)} and {\bf [L]} agree.\\
On the other hand {\bf [T3],(3.2.6)} implies that:
$$  \Delta(E,\chi) =\ve(\chi|.|_E^{1/2},\psi_{E|F},d_Ex)\quad \textrm{for $d_E=$ the $\psi_{E|F}$-selfdual Haar-measure on $E^+$}.$$
From the explicit definition {\bf (0)} it is obvious that:
\begin{itemize}
 \item $\Delta(E,1)=1$ and $\Delta(E,\chi)=\Delta(E^g,\chi^g)$ for all $g\in G_F$, where $E^g=g^{-1}(E)$ and $\chi^g=\chi\circ g,$ hence the basic properties {\bf 1.(1), 1.(2)} to be expected of a function $\Delta$ are fulfilled.
\end{itemize}
Since the absolute Galois group $G_F$ is profinite solvable, our {\bf Theorem 3.1} directly applies to the question if $\Delta(E,\chi)$ on $R_1(\le G_F)$ will extend to $R(\le G_F)= \bigsqcup_{E|F} R(G_E).$\\

On the other hand $W_F$ is not profinite but comes as the extension of $\Z$ with the profinite group $W_F^0=G_F^0=$ inertia subgroup. Thus we have a continuous embedding $W_F\hookrightarrow G_F$ with dense image and identical commutator subgroups $[W_F,W_F]=[G_F,G_F].$ Since $W_F/W_F^0\cong \Z,$ the group $W_F$ admits representations which are not semisimple. But {\bf our aim is to extend the local root numbers $\Delta(E,\chi)$ only to semisimple representations of the Weil groups.}
Therefore from now we think of
\begin{itemize}
 \item $ R(W_F):=$ the free $\Z$-module over the equivalence classes of irreducible smooth representations of $W_F.$
\end{itemize}
According to {\bf [D]\S4.10} every irreducible smooth representation of $W_F$ comes up to unramified twist as the restriction of an irreducible representation of $G_F.$ Therefore the tensor product of irreducibles will be semisimple such that $R(W_F)$ becomes a commutative ring, and the restriction of representations induces an injection of rings $R(G_F)\hookrightarrow R(W_F)$. Also the induction $Ind_{E|F}:= Ind_{W_E}^{W_F}$ takes $R(W_E)$ into $R(W_F)$ if $W_E< W_F$ is an open subgroup of finite index. (For these and further details one may also consult {\bf 28.} of the textbook {\bf [BH]}).
Because of $|N_{E|F}|^s_F = |.|^s_E$ the inductions $Ind_{E|F}$ will commute with unramified twisting such that
the modifications Brauer 3 and Brauer 4 apply to $\OO=W_F$ as well and refer now to the relative Weil groups $W_{K|F} = W_F/[W_K,W_K]$ for finite normal extensions $W_K=N < \OO=W_F.$\\
{\bf Brauer 3:} {\it The Brauer map $\vf_{W_{K|F}}:R_+(\le W_{K|F})\rightarrow R(W_{K|F})$ induces a surjection
$$  \vf_{K|F}:R_+(W_{K|K}\le W_{K|F}) \twoheadrightarrow R(W_{K|F}),$$
where $R_+(W_{K|K}\le W_{K|F}) =: R_+(K\ge F)$ is the free abelian group over all $F$-conjugacy classes $[E,\chi_E]$ consisting of a subfield $E|F$ in $K$ and a quasicharacter $\chi_E:E^\times\rightarrow\fC^\times$, and
$$  \vf_{K|F}([E,\chi_E]) = Ind _{E|F}(\chi_E).$$}
{\bf Brauer 4:} {\it If $\rho\in R(W_{K|F})$ is of dimension $0$, then it admits a presentation
$$  \rho = \sum n_i Ind_{E_i|F}(\chi_{E_i} -  1_{E_i})$$
for subfields $E_i\subseteq K.$}\\
The projector $\Phi_C$ of Corollary {\bf 2.3} turns now into a projector
$$\Phi_{W_{K|K}}:R_+(\le W_{K|F})\twoheadrightarrow R_+(W_{K|K}\le W_{K|F})\quad\textrm{such that}\quad \vf_{W_{K|F}} = \vf_{K|F}\circ \Phi_{W_{K|K}},$$
and the relative {\bf criterion 1.4} deals with extending a function $\Delta$ from $R_1(W_K\le W_F) = \bigsqcup_{K\ge E\ge F} (E^*)$ (for quasicharacters) onto $R(W_K\le W_F) = \bigsqcup_{K\ge E\ge F} R(W_{K|E}).$\\
Our modest aim in this last section is to translate the sufficient conditions {\bf 3.(1) - 3.(3)} into conditions for the local root numbers $\Delta(E,\chi)=\Delta(\chi,\psi_{E|F}).$ Thus we have to translate pairs $(H,\chi)$ for $H\le\OO$ and $\chi\in H^*$ into pairs $(E,\chi)$ for finite extensions $E|F$ and continuous quasicharacters $\chi:E^\times\rightarrow \fC^\times.$\\
{\bf First of all}, in all of {\bf 3.(1) - (3)} we have ingredients $B\le\OO$ which means we go from our base field $F$ to some finite extension $F'|F$ such that $B=W_{F'}\le\OO=W_F.$ Then up to this change of base field the conditions look all the same. {\bf Expecting now that the verification of 3.(1) - (3) will not depend on taking $F'$ or $F$ as our base field,} we will restrict here to expressing the conditions relative to the base field $F.$ In other words we restrict to cases $B=\OO,$ $F'=F,$
and $\chi\in B^*=\OO^*$ is via class field theory nothing else than a quasicharacter of $F^\times.$ Then:\\
{\bf in order to check the extendibility of $\Delta(E,\chi)$ to local root numbers $\ve_\Delta(E,\rho)$ for all $\rho\in R(W_E)$ we have to verify the following arithmetic conditions:}
\begin{itemize}
 \item[{\bf (1)}] {\bf Davenport Hasse = [L], First Main Lemma:} Let $K|F$ be a cyclic extension of prime degree $\ell$ and let $S(K|F)$ be the characters of $F^\times/N_{K|F}(K^\times) \cong G_{K|F}$. Then:
 $$ \Delta(K,\chi\circ N_{K|F})\cdot\prod_{\mu\in S(K|F)}\Delta(F,\mu) =\prod_{\mu\in S(K|F)} \Delta(F,\chi\mu)$$
 for all quasicharacters $\chi$ of $F^\times.$
 \item[{\bf (2)}] {\bf Heisenberg identity = [L], Second Main Lemma:} Let $K|F$ be an abelian extension of type $G_{K|F}\cong \Z/\ell \times \Z/\ell$ for some prime number $\ell,$ and let $\eta: K^\times/I_F K^\times\rightarrow \fC^\times$ be a $G_{K|F}$-invariant quasicharacter which is non-trivial on the norm-1-subgroup $\wh H^{-1}(G_{K|F},K^\times) =K^\times_F/I_F K^\times \cong \Z/\ell.$ Consider subextensions $E|F$ of degree $\ell$ inside $K|F$ and quasicharacters $\eta^E$ of $E^\times$ such that $\eta^E\circ N_{K|E} =\eta.$ Then:\\
 {\it the expression $\Delta(E,\eta^E)\cdot\prod_{\mu\in S(E|F)} \Delta(F,\mu)$ does not depend on the choice of $(E,\eta^E).$}
 \item[{\bf (3)}] {\bf Generalized Davenport Hasse = [L], Third and Fourth Main Lemma:} Let $E|F$ be an extension which is minimal (in the sense that there are no proper intermediate fields), and let $K|F$ be the Galois closure of $E|F.$ {\bf We assume here $K\ne E$} because otherwise we are in the situation of {\bf (1)}. Then $G_{K|F}$ is a type-III-group and we have a unique {\bf normal} subextension $L|F$ such that $C=G_{K|L}$ is {\bf abelian}, and moreover: $L\cdot E =K,$  $L\cap E =F.$
 Now the arithmetic condition is:
 $$ \Delta(E,\chi\circ N_{E|F})\cdot\prod_{\mu\in S(K|L)/H}\Delta(F_\mu,\mu') = \prod_{\mu\in S(K|L)/H}\Delta\left(F_\mu,\;(\chi\circ N_{F_\mu|F})\cdot\mu'\right)$$
for all quasicharacters $\chi$ of $F^\times$ and characters $\mu$ of $L^\times/N_{K|L}(K^\times) \cong G_{K|L}$ where we restrict to a system of representatives with respect to the action of the maximal subgroup $H =G_{K|E}\cong G_{L|F}$ on $L^\times.$ And $F_\mu|F$ denotes the subfield of $L$ such that $H_\mu=G_{K|E_\mu} \cong  G_{L|F_\mu}$ is the stabilizer of $\mu$ under that action. Finally with respect to $G_{K|F_\mu}= G_{K|E_\mu}\cdot G_{K|L}$ the character $\mu'$ of $(F_\mu)^\times/N_{K|F_\mu}(K^\times)$ is defined via: $\mu'\circ N_{E_\mu|F_\mu}\equiv 1$ and
 $\mu'\circ N_{L|F_\mu}=\mu.$

 {\bf Remark:} From the displayed formula we may separate the trivial character $\mu_0\equiv 1$ such that $F_{\mu_0}=F$ and $\mu_0'$ is the trivial character of $F^\times.$ Thus for $\mu=\mu_0$ we obtain the factors $\Delta(F,1)=1$ on the left side, and $\Delta(F,\chi)$ on the right.
\end{itemize}
A detailed translation of {\bf 3.(1) - 3.(3)} (for the special case $B=\OO$) into the arithmetic conditions {\bf (1) - (3)} uses only standard facts from local class field theory. The underlying relations in $R(W_F)$ are the following:
\begin{itemize}
 \item[{\bf (R1)}] $Ind_{K|F}(\chi\circ N_{K|F}) = \sum_{\mu\in S(K|F)}\chi\mu,$
 \item[{\bf (R2)}] $Ind_{E|F}(\eta^E) =\frac{1}{\ell} Ind_{K|F}(\eta)$ for all pairs $(E,\eta^E),$
 \item[{\bf (R3)}] $Ind_{E|F}(\chi\circ N_{E|F})=\sum_{\mu\in S(K|L)/H} Ind_{F_\mu|F}((\chi\circ N_{F_\mu|F})\mu').$
\end{itemize}

\centerline{{\bf  Unramified twisting}}

Here we want to see that {\bf the arithmetic conditions {\bf (1)-(3)} are stable under unramified twisting,} that means if (1) or (3) are valid for a certain $\chi\in {F^\times}^*,$ then also for the unramified twists $\chi|.|_F^s$ and similarly if (2) holds for a certain $\eta\in(K^\times/I_FK^\times)^*$ then also for the unramified twists $\eta|.|_K^s.$ Thus from each  equivalence class of quasicharacters modulo unramified twisting we need to check only one representative and for this we may choose a character of finite order. The argument is similar as in {\bf [T2], end of \S2}. Because of the explicit formula {\bf (0)} and using that $\chi$ and $\chi|.|_E^s$ have the same restriction to $U_E$ we obtain:\\
{\bf (4)}
$$  \frac{\Delta(E,\chi)}{\Delta(E,\chi|.|_E^s)}=\frac{1}{|c|_E^s}=q_E^{s\cdot(a_E(\chi)-\ell(\psi_{E|F}))},$$
for all quasicharacters $\chi$ of $E^\times.$
Now the argument is based on the following\\

{\bf Lemma:} {\it For finite separable extensions $E|F$ and semisimple virtual representations $\rho_E\in R(W_E),$ the homomorphic function
$$  (E,\rho_E)\mapsto N(E,\rho_E):=q_E^{a_E(\rho_E)- dim(\rho_E)\ell(\psi_{E|F})}\quad \in \R_+^\times$$
is inductive in the unconditional sense:
$$  N(E,Ind_{K|E}(\rho_K)) = N(K,\rho_K).$$}
As a generalization of $a_E(\chi)$ from above we have used here the exponential Artin conductor $a_E(\rho_E).$ (In terms of {\bf [S],chap.VI,\S2} it is the intertwining number $f(\rho_E)$ with the Artin representation).\\

{\bf Proof:} Obviously it is enough to deal with the case $E=F.$ Then using the well known formulas (cf. {\bf [S]} p.101 below, and chap. III Prop.7 resp.):
$$ a_F(Ind_{K|F}(\rho_K)) = f_{K|F}(d_{K|F}\cdot dim(\rho_K)+a_K(\rho_K)),\qquad e_{K|F}\ell(\psi_F) = d_{K|F} +\ell(\psi_{K|F}),$$
where $d_{K|F}$ denotes the differental exponent, we obtain:
$$  a_F(Ind_{K|F}(\rho_K)) - [K:F]dim(\rho_K)\ell(\psi_F)= $$
$$ = f_{K|F}\left(d_{K|F}\cdot dim(\rho_K) + a_K(\rho_K)\right) - f_{K|F}dim(\rho_K)(d_{K|F} +\ell(\psi_{K|F}))=$$
$$= f_{K|F}\left(a_K(\rho_K) - dim(\rho_K)\ell(\psi_{K|F})\right).$$
{\bf qed.}

Now for {\bf the invariance of (1) and (3) resp. under unramified twisting} we check the quotient relations which via {\bf (4)} rewrite as
$$  N(K,\chi\circ N_{K|F})^s =\prod_{\mu\in S(K|F)} N(F,\chi \mu)^s,\qquad
N(E,\chi\circ N_{E|F})^s =\prod_{\mu\in S(K|F)/H} N(F_\mu,(\chi\circ N_{F_\mu|F})\mu')^s$$
resp. But due to the Lemma these equalities are direct consequences of the underlying relations (R1), (R3) in $R(W_F).$\\
Similarly concerning {\bf (2)} we have to check that
$N(E,\eta^E)^s$
does not depend on the choice of $(E,\eta^E).$ Again this is clear from the Lemma because  $Ind_{E|F}(\eta^E)$ does not depend on that choice.\\
Once the extendibility of $\Delta$ has been verified the Lemma will imply that {\bf (4)} extends to
$$  \frac{\ve_\Delta(E,\rho)}{\ve_\Delta(E,\rho\otimes|.|_E^s)} = N(E,\rho)^s$$
for all $\rho\in R(W_E).$\\

\centerline{{\bf A simplification in the case of odd order}}
Our $\Delta =\Delta(E,\mu)$ as in {\bf (0)} satisfies the functional equation
$$  \Delta(E,\mu)\cdot \Delta(E,\mu^{-1}) = \mu(-1).$$
If now $\mu$ is of odd order then it can be recovered as a power of $\mu^2$, hence
$\Delta(E,\mu^{-1}) = \Delta(E,\mu)^{-1}$ in that case. {\bf Therefore if $\ell\ne 2$, or if (in case (3)) the Galois closure of $E|F$ is of degree $[K:F]$ odd,} then we are left with the arithmetic conditions:\\
{\bf (1odd)}
$$\Delta(K,\chi\circ N_{K|F}) =\prod_{\mu\in S(K|F)} \Delta(F,\chi\mu),$$
{\bf (2odd)}\\
{\it the expression $\Delta(E,\eta^E)$ does not depend on the choice of $(E,\eta^E),$}\\
{\bf (3odd)}
$$ \Delta(E,\chi\circ N_{E|F})  = \prod_{\mu\in S(K|L)/H}\Delta\left(F_\mu,\;(\chi\circ N_{F_\mu|F})\cdot\mu'\right),$$
if we take the system of representatives $\mu$ for $S(K|L)/H$ such that it is stable under $\mu\mapsto \mu^{-1}.$ {\bf Actually we are free to choose such a system as we like, because our function $\Delta$ has the invariance-property 1.(2).}
Note here that all $\mu$ are of odd order and the acting group $H$ is of odd order too. Therefore we will have $\mu\ne \mu^{-1}$ and moreover $\mu,\;\mu^{-1}$ cannot be $H$-equivalent, if $\mu\not\equiv 1.$\\

Langlands {\bf [L]} divides the condition {\bf (3)} into two subcases: {\bf the Third and the Fourth Main Lemma of tame and of wild ramification resp.} which are dealt with in chapters 13 and 14 of loc.cit. As an introduction we consider\\
{\bf Type-III-groups as local Galois groups}\\
For the group theoretical background we rely on Appendix 2. And we keep the notations $K\supseteq L, E \supseteq F$ from {\bf (3)}. Thus let
$$  H= G_{K|E} <  G= G_{K|F},\qquad H\cong G/C =G_{L|F}$$
be a type-III-group realized as local Galois group, where we exclude now the degenerate case $H=\{1\},$ which means $G=C$ is cyclic of prime order $\ell$. Equivalently $E|F$ is a minimal extension which is not Galois and $K|F$ is its Galois closure. Due to Lemma 2 of {\bf A2} this actually means that $E|F$ is anti-Galois in the sense that the roots of a minimal polynomial generate pairwise different extensions $E_i|F$ precisely one of them: $E|F$. Then $E|F$ cannot be unramified and it cannot have mixed ramification because this would imply existence of an intermediate field. Thus the only cases are:
\begin{itemize}
 \item $E|F$ is totally tamely ramified of prime degree $[E:F]=\ell \ne p$. Then $L|F$ is the inertia subfield in $K|F$, the extension $K|L$ is cyclic and totally ramified of degree $\ell$ and $[K:E]=[L:F]$ is the order of $q_F\in (\Z/\ell)^\times.$ The degenerate case happens for $q_F\equiv 1\;(mod\;\ell).$ Considering $C=G_{K|L}\twoheadleftarrow U_L/U_L^1$ as a 1-dimensional $\F_\ell$-space we see that the inner action $Int: G_{L|F}\rightarrow GL_{\F_\ell}(C)\cong \F_\ell^\times$ will be injective. It can be interpreted as the map $\sigma^i \mapsto q_F^i\;(mod\;\ell),$ where $\sigma\in G_{L|F}$ denotes the Frobenius.
 \item $E|F$ is totally wildly ramified of prime power degree $[E:F]= p^s.$ Then $L|F$ is the maximal tamely ramified subextension in $K|F$ and $C=G_{K|L}\twoheadleftarrow U_L/(U_L)^p$ is a simple $\F_p[G_{L|F}]$-module which is realized as quotient of $U_L/(U_L)^p.$ And the inner action $Int: G_{L|F}\rightarrow GL_{\F_p}(C)$ must be faithful. More precisely with respect to the lower numeration of ramification subgroups we have:
 $$ C= C_i = G_i\quad\textrm{for}\; i=1,...,t,\qquad C_{t+1}= G_{t+1} =\{1\},\qquad C=C_t=C^t\twoheadleftarrow U_L^t/U_L^{t+1},$$
 where $U_L^t \not\subseteq (U_L)^p.$ {\bf [L]}p.172 puts $q:=\#C=[E:F]= dim(Ind_{E|F}(\chi\circ N_{E|F})).$
\end{itemize}
{\bf Proof: If $E|F$ is tamely ramified} then it is well known that the Galois closure $K|F$ will be obtained by adjoining sufficiently many roots of unity of order prime to $p$ that means $K|E$ will be an unramified extension. Let $F_0|F$ be the inertia subfield of $K|F.$ From Lemma 1 of {\bf A2} the inclusion $C=G_{K|L}\le G_{K|F_0}$ will follow because $C$ is contained in every non-trivial normal subgroup. Therefore we have $L|F_0.$ {\bf Moreover $K|F_0$ is normal and totally tamely ramified, hence it will be cyclic}. Then of course $C=G(K|F_0)$ because it is the only abelian normal subgroup, hence $L=F_0$ and $[K:L]=[E:F]=\ell$ must be a prime because $C$ is cyclic and minimal. Therefore the Galois closure of $E|F$ is obtained by adjoining the $\ell$-th roots of unity that means $[L:F]$ must be the order of $q_F\in (Z/\ell)^\times$, which implies $q_L\equiv 1\;(mod\;\ell).$\\
{\bf If $E|F$ is totally wildly ramified} then let $F_1|F$ be the maximal tamely ramified subextension in $K|F$.  Again $C=G_{K|L}\le G_{K|F_1}$ will follow because $C$ is contained in every non-trivial normal subgroup. Thus $C$ and $G_{K|F_1}$ are both normal $p$-subgroups in $G=G_{K|F}.$ But from Lemma 1 of {\bf A2} we see that $C$ is the only normal $p$-subgroup in $G$ and therefore $L=F_1.$ Thus we have: $C=C_1 =G_1$ and this implies $G_i\le C$ for $i\ge 1$ and therefore:
$$   G_i= G_i\cap C =C_i\quad\forall\; i\ge 1,\qquad C=G_{K|L}\twoheadleftarrow U_L^t/U_L^{t+1}.$$
The second implication follows from $G_i=C$ if $G_i\ne \{1\},$ because $C$ is minimal normal. And $C=C_t$ implies $C=C^t$ because lower and upper numeration for $C=G_{K|L}$ agree up to the first jump, and moreover $C_{t+1}=\{1\}$ such that by class field theory $C$ must be a quotient of $U_L^t/U_L^{t+1}.$\\
Again by class field theory this must be a morphism of $\F_p[G_{L|F}]$-modules where the action of $G_{L|F}$ on $C$ is via conjugation. Finally $C$ must be simple because a proper $\F_p[G_{L|F}]$-submodule $C_0$ is nothing else than an abelian normal subgroup of $G=G_{K|F}.$ And from Lemma 1 of {\bf A2} we also see that the inner action must be faithful.\\

In the {\bf tame case} this yields quite a simplification of {\bf (3)}. Indeed $S(K|L)$ consists of the trivial character $\mu_0$ and $(\ell -1)$ characters $\mu\ne\mu_0$ which are faithful because $L^\times/N_{K|L}(K^\times) \cong G_{K|L}=C$ is cyclic of prime order $\ell.$ And $H=G_{K|E}\cong G_{L|F}\hookrightarrow \F_\ell^\times$ is cyclic too and acts (via conjugation) faithful on $C$ and on the dual $C^*$ resp. Thus:\\
$s\in H$, $\mu\ne \mu_0 \in C^*$ and $\mu^s =\mu$,  implies $(\mu^\nu)^s=\mu^\nu$ for all powers of $\mu$, and this implies $s=1.$ This means we have:
$H_\mu =\{1\},$ $E_\mu =K$ and $F_\mu =L,$ hence $\mu' =\mu$, if $\mu\ne\mu_0.$\\
Now taking into account the last remark of {\bf (3)} we obtain the {\bf arithmetic condition:} (see {\bf [L]}, Lemma 13.3 plus remark at the end of p.163)\\
{\bf (3 tame)}
$$  \Delta(E,\chi\circ N_{E|F})\cdot\prod_{[\mu]\ne\{\mu_0\}}\Delta(L,\mu)= \Delta(F,\chi)\cdot \prod_{[\mu]\ne\{\mu_0\}}\Delta(L,(\chi\circ N_{L|F})\mu)$$
should hold for any character $\chi$ of $F^\times,$ if $E|F$ is ramified of prime degree $\ell \ne p,$ and $L|F$ is the inertia subfield of its Galois closure $K|F.$ The products are over representatives for the $\frac{\ell -1}{[L:F]}$ different $H\cong G_{L|F}$-orbits $[\mu]$ of non-trivial characters $\mu\in S(K|L)$ where each orbit is of length $\#H =[L:F].$ And the product on the left vanishes if $[K:F]=\ell\cdot [L:F]$ is odd.\\

.\\

{\bf Appendix 1. Proving Brauer 3 and Brauer 4}\\

First of all we note that both theorems can be reduced to the case where $N=C$ is a commutative normal subgroup that means $[N,N]=\{1\}$. For this replace $H,$ $\OO$ by  $\ol H= H/_{[N,N]}$ and $\ol\OO =\OO/_{[N,N]}$ and observe that $R_+(N\le\OO) =R_+(\ol N\le\ol\OO)$ which brings us down to consider $\vf_{\ol\OO}:R_+(\ol N\le\ol\OO)\rightarrow R(\ol\OO)$  where $\ol N\le\ol\OO$ is an abelian normal subgroup.\\
Now we are going to give proofs of Brauer 3 and Brauer 4 {\it in the case where $N=C$ is a commutative normal subgroup}, hence $\vf_{C\le\OO}:R_+(C\le\OO)\rightarrow R(\OO).$\\
{\bf Brauer 3:} The proof is based on Corollary {\bf 2.3} constructing a {\bf projector}
$$ \Phi_C: R_+(\le\OO) \twoheadrightarrow R_+(C\le\OO) \;\subset\; R_+(\le\OO),\quad\textrm{such that}\quad \vf_\OO=\vf_{C\le\OO}\circ\Phi_C.$$
Then using Brauer 1 for $\rho\in R(\OO)$ we obtain
$$  \rho= \vf_\OO(\sum_i n_i[H_i,\chi_i]) = \vf_{C\le \OO}\left(\sum_i n_i\Phi_C([H_i,\chi_i])\right).$$
Since all our subgroups are of finite index and all characters are of finite order, the construction of $\Phi_C$ does not depend on taking $\OO$ finite or profinite.

{\bf Brauer 4:} First let $\OO$ be finite, where we can use results of {\bf [B]}. We consider $\rho$ an irreducible representation of $\OO$ and prove the assertion for the dimension-0-representation $\rho - dim(\rho)\cdot 1_\OO.$ The restriction of $\rho$ to $C$ yields
$$  \rho|_C = \frac{dim(\rho)}{(\OO:S)} \bigoplus_{g\in S\backslash \OO} \chi_C^g,$$
where $\chi_C$ is a certain character of $C$ and $S:=\{g\in \OO\;|\; \chi_C^g =\chi_C\}$ is its stabilizer. Moreover we have then:
$$  \rho\cong Ind_S^\OO(\rho_S),\qquad\textrm{where}\quad \rho_S|_C = \frac{dim(\rho)}{(\OO:S)}\cdot \chi_C$$
that means the representation $\rho_S$ is central on $C$.  Now for $\rho_S\in R(S)$ we consider Boltje's canonical Brauer presentation {\bf [B],(2.6)}
$$   \rho_S \cong \sum_{[H,\chi]} \alpha_{[H,\chi]}(\rho_S) Ind_H^S(\chi),$$
with integral coefficients $\alpha_{[H,\chi]}(\rho_S)$,
where according to {\bf loc.cit. item k)(iii) on p.40}, only subgroups $(H,\chi)\ge (C,\chi_C)$ (that means $\chi|_C=\chi_C$) will occur because $\rho_S$ is central on $C.$ Moreover this canonical presentation has the property {\bf loc.cit. (2.23)}:
$$  dim(\rho_S)\cdot 1_S \cong \sum_{[H,\chi]} \alpha_{[H,\chi]}(\rho_S) Ind_H^S(1_H),$$
such that subtracting this from the first equivalence we have:
$$   \rho_S - dim(\rho_S)\cdot 1_S \cong \sum_{[H,\chi]} \alpha_{[H,\chi]}(\rho_S) Ind_H^S(\chi - 1_H).$$
Thus we are done in the particular case where $S=\OO.$ Otherwise
we induce this equivalence from $S$ to $\OO$ and obtain:
$$  \rho - dim(\rho)\cdot 1_\OO +\left(dim(\rho)\cdot 1_\OO - dim(\rho_S) Ind_S^\OO(1_S)\right)\cong \sum_{[H,\chi]} \alpha_{[H,\chi]}(\rho_S) Ind_H^\OO(\chi - 1_H).$$
But $Ind_S^\OO(1_S)= Ind_{S/C}^{\OO/C}(1_{S/C})$ is a representation of $\OO/C$ and therefore\\ $dim(\rho)\cdot 1_\OO - dim(\rho_S) Ind_S^\OO(1_S)$ is a virtual representation of $\OO/C$ which is of dimension $0$. Here we may apply {\bf Brauer 2} with respect to the finite factor group $\OO/C.$ (A subgroup of $\OO/C$ has always the form $H/C$ for $H\ge C$). Then substituting the corresponding expression for $dim(\rho)\cdot 1_\OO - dim(\rho_S) Ind_S^\OO(1_S)$, we get our assertion for $\rho - dim(\rho)\cdot 1,$ if $\rho$ is irreducible.\\
If now $\rho$ is any virtual representation of $\OO$, of dimension $0$, then using that the irreducible representations form a $\Z$-basis of the Grothendieck group $R(\OO)$ we have:
$$ \rho = \sum_i n_i\rho_i,\qquad 0 = \sum_i n_i dim(\rho_i).$$
Therefore $\rho = \sum_i n_i(\rho_i - dim(\rho_i)\cdot 1_\OO),$ and to each term we may apply our previous result. {\bf qed.}  If $\OO=\varprojlim \OO/U$ is profinite and $N\le\OO$ is an open normal subgroup then 
$$  \vf_{N/U\cap N\le\OO/U}: R_+(N/U\cap N\le\OO/U)\twoheadrightarrow R(\OO/_{[N,N]U})$$
brings us back to the case of finite groups, and the result follows from $\vf_{N\le\OO}=\varinjlim_U \vf_{N/U\cap N \le \OO/U}.$\\

{\bf Appendix 2 on type-III-groups}\\

The groups to be considered here are special examples of finite solvable groups $G$ with trivial Frattini subgroup:
$\Phi(G)=\textrm{intersection of all maximal subgroups of $G$}=\{1\}.$ For readers with some experience in group theory these remarks are more or less obvious.\\
By definition a finite group $G$ will be called of type III if
\begin{itemize}
 \item  $G$ is solvable
 \item $G$ has a proper subgroup $H< G$ such that
 \begin{itemize}
  \item $H$ is a maximal subgroup,
  \item $\bigcap_{g\in G} gHg^{-1} =\{1\}.$
 \end{itemize}
\end{itemize}
(For local Galois groups $G=G_{K|F}$ this means that $K|F$ is the Galois closure of a minmal extension $E|F$ corresponding to $H=G_{K|E}< G$; see {\bf [L] Lemma 3.2}).\\
{\bf Remark:} Note that $H=\{1\}$ is a maximal subgroup if and only if $G$ is cyclic of prime order. This is the degenerate case which we will not exclude.\\
Let $G$ be a type-III-group. Since $G$ is solvable there must be a non-trivial commutative normal subgroup $C\le G.$\\

{\bf Lemma 1:} {\it Let $G$ be a type-III-group and $H<G$ a corresponding maximal subgroup. Then any nontrivial commutative normal subgroup $C\le G$ must be a complement of $H$ in $G$, and actually there exists only one such $C$. In particular, $C$ will be contained in every non-trivial normal subgroup of $G.$ (The degenerate case is $H=\{1\}$ and $C=G.$) Therefore $C$ must be $\ell$-elementary:
$$ C\cong \Z/\ell \times \cdots \times \Z/\ell$$
for a certain prime $\ell =\ell(G)$ such that $(G:H)=\# C$ is a power of $\ell.$  At the same time $C$ is the unique {\bf normal} $\ell$-subgroup of $G$, hence the $\ell$-Sylow-subgroups of $G$ and of $H$ must have intersection $C$ and $\{1\}$ resp. And $C$ is its own centralizer that means $h\in H \mapsto Int(h)\in GL_{F_\ell}(C)$ is an embedding by inner automorphisms which induces $G=H\ltimes C\hookrightarrow GL_{\F_\ell}(C)\ltimes C.$}\\
{\bf $C$ can be also characterized as}
$C =F(G)$ the Fitting subgroup (=maximal nilpotent normal subgroup) of $G.$ This follows because the Fitting subgroup will be abelian if the Frattini subgroup is trivial. And it shows that in the non-degenerate case $G$ is not nilpotent.\\

{\bf Proof:} Let $C$ be a non-trivial commutative normal subgroup. We will have $H\cdot C=G$  because $H$ is maximal and cannot contain the normal subgroup $C$. But $H\cap C$ is normal in $HC$, hence $H\cap C=\{1\}$ and $G=H\ltimes C$ is semidirect.\\
Thus we have seen that all nontrivial commutative normal subgroups must have the same order $\# C=(G:H),$ and this implies that $C$ must be unique. Then it is also clear that $C$ is the absolutely minimal {\bf normal} subgroup in the solvable group $G,$
hence $C\cong \Z/\ell \times \cdots \times \Z/\ell$ for a certain prime $\ell=\ell(G).$\\
Let $B$ be an $\ell$-subgroup of $G$ which is normal. Then the center $Z(B)\ne \{1\}$ is a normal abelian subgroup of $G,$ hence $Z(B)=C$ and $B=(B\cap H)\ltimes C$ is actually a direct product because $C$ is the center of $B.$ Therefore $B\cap H$ is not only normal in $H$ but also in $G=H\ltimes C.$ However by the very definition $H$ cannot contain normal subgroups of $G$, and therefore $B\cap H=\{1\}$ and $B=C.$ Obviously $C$ is then the intersection of all $\ell$-Sylow groups $P$ of $G.$ On the other hand the $\ell$-Sylow groups $P_H$ come as $P_H = H\cap P$ for all $P$ and therefore the intersection of all $P_H$ is $H\cap C =\{1\}.$\\
Similarly the centralizer $\wt C$ of $C$ has the form $\wt C =(\wt C \cap H)\times C$ where $\wt C\cap H$ is normal in $G,$ hence $\wt C\cap H =\{1\}.$
{\bf qed.}\\

{\bf Lemma 2:} {\it Let $G> H$ be a type-III-group and assume that it is also type-III with respect to $G> H'.$ Then we will have $H' = c_0H c_0^{-1}$ for some $c_0\in C,$ and in the non-degenerate case $H\ne\{1\}$ the conjugates $\{cHc^{-1}\;|\; c\in C\}$ will be all different.}\\

{\bf Proof:} We assume $H\ne\{1\}.$ By assumption all statements of Lemma 1 also hold for $H'$ instead of $H,$ in particular $(G:H')=(G:H)=\#C.$ We consider the inclusion
$$ H\subset G= C\rtimes H',\qquad h = c(h)\cdot h'(h).$$
Replacing $h$ by $h_1h_2$ we conclude:
$$  c(h_1h_2) = c(h_1)\cdot h'(h_1)c(h_2)h'(h_1)^{-1} = c(h_1)\cdot h_1c(h_2)h_1^{-1},\qquad h'(h_1h_2)=h'(h_1)\cdot h'(h_2),$$
because $C$ is an abelian normal subgroup. Thus $\{h\mapsto c(h)\} \in Z^1(H,C)$ is a 1-cocycle for the $H$-module $C$ (where the action is defined by conjugation). We are going to see that $H^1(H,C)=\{1\},$ hence $c(h)= hc_0h^{-1}\cdot c_0^{-1}$ is a coboundary for some fixed $c_0\in C$, and therefore:
$$   h = c(h)\cdot h'(h) = (hc_0h^{-1}\cdot c_0^{-1})\cdot h'(h),\qquad c_0hc_0^{-1} =h'(h),\quad \forall \;h\in H,$$
{\bf Sublemma:} {\it If $G=C\rtimes H$ is a proper type-III-group, then $H^1(H,C)=\{1\}.$}\\
{\bf Proof:} Since $H$ is also solvable, there exists a minimal normal subgroup $H_0\le H$ of type
$H_0 \cong \Z/\ell_0\times \cdots \Z/\ell_0,$ for some prime $\ell_0$. And $\ell_0\ne \ell$ because the $\ell$-Sylow subgroups of $H$ have intersection $\{1\}.$  Furthermore $H_0$ normal in $H$ implies $C^{H_0}$(=invariants under conjugation) normal in $G,$ hence the only possibilities are $C^{H_0}=\{1\}$ or $C^{H_0}=C,$ where the second case cannot happen.  Thus we must have $C^{H_0}=\{1\}.$ Finally we look at the inflation-restriction sequence:
$$  1\rightarrow H^1(H/H_0, C^{H_0}) \rightarrow H^1(H,C) \rightarrow H^1(H_0,C).$$
Then $H^1(H_0,C)=\{1\}$ because the order of $H_0$ is prime to $\ell,$ and therefore together with $C^{H_0}=\{1\}$ we obtain $H^1(H,C)=\{1\}.$  {\bf qed.}\\
Finally assume $cHc^{-1} =H.$ Then we have $chc^{-1}h^{-1}\in H\cap C=\{1\}$ for all $h$, hence $c$ will be in the centralizer $Z_C(H).$ However $Z_C(H) < C$ will be a normal subgroup of $G=H\cdot C$ and therefore $Z_C(H)=\{1\},$ hence $cHc^{-1}=H$ implies $c=1.\quad$ {\bf qed.}\\

And as a converse we note:\\
{\bf Lemma 3:} {\it Let $G$ be a solvable group which comes as a semidirect product $G=H\ltimes C.$ Moreover assume that $C\ne\{1\}$ is minimal normal in $G$ (hence commutative) and let $K:=\bigcap_{g\in G} gHg^{-1}.$ Then:\\
{\bf (i)} $H< G$ is a maximal subgroup.\\
{\bf (ii)} $K= H\cap\wt C$ is the intersection of $H$ with the centralizer of $C.$\\
Therefore $H<G$ will be a type-III-group if and only if $H\cap\wt C=\{1\}$ or equivalently: $\wt C =C$ or equivalently: $H$ does not contain normal subgroups of $G$.}\\
{\bf Proof:} {\bf (i)}:  Assume $H< U< G.$ Then $G=HC$ implies $G=UC$ and $U\cap C$ is a normal subgroup of $U,$ hence also of $UC=G.$   On the other hand: $U/U\cap C \cong UC/C \cong H$ and therefore $\#(U\cap C)=(U:H)$ is a proper divisor of $(G:H)=\#C.$ Therefore $C$ minimal implies $U\cap C=\{1\}$, hence $U=H$ which is a contradiction.\\
{\bf (ii)}: $K$ is by definition a normal subgroup of $G$ which sits in $H.$ Therefore $[K,C]\le K\cap C\le H\cap C=\{1\},$ hence $K\le H\cap\wt C.$ On the other hand $H\cap\wt C$ is clearly normal in $G=HC,$ and therefore $H\cap\wt C\le K.$  {\bf qed.}\\

\centerline {\bf \large References}

[B] R. Boltje, A canonical Brauer induction formula, Asterisque, tome 181-182(1990),p.31-59.\\

[BH] C. Bushnell and G. Henniart, The Local Langlands Conjecture for GL(2), Grundlehren der mathematischen Wissenschaften, Vol.335 (Springer-Verlag, Berlin Heidelberg 2006).\\

[D] P. Deligne, Les constantes des equations fonctionelle de fonction L, Springer Lecture Notes {\bf 349} (1974) pp. 501-597.\\

[K] H. Koch, Extendible Functions, C.I.C.M.A.-preprint, Concordia University of Edmonton, Dept. of Math., September 1990\\

[L] R.P. Langlands, On the functional equation of the Artin L-functions (1970),\\   https://publications.ias.edu/sites/default/files/a-ps.pdf\\

[M] J. Martinet, Character theory and Artin L-functions, in: Algebraic Number Fields, edited by A. Fr\"ohlich, Academic Press London e.a. 1977, pp. 1-87.\\

[S] J.-P. Serre, Linear Representations of Finite Groups (Translation of the French edition), New York 1977 (Corrected fourth printing, 1993)\\

[T1] J. Tate, Fourier analysis in number fields and Hecke's zeta-functions (Thesis, Princeton 1950).{\it Algebraic Number Theory (Proc. Instructional Conf. Brighton, 1965)} (Thompson, Washington, D.C.), pp. 305-347.\\

[T2] J. Tate, Local constants, in Algebraic Number Fields,  edited by A. Fr\"ohlich, Academic Press London e.a. 1977, pp. 89-131.\\

[T3] J. Tate, Number Theoretic Background, in Automorphic Forms, Representations, and L-Functions, PSPM Vol.33(1979), part 2, pp.3-26,  AMS 1979\\

[W] A. Weil, Basic Number Theory (Grundlehren der mathematischen Wissenschaften, Band 144), Berlin Heidelberg New York 1967\\

.\\

Helmut Koch\\
Oybiner Str. 12\\
Dresden D-01324, Germany\\

Ernst-Wilhelm Zink\\
Institut f\"ur Mathematik\\
Humboldt-Universit\"at zu Berlin\\
Unter den Linden 6\\
Berlin D-10099, Germany\\
zink@mathematik.hu-berlin.de

\end{document}